\documentclass[10pt]{article}
\usepackage{amssymb, amsfonts}
\usepackage[english]{babel}
\usepackage{hyperref}
\usepackage{stmaryrd}
\usepackage{mathtools}
\input{liemacs10.sty}

\DeclareMathAlphabet{\Ma}{U}{msa}{m}{n}
\DeclareMathAlphabet{\Mb}{U}{msb}{m}{n}
\DeclareMathAlphabet{\Meuf}{U}{euf}{m}{n}

\def\got#1{\Meuf{#1}}
\def\br#1.{\llbracket #1 \rrbracket}
\def\lbr{\llbracket}
\def\rbr{\rrbracket}

\newcommand{\bc}{\mathbf{c}}
\newcommand{\bd}{\mathbf{d}}
\newcommand{\cJ}{\mathcal{J}}
\newcommand{\PSU}{\mathop{{\rm PSU}}\nolimits}

\newcommand{\ran}{\mathop{{\rm ran}}\nolimits}
\renewcommand{\phi}{\varphi}

\renewcommand{\baselinestretch}{1.2}
\normalsize
 \def\ot #1.{{\got{#1}}}

\newcommand{\Osc}{\mathop{{\rm Osc}}\nolimits}
\renewcommand{\mlabel}{\label}

\def\cR{\mathcal{R}}
\def\al#1.{{\cal #1}}
\def\wt{\widetilde}
\def\slim{\mathop{\hbox{\rm s-lim}}}
\def\s #1.{_{\smash{\lower2pt\hbox{\mathsurround=0pt $\scriptstyle #1$}}\mathsurround=5pt}}
\def\XP#1!{\renewcommand{\baselinestretch}{.7}\marginpar{{\footnotesize
$\leftarrow$#1}\hfil}\renewcommand{\baselinestretch}{1.2}}
\def\ccr #1,#2.{\overline{\Delta(#1,\,#2)}}

\def\b#1.{{\bf #1}}

\def\hlf{\frac{1}{2}}
\def\un{\1}
\def\bn{{\bf n}}

\begin{document}


\title{Crossed products of $C^*$-algebras for singular actions with spectrum conditions}
\author{Hendrik Grundling, Karl-Hermann Neeb}

\maketitle

\begin{abstract}
We  analyze existence of crossed product constructions of Lie group actions on \break
$C^*$-algebras which are singular. 
These are actions where the group need not be locally compact, or the action need not be strongly continuous.
In particular, we consider the case where spectrum conditions are required for the implementing
unitary group in covariant representations of such actions.
The existence of a crossed product construction is guaranteed by the
existence of ``cross representations''. For one-parameter automorphism groups, we prove
that the existence of cross representations is stable with respect to a large set of perturbations of the action,
and we fully analyze the structure of cross representations of inner actions on von Neumann algebras.
For one-parameter automorphism groups we study the cross property for covariant representations, where
the generator of the implementing unitary group is positive. In particular, we find that if the
Borchers--Arveson minimal implementing group is cross, then so are all other implementing groups.
We study a smoothing phenomenon for one-parameter actions on Lie groups, and display the usefulness of
cross representations for this context.
For higher dimensional Lie group actions,  we consider a class of spectral conditions which include
the ones occurring in physics, and is sensible also for non-abelian or for infinite dimensional Lie groups.
We prove that the cross property of a covariant representation is fully determined by the cross property
of a certain one-parameter subsystem. This greatly simplifies the analysis of the existence of cross
representations, and it allows us to prove the cross property for several examples of interest to
physics. We also consider non-abelian extensions of the Borchers--Arveson theorem. There is a full
extension in the presence of a cyclic invariant vector, but otherwise one needs to
determine the vanishing of lifting obstructions.
\end{abstract}

\tableofcontents

\section{Introduction}

In this paper, we interweave two  strands of enquiry, which are singular actions
${\alpha \: G \to \Aut(\cA)}$ of topological groups $G$ on $C^*$-algebras $\cA,$
and spectral conditions for covariant representations.
A {\it singular action} $\alpha$ is an action for which
either $G$ is not locally compact or the  $G$-action on $\cA$
is not strongly continuous. For these actions the usual construction of a crossed
product \break 
$C^*$-algebra $\cA \rtimes_\alpha G$ fails.
Singular actions are abundant in physics and arise naturally in
mathematics in important examples.
For instance, the field $C^*$-algebra for bosonic field theories  is usually
chosen to be either the Weyl algebra or the resolvent algebra (\cite{BG08}), and then
non-constant one-parameter symplectic groups
produce one-parameter automorphism groups on these algebras which are not
strongly continuous (cf.\ \cite[Exs.~2.8/9]{BGN17}).
Singular actions of the non-locally compact group
$C^\infty(M,K)$ of local gauge transformations arise also in
quantum gauge theories. In fact, Borchers pointed out
on many occasions (\cite{Bo83, Bo87, Bo96}), that in quantum field theory
the natural actions $\alpha \: G \to \Aut(\cA)$
are usually not strongly continuous, and the groups $G$
need not be locally compact. E.g., the group $G$
may be an infinite dimensional Lie group,
such as the group $C^\infty(M,K)$ of gauge transformations,
or the group $\Diff(M)$ of diffeomorphisms of a compact manifold.
It may also be a topological subgroup of the unitary group in
a von Neumann algebra which carries no manifold structure.
We refer to \cite{BGN17} for a review of covariant
representations for singular actions.

For singular actions, many of the usual mathematical tools break down.
In particular, the non-existence of a $C^*$-crossed product shows that
there is no good global
structure theory for their covariant representations.
We addressed this problem in a previous work (cf.~\cite{GrN14}).
The difficulty arising of $G$ not being locally compact
can be overcome by imposing
regularity conditions on the considered unitary representations
$U \: G \to \U(\cH)$. A~particularly natural requirement is that
the class of specified representations of $G$ is in one-to-one
correspondence with the non-degenerate representations of a
$C^*$-algebra~$\cL$, called a {\it host algebra of~$G$}
(see Definition~\ref{def:2.1a}, or \cite{Gr05}).
If $G$ is locally compact, then $\cL$ will typically be a quotient of
$C^*(G)$, and if $G$ is an infinite dimensional Lie group, we shall
see other natural examples of host algebras.

In \cite{GrN14}, we introduced a crossed product construction which is possible
for a subclass of singular actions (cf.~\cite{GrN14}), relative to a choice of
host algebra $\cL$, which controls the unitary representations allowed for covariant representations.
This is particularly useful in that the existence of such a
``crossed product host'' brings with it
 a good structure theory for states and a
subclass of representations associated with it, and offers
 tools such as direct integral decompositions.
Depending on the choice of the class of representations required,
controlled by a host algebra $\cL$, these crossed product hosts
are unique, but the existence problem is much harder and depends on the
existence of  {\it cross representations}  (cf.~Definition~\ref{crossrepDef}).
Whilst in \cite{GrN14} we provided many examples of  cross representations
where the usual crossed product does not exist, in general they are not easy
to find for an arbitrary singular action, so here we want to extend this analysis.
 In particular, one should be able to
deal with the main covariant representations of physical interest. Our first
problem is to extend the small set of known systems in \cite{GrN14},
and so we analyze the cross property with respect to  perturbations of the action.
This allows us to show that, for the Fock representation of a bosonic system
equipped with the dynamics of a positive one-particle Hamiltonian,
the second quantized covariant representation is cross (Example~\ref{NoGap}).
For the special case of a singular action
 $\alpha \: G \to \Aut(\cA)$, where $\cA$ is a von Neumann algebra
 and $\alpha$ is inner, we  fully characterize
 the cross representations in Subsection~\ref{InnerW}. These will restrict to
 cross representations on invariant $C^*$-subalgebras, but such subalgebras may have
 cross representations which do not extend to cross representations of the von Neumann algebra.

The second strand in this paper concerns  covariant representations
of singular actions satisfying spectral conditions. Here $G$ is a Lie group
and for some elements $X$ in its  Lie algebra, we require
that the selfadjoint
generator $H_X$ of the corresponding unitary one-parameter group
$U(\exp tX) = e^{it H_X}$ is bounded from below.
As the concrete lower spectral bounds are
specified in terms of a closed convex subset $C$ in the dual
of the Lie algebra of $G$, we speak of a {\it $C$-spectral condition}
(Definition~\ref{defSpec}). Of particular importance are
representations which are {\it semibounded}.
This is a stable version of the positive spectrum condition
which means that the operators $H_X$ are uniformly bounded from below for
all $X$ in some non-empty open subset of the Lie algebra of~$G$.

Special tools are available for such representations,
which allow us to further analyze the cross property for these.
Covariant representations of this type are of profound interest in quantum physics, as the
generator of time evolution is almost always required to be a positive
selfadjoint operator.
This has a natural generalization to higher dimensional abelian groups,
e.g.\ in relativistic quantum field
theory, and non-abelian extensions have also been studied
(\cite{Ne99, Ne10, Ne12, Ne14b,  Ne17}).
Lie group representations
satisfying spectral conditions are of fundamental importance in
physics~\cite{SW64, Bo87, Bo96, Ot95, H92, LM75},
harmonic analysis~\cite{Ol82, Ol90, HO96, Ne99, Ne10, Ne12, Ne14b,  Ne17}
and operator algebras (cf.~\cite[Ch.~8]{Pe89})
and have some nice structural properties
(cf. \cite{Ne99, Ne10}). The paper
\cite{Bo84} seems to be the first one where the spectrum
condition is studied for group actions which are not strongly
continuous, which is in our focus here.
We already considered the cross property for a few of these
 systems (cf.\ Examples~6.11 and 9.1--9.3 in  \cite{GrN14}),
but here we want to pursue the general analysis of this question,
because new mathematical tools such as
smoothing operators for infinite dimensional Lie groups
have recently become available (\cite{NSZ17}).

The setting of a discontinuous action of $G = \R$ is the easiest
to study the cross property
for positive covariant representations of singular actions,
and so we start this strand with that.
Here the generator of the implementing unitary group is positive.
This will turn out to
be important, as we will show that for $C$-spectral covariant representations
on higher dimensional Lie groups,
the cross property for the full system can be tested
on  certain one-parameter subsystems.
The Borchers--Arveson Theorem provides a special ``minimal'' inner implementing
unitary group for positive covariant representations, and we show that if the
 Borchers--Arveson inner covariant representation is cross, then so are all the
 other positive covariant representations, keeping the representation
of the algebra itself fixed.

As the generalization of the positive spectrum condition
to non-abelian Lie groups is unfamiliar to physicists, we present it in detail,
including some standard tools such as complex Olshanski semigroups,
generalizing the complex upper half plane.
We explore interesting properties of covariant representations satisfying $C$-spectral
 conditions, and we obtain an important  result
that reduces the question about the cross condition to that of the representation restricted
to some one-parameter subgroup (Theorem~\ref{CrossSpecParam2}).
This is extremely useful, and allows us to show that
the Fock representation for a bosonic system
is cross for the symplectic action of the conformal
group $\SO(2,d)$ (which contains the Poincare group) on the Weyl algebra
where the action is defined by a semibounded representation of $\SO(2,d)$.
We also obtain a generalization to the non-abelian case of the Borchers--Arveson Theorem
for when there is an invariant cyclic vector, and
in the general case there are lifting obstructions that we describe in terms of
central extensions.

In more detail, the  paper is structured as follows.
After establishing notation, we give a brief review of the results on crossed product hosts
from \cite{GrN14} which we will need. A few new results (e.g. on the Laplacian) are added with proof.
In Subsection~\ref{CrossPertb} we analyze the question on whether cross representations are stable
under perturbation and how their associated crossed product hosts are related. We obtain conditions
that allow us to establish the cross property for the important physics example
consisting of the Fock representation of the Weyl algebra
equipped with the dynamics of a positive one-particle Hamiltonian (Example~\ref{NoGap}).
We move to the special case of an inner singular action on a $W^*$-algebra in
Subsection~\ref{InnerW},
analyze its normal cross representations   and fully characterize these
in Theorem~\ref{W-CrossRepsChar}.
Our assumption of an inner action is not restrictive if there are (positive) spectrum conditions,
as the Borchers--Arveson Theorem ensures that such actions
must be inner.

We then turn to spectral conditions for covariant representations and start
with  the one-parameter case, where we can use the  generator
of the implementing unitary group  to concretely analyze the cross
property (Section~\ref{RepPosCros}).
One finds that if the Borchers--Arveson minimal inner covariant representation is cross,
then so is the original one
 (Theorem~\ref{CrossRepsSame}).
In the following we describe how to factor out an ideal from a crossed product host to obtain one
which only allows positive covariant representations (Proposition~\ref{Cideal}),
and  how to build a positive covariant cross representation from a covariant cross
representation which does not satisfy the positive spectral condition (Proposition~\ref{CrossConstrain}).

The analysis of positive covariant representations for one-parameter automorphism
groups of groups (rather than algebras) has been well-studied in mathematics
(cf.~\cite{Ne14, NS14}), and  we consider how
the current analysis connects with that area in Section~\ref{CRATG}.
Here one would choose for the algebra with the
singular action the discrete group algebra $C^*(G_d)$. This leads to an interesting phenomenon,
where in a positive cross representation $(\pi, U)$ for the action
${\alpha \: \R \to \Aut(G)}$,
 the positive unitary one-parameter group $(U_t)_{t \in \R}$
can ``regularize'' the representation $\pi$ of
the group~$G$. It defines a crossed product host which only allows
continuous representations of~$G$.
A~typical example of this is the Fock representation of
the Weyl algebra,
where the one-parameter group generated by the number operator produces a crossed product host which only
allows regular cross representations. This leads naturally to the study of smoothing operators
(cf.~\cite{NSZ17}) where this phenomenon is studied in greater depth, and in
Subsection~\ref{subsec:smop}
we recall some results from \cite{NSZ17} and  establish the link with
crossed product hosts (Theorem~\ref{thm:5.6}).

We next generalize the positive spectral condition from $\R$ to general Lie groups (possibly
non-abelian or infinite dimensional) to define  ``$C$-spectral representations''
(Section~\ref{SpecCondRev}).
For brevity, we  only study finite dimensional Lie groups, possibly non-abelian.
As  host algebras for  $C$-spectral representations
are not well-known in physics, we give two constructions,  one via smoothing operators
in Subsection~\ref{smoothHost}, and another one via the Olshanski semigroup,
defined in Subsection~\ref{OlshHost}. There are some interesting
restriction properties for these, listed in
Proposition~\ref{ResCspec} and Corollary~\ref{subgpalgcont},
not shared by the usual group algebra $C^*(G)$.
These host algebras are now used to analyze the cross property
for covariant $C$-spectral representations (Subsection~\ref{CrossParam}).
This culminates in the result that  the property of whether a $C$-spectral representation is cross
 is fully
determined by whether its restriction to some one-parameter subgroup
is cross or not (Theorem~\ref{CrossSpecParam2}).
Thus the previous analysis of the cross property for the one-parameter case applies, and this is so
useful that we can immediately  establish the cross property for
the Fock representation on the Weyl algebra, for the symplectic actions of either the
translation group on Minkowski
space or  the conformal group (which contains the Poincar\'e group),
 where the spectral condition is the
usual one from physics, i.e.\ the
joint spectrum of the generators of translation is contained in the closed forward
light cone (Example~\ref{FockPosCross}).
Unfortunately we do not know whether the restriction of this representation to
the Poincar\'e group is also cross.
In the case that a conventional crossed product exists, we show
that our crossed product hosts are
quotients of the conventional one.

Given the generalization of the spectral condition to non-abelian groups, we need to check how
to adapt tools from the abelian case. One natural question is whether one can
obtain a version of the Borchers--Arveson Theorem, and indeed we find
 that if we have a $C$-spectral covariant representation containing a cyclic
$G$-invariant vector,
then there is an inner $C$-spectral representation such that on the one-parameter
subgroups corresponding to semibounded generators,
it restricts to the Borchers--Arveson minimal groups
(Proposition~\ref{BA-thmGen}). If there is no
cyclic invariant vector, the general problem is
further explored in Subsection~\ref{BATnonab}.

Finally, we conclude with a section containing further examples, first, in
 Subsection~\ref{subsec:8.1} a full crossed product host with a non-abelian spectrum condition.
It is constructed for the group of symplectic transformations acting on the Weyl algebra
of a finite dimensional symplectic space.
For the next example we consider the following.
The translation action, whilst strongly continuous on $C_0(\R)$ (hence cross
in all covariant representations),
is singular on the larger $C^*$-algebra
$C_b(\R)$ and so it is a natural question as to whether the
representation on $L^2(\R)$ is cross. In Subsection~\ref{transCb}
we show that this is not the case. In Subsection~\ref{CPHinfdim} we finally
give two examples of  hosts for infinite dimensional Lie groups, one for the
Virasoro algebra, and one for the twisted loop group. In both cases one proves that
for a certain one-parameter group acting on these, a positive covariant representation
is cross using smoothing operators, hence the crossed product host constructed from
the discrete group algebra is a host algebra for the respective central extensions of the groups.

\section{Basic concepts and notation}

Below we will need the following.
\begin{defn}
\mlabel{def:1.1c}
\begin{itemize}
\item[(i)]
For a $C^*$-algebra $\cA$, we write $M({\cal A})$ for the
multiplier algebra of ${\cal A}$. If
${\cal A}$ has a unit, $\U({\cal A})$ denotes  its unitary group.
There is an injective morphism of $C^*$-algebras
$\iota_{\cal A} \: {\cal A} \to M({\cal A})$ and we will just write
${\cal A}$ for its
image in $M({\cal A})$. Then ${\cal A}$ is dense
in  $M({\cal A})$ with respect to
the {\it strict topology},
which is the locally convex topology defined by the seminorms
$$ p_A(M) := \|M \cdot A\| + \|A \cdot M\|,
\qquad A\in {\cal A},\; M\in M({\cal A})$$
(cf.\ \cite[Prop.~3.5]{Bu68} and \cite[Prop.~2.2]{Wo95}).
\item[(ii)]
Let $\cA$ and $\cL$ be $C^*$-algebras and
$\phi \: \cA \to M(\cL)$ be a morphism of $C^*$-algebras. We
say that $\phi$ is {\it non-degenerate} if ${\rm span}(\phi(\cA)\cL)$
 is dense in $\cL$ (cf.\ \cite{Rae88}).
A representation $\pi:\cA\to{\cal B}({\cal H})$ is called {\it non-degenerate}
if $\pi(\cA){\cal H}$ is dense in the Hilbert space ${\cal H}$.
\end{itemize}

\end{defn}
If $\phi \: \cA \to M(\cB)$ is a morphism of $C^*$-algebras
which is non-degenerate, then we  write $\tilde\phi \: M(\cA) \to M(\cB)$ for its
uniquely determined extension to the multiplier algebras
(cf.\ \cite[Prop.~10.3]{Ne08}).
\medskip

For a complex Hilbert space ${\cal H}$, we write $\Rep({\cal A},{\cal H})$ for the
set of non-degenerate representations of ${\cal A}$ on ${\cal H}$,
and $\ot S.(\cA)$ for the set of states of $\cA.$
To avoid set--theoretic subtleties, we will express our results below
concretely, i.e., in terms of $\Rep({\cal A},{\cal H})$ for given Hilbert spaces
$\al H..$
We have an injection
$$ \Rep({\cal A}, {\cal H}) \into \Rep(M({\cal A}),{\cal H}), \quad \pi \mapsto \tilde\pi
\quad \hbox{ with } \quad \tilde\pi \circ \iota_{\cal A} = \pi, $$
which identifies a non-degenerate representation $\pi$ of
${\cal A}$ with the representation $\tilde\pi$ of its multiplier algebra
which extends $\pi$ on the same Hilbert space.  The representations of $M(\cA)$ on $\cH$ arising
from this extension process are characterized as those representations which are
continuous with respect to the strict topology on $M(\cA)$
and the strong operator topology on $\cB(\cH)$, or equivalently
by non-degeneracy of their restriction to
$\cA$ (cf.\ \cite[Prop.~10.4]{Ne08}).
We will refer to $\tilde\pi$ as the {\it multiplier extension}
of $\pi$.
It can be obtained by
\[
\wt{\pi}(M)=\slim \pi(M E_\lambda)\quad
\mbox{ for } \quad M\in M(\al A.)
\]
where $(E_\lambda)_{\lambda\in\Lambda}$ is any approximate
identity  of $\al A..$
\medskip

\begin{rem} \mlabel{rem:factorize}
Below we will use the notation $\br S.:=\overline{\rm span}(S)$, where $S\subset Y$
 and $Y$ is a Banach space.
If $\cB$ is a $C^*$-algebra, and  $X$ is a left Banach $\cB\hbox{--module,}$
then the closed span of $\cB X$ satisfies $\lbr \cB X\rbr=\cB X=
 \{ Bx \mid B \in \cB, x \in X\}$ (cf.~\cite[Th.~II.5.3.7]{Bla06} or \cite[Th.~5.2.2]{Pa94}).
 In particular it implies that if $\phi \: \cA \to M(\cL)$ is non-degenerate, then
 $\cL=\phi(\cA)\cL$, and if $\pi:\cA\to \cB(\cH)$ is a non-degenerate representation, then
 $\lbr \pi(\cA)\cH\rbr=\pi(\cA)\cH$.
\end{rem}

For a topological group $G$ we will write
$\Rep(G,{\cal H})$ for the set of all (strong operator)
continuous unitary representations of
$G$ on ${\cal H}\,.$ 
Moreover $G_d$ will denote the group $G$ equipped with the discrete topology.

\begin{defn}\mlabel{def:1.1}
(i) We write  $(\cA, G, \alpha)$ for a triple, where
$\cA$ is a $C^*$-algebra, $G$ a topological group and
$\alpha \: G \to \Aut(\cA)$ is a homomorphism.
We call $\alpha$ {\it strongly continuous} if for every $A \in \cA$, the
orbit map $\alpha^A \: G \to \cA,$ $g \mapsto \alpha_g(A)$ is continuous.
If $\alpha$ is strongly continuous,
we call $(\cA, G, \alpha)$ a {\it $C^*$-dynamical system} (cf.~\cite{Pe89},
\cite[Def.~2.7.1]{BR02}), or say that the action is  strongly continuous.
Unless otherwise stated, we will not assume that $\alpha$ has this property
and simply speak of the triple $(\cA, G, \alpha)$ as a
{\it $C^*$-action}.
The {\it usual case} will mean that the action is  strongly continuous, and the group
$G$ is locally compact. A {\it singular} action is one which is not the
usual case.

Given any $C^*$-action $(\cA, G, \alpha)$, we can always define the strongly continuous part of it by
\[
{\cal A}_c:=\{A\in{\cal A}\,\mid\, g\mapsto\alpha_g(A)\quad\hbox{is norm continuous}\}
\quad \mbox{ and } \quad
\alpha^c_g:=\alpha_g\restriction{\cal A}_c.
\]
In the case where $\cA$ is a von Neumann algebra or a $W^*$-algebra and the orbit maps
of elements $A \in  \cA$ are $\sigma(\cA,\cA_*)$-continuous, we  call
$(\cA, G, \alpha)$ a {\it $W^*$-dynamical system}.

(ii) A {\it covariant representation of $(\cA, G, \alpha)$}
is a pair $(\pi,U)$, where
$\pi \: \cA \to  \cB(\cH)$ is a non-degenerate representation of $\cA$
on the Hilbert space $\cH$ and $U \: G \to \U(\cH)$ is a
continuous unitary representation satisfying
\begin{equation}
  \label{eq:covar}
U(g)\pi(A)U(g)^* = \pi(\alpha_g(A)) \quad
\mbox{ for } \quad g \in G, a \in \cA.
\end{equation}
For a fixed Hilbert space $\cH$, we write
${\rm Rep}(\alpha,\cH)$ for the set of covariant representations
$(\pi,U)$ of $(\cA, G, \alpha)$ on $\cH$.
In the case that we are given a concrete von Neumann algebra
or non-degenerate $C^*$-algebra
$\cM\subseteq \cB(\cH)$ and a unitary representation $U:G\to\U(\cH)$
which produces an action by $\alpha_g(M)=\Ad U_g(M) := U_g M U_g^*$
on $\cM$, we will also
speak of the pair $(\cM,U)$ as a covariant representation for $\alpha$.


\end{defn}
\begin{rem}
Given  $(\cA,G, \alpha)$ where $\alpha$ need not be strongly continuous, and
a covariant representation  $(\pi,U)\in{\rm Rep}(\alpha,\cH)$, then the map $U \: G \to \U(\cH)$
is strong operator continuous.
Then it follows from \cite[7.4.2]{Pe89} that
$(\pi(\cA)'',G, \beta)$, defined by
 $\beta(g)=\Ad(U(g))$ is a  $W^*$-dynamical system.

\end{rem}

\section{Cross representations}

\subsection{Review}
\label{Cross-review}

In \cite{GrN14} we examined the question of how to construct crossed products for
$C^*$-actions $(\cA, G, \alpha)$ which may be singular (i.e.\ $\alpha$
need not be strongly continuous
or $G$ need not be locally compact). In the following sections we want to examine such crossed products
 where we require a spectral condition for the covariant representations
produced by the crossed product.

We review the basic facts regarding cross representations, extracted from \cite{GrN14},
which is where proofs of the material below can be found. There are a few new results at the end,
which will be proven here.

First, we need to generalize  the concept of a group algebra (\cite{Gr05}).
\begin{defn} \mlabel{def:2.1a} Let $G$ be a topological group.
A {\it host algebra for $G$} is a pair
$({\cal L}, \eta)$, where  ${\cal L}$ is a $C^*$-algebra and
$\eta \: G \to \U(M({\cal L}))$ is a group homomorphism
such that:
\begin{itemize}
\item[\rm(H1)] For each non-degenerate representation $(\pi, {\cal H})$
of $\cL$, the representation $\tilde\pi \circ \eta$ of $G$ is continuous.
\item[\rm(H2)] For each complex Hilbert space
${\cal H}$, the corresponding map
$$ \eta^* \: \Rep({\cal L},{\cal H}) \to \Rep(G, {\cal H}), \quad
\pi \mapsto \tilde\pi \circ \eta $$
is injective.
\end{itemize}
We write $\Rep(G,{\cal H})_\eta$ for the range of $\eta^*$,
and its elements are called {\it $\cL$-representations of $G$} on~$\cH$.
If $U \: G \to \U(\cH)$ is a unitary representation of $G$ in the range of
$\eta^*$, we write $U_\cL$ for the unique corresponding representation of $\cL$
with $\tilde U_\cL \circ \eta = U$.
\\[3mm]
We call $({\cal L}, \eta)$  a {\it full host algebra} if, in addition, we have:
\begin{itemize}
\item[\rm(H3)] $\Rep(G,{\cal H})_\eta = \Rep(G,{\cal H})$
for each Hilbert space~${\cal H}$.
\end{itemize}
\end{defn}
Thus by (H2) and (H3), a full host algebra, when it exists, carries precisely the continuous
unitary representation theory of $G$, and if it is not full, it carries some subtheory of the
continuous unitary representations of $G$. In particular, if we want to impose
 a spectral condition, then we will specify a host algebra which is not full.
 In this case we define:\\[2mm]
If $({\cal L}, \eta)$ satisfies (H1) and (H2) and if for each Hilbert space~${\cal H}$
we specify a subset $\Rep_c(G,{\cal H})\subseteq\Rep(G,{\cal H}),$ we will
call $({\cal L}, \eta)$  a {\it full host algebra for the class of representations in all}
$\Rep_c(G,{\cal H})$ if we also have:
\begin{itemize}
\item[\rm(H3')] $\Rep(G,{\cal H})_\eta = \Rep_c(G,{\cal H})$
for each Hilbert space~${\cal H}$.
\end{itemize}
If $G$ is locally compact,
then let
$\cL=C^*(G)$ with the canonical map
$\eta \: G \to \U(M(C^*(G)))$ which is strictly continuous,
(i.e.\ continuous for the strict topology on $\U(M(C^*(G)))).$
This defines on $C^*(G)$ the structure of a full host algebra for
the class of continuous representations of $G$ (cf.\ \cite[Sect.~13.9]{Dix77}).

 If $\cL$ is a host algebra of $G$ and
$\cI \subeq \cL$ is a closed ideal, the quotient map $q \: \cL \to \cL/\cI$
induces a surjective homomorphism
$\tilde q \: M(\cL) \to M(\cL/\cI)$, and
$\tilde q \circ \eta_G \: G \to M(\cL/\cI)$ defines on
$\cL/\cI$ the structure of a host algebra for $G$.
By construction we then have a morphism from the host algebra
$\cL$ (with respect to~$G$) to the host algebra $\cL/\cI$ (with respect to $G$).
For example, if $G=\R$, $\cL=C^*(G)$, and $\cI$ consists of the $C^*$-subalgebra of $\cL$
generated by those $L^1\hbox{-functions}$ whose Fourier transforms are supported in ${(-\infty,0)}$,
then $\cL/\cI\cong C_0([0,\infty))$ will be a full host algebra for the
positive continuous unitary representations of $\R$.

We are now ready to define the main structure we are interested in:
\begin{defn} \mlabel{def:a.5} Let $G$ be a topological group,
and let $(\cL,\eta)$ be a host algebra for $G$
and $(\cA,G, \alpha)$ be a $C^*$-action.
We call a triple $(\cC, \eta_\cA, \eta_\cL)$
a {\it crossed product host} for $(\alpha,\cL)$ if
\begin{itemize}
\item[\rm(CP1)] $\eta_\cA \: \cA \to M(\cC)$ and
$\eta_\cL \: \cL \to M(\cC)$ are morphisms of $C^*$-algebras.
\item[\rm(CP2)] $\eta_\cL$ is non-degenerate, i.e.
$\eta_\cL(\cL)\cC$ is dense in $\cC$ (cf. Definition~\ref{def:1.1c}).
\item[\rm(CP3)] The multiplier extension $\tilde\eta_\cL \: M(\cL)
\to M(\cC)$ satisfies in $M(\cC)$ the relations
\[ \tilde\eta_\cL(\eta(g))\eta_\cA(A)\tilde\eta_\cL(\eta(g))^*
= \eta_\cA(\alpha_g (A))\quad \mbox{for all} \quad A \in \cA,\;\hbox{and}\; g \in G.\]
\item[\rm(CP4)] $\eta_\cA(\cA) \eta_\cL(\cL) \subeq \cC$
and $\cC$ is generated by this set as a $C^*$-algebra.
\end{itemize}
We call the crossed product host for $(\alpha,\cL)$
{\it full } if it also satisfies
\begin{itemize}
\item[\rm(CP5)] For every covariant representation
$(\pi, U)$ of $(\cA,\alpha)$ on $\cH$ for which $U$ is an $\cL$-representation
of $G$, there exists a unique
representation $\rho \: \cC \to  \cB(\cH)$ with
\[ \rho(\eta_\cA(A)\eta_\cL(L)) = \pi(A) U_\cL(L)
\quad \mbox{ for } \quad A \in \cA, L \in \cL.\]
(Note that $(\cL,\eta)$ need not be full for
 $(\cC, \eta_\cA, \eta_\cL)$ to be full).
\end{itemize}
Two crossed product hosts $(\cC^{(i)}, \eta_\cA^{(i)}, \eta_\cL^{(i)})_{i=1,2}$
are said to be
{\it isomorphic} if there is an isomorphism $\Phi:\cC^{(1)}\to\cC^{(2)}$ such that
  ${\big(\Phi(\cC^{(1)}),\tilde\Phi\circ\eta_\cA^{(1)},\tilde\Phi\circ\eta_\cL^{(1)}\big)}
=(\cC^{(2)}, \eta_\cA^{(2)}, \eta^{(2)}_\cL)$.
\end{defn}

\begin{ex} \mlabel{ex:1.1}
Consider the usual case, i.e.\ a strongly continuous homomorphism
 $\alpha \: G \to \Aut(\cA)$, and
 $G$ is locally compact and we take  $\cL=C^*(G)$ as a full host algebra for $G$.
Then the crossed product algebra $\cC=\cA \rtimes_\alpha G$
is a full crossed product host for $(\alpha,\cL)$ in the sense above,
by the following reasoning.

From the usual construction (cf.  \cite[Thm.~7.7]{Pe89})
we already have the two homomorphisms
$\eta_\cA \: \cA \to M(\cC)$ and
$\eta_G \: G \to \U(M(\cC))$ (strictly continuous) such that
 \[
 \eta_\cA(\alpha_g A) = \eta_G(g) \eta_\cA(A) \eta_G(g)^*
\quad \mbox{ for } g \in G, A \in \cA\]
(also see \cite[Prop.~3]{Rae88}). Then the
strict continuity of $\eta_G \: G \to \U(M(\cC))$ leads by integration
to a morphism
$L^1(G) \to M(\cC)$ of Banach $*$-algebras. It therefore  extends
to a morphism $\eta_\cL \: C^*(G) \to M(\cC)$ of $C^*$-algebras
which is non-degenerate, and
$\eta_G = \tilde\eta_\cL \circ \eta$ follows.
It is easy to verify conditions (CP1)--(CP5),
hence the usual crossed product $\cC=\cA \rtimes_\alpha G$ is a
full crossed product host.
\end{ex}

A property of central importance for a crossed product host is that
it carries the covariant  $\cL$-representations of $(\cA, G, \alpha)$:

\begin{defn}
Given $(\cA, G, \alpha)$, where the action
$\alpha \: G \to \Aut(\cA)$ need not be
strongly continuous, assume a
host algebra $(\cL,\eta)$ for $G$. Then
a covariant representation
$(\pi, U)$ of $(\cA, G, \alpha)$ is called an
{\it $\cL$-representation} if $U$ is an $\cL$-representation
(cf.\ Definition~\ref{def:2.1a}). We write
$\Rep_\cL(\alpha, \cH)$ for the set of covariant
$\cL$-representations of $(\cA, G, \alpha)$ on $\cH$.
\end{defn}

\begin{thm} {\rm(\cite[Thm.~4.5]{GrN14})}
\mlabel{Thm-HostCross}
Let $(\cC, \eta_\cA, \eta_\cL)$ be a crossed product host
for $(\alpha,\cL)$, and recall the homomorphism
$\eta_G:=\tilde\eta_\cL\circ\eta:G\to \U(M(\cC))$.
Then, for each Hilbert space $\cH$, the map
$$ \eta^*_\times \: \Rep(\cC,{\cal H}) \to \Rep(\alpha, \cH),  \quad\hbox{given by}\quad
\eta^*_\times(\rho):=\big(\tilde\rho \circ \eta_\cA, \tilde\rho \circ \eta_G\big)$$
is injective, and its range
$\Rep(\alpha,{\cal H})_{\eta_\times}$ consists of $\cL$-representations
of $(\cA, G, \alpha)$.
If $\cC$ is full, then we also have
$\Rep(\alpha,{\cal H})_{\eta_\times} = \Rep_\cL(\alpha,{\cal H})$.
\end{thm}

\begin{thm} {\rm(Uniqueness Theorem; \cite[Thm.~4.8]{GrN14})}
Let $(\cL,\eta)$ be a host algebra for the topological group $G$,
and $(\cA, G, \alpha)$ be a $C^*$-action.
Let $(\cC, \eta_\cA, \eta_\cL)$ and
$(\cC^\sharp, \eta_\cA^\sharp, \eta_\cL^\sharp)$ be
crossed product hosts for $(\alpha,\cL)$,  such that
$\Rep(\alpha,{\cal H})_{\eta_\times} = \Rep(\alpha,{\cal H})_{\eta^\sharp_\times}$
for every Hilbert space $\cH$.
Then there exists a unique
isomorphism $\phi \: \cC \to \cC^\sharp$ with
$\tilde\phi \circ \eta_\cA = \eta_\cA^\sharp$ and
$\tilde\phi \circ \eta_\cL = \eta_\cL^\sharp$.
In particular, full crossed product hosts for $(\alpha,\cL)$ are isomorphic.
\end{thm}

\begin{thm} \mlabel{thm:5.1} {\rm(\cite[Thm.~5.1]{GrN14})}
Let $(\cL,\eta)$ be a host algebra
for the topological group $G$
and $(\cA,G, \alpha)$ be a  $C^*$-action.
\begin{itemize}
\item[\rm(a)] Let
$(\cC, \eta_\cA, \eta_\cL)$ be a  crossed product host for $(\alpha,\cL)$.
Then for the (faithful) universal representation $(\rho_u, \cH_u)$ of $\cC$,
the corresponding covariant $\cL$-representation
$(\pi, U)$ of $(\cA, G, \alpha)$ satisfies
\begin{eqnarray*}
 \rho_u(\eta_\cA(A)\eta_\cL(L)) &=& \pi(A) U_\cL(L)
\quad \mbox{ for } \quad A \in \cA, L \in \cL, \\[1mm]
\eta^*_\times(\rho_u)=\big(\tilde\rho_u \circ \eta_\cA, \tilde\rho_u \circ \eta_G\big)
&=& {(\pi,U)}\quad\hbox{ and  }\quad\rho_u(\al C.)=C^*\big(\pi(\cA) U_\cL(\cL)\big)\,.
\end{eqnarray*}
\item[\rm(b)] Conversely, let ${(\pi,U)}\in\Rep_\cL(\alpha,{\cal H})$ and put $\al C.:=C^*\big(\pi(\cA) U_\cL(\cL)\big)$.
Then $\pi(\cA)\cup U_\cL(\cL)\subset M(\cC)\subset  \cB(\cH)$, and we
obtain morphisms ${\eta_\cA \: \cA \to M(\cC)}$ and
$\eta_\cL \: \cL \to M(\cC)$ determined by
$\eta_\cA(A)C:=\pi(A)C$ and $\eta_\cL(L)C:=U_\cL(L)C$ for
$A\in\cA$, $L\in\cL$ and $C\in\cC$.
Moreover,  the following are equivalent:
\begin{itemize}
\item[\rm(i)]$ (\cC, \eta_\cA, \eta_\cL)$ is  a  crossed product host.
\item[\rm(ii)] $\pi(\cA) U_{\cL}(\cL) \subeq
{U_{\cL}(\cL) \cB(\cH)}$.
\item[\rm(iii)] For every  approximate identity
$(E_j)_{j\in J}$ of $\cL$ we have
\[ \| U_{\cL}(E_j) \pi(A) U_{\cL}(L)
- \pi(A) U_{\cL}(L) \| \to 0 \quad \mbox{ for } \quad
A \in \cA, L \in \cL.\]
\item[\rm(iv)] There exists an approximate identity
$(E_j)_{j\in J}$ of $\cL$ such that
\[ \| U_{\cL}(E_j) \pi(A) U_{\cL}(L)
- \pi(A) U_{\cL}(L) \| \to 0 \quad \mbox{ for } \quad
A \in \cA, L \in \cL.\]
\end{itemize}
\item[\rm(c)] Let $(\cC, \eta_\cA, \eta_\cL)$ be a  crossed product host
and $\Phi:\cC\to\cC/{\cal J}$ be a quotient map, where
${\cal J}$ is a closed two-sided ideal.
Then  ${\big(\cC/{\cal J},\tilde\Phi\circ\eta_\cA,\tilde\Phi\circ\eta_\cL\big)}$
is a crossed product host.
\end{itemize}
\end{thm}

The preceding theorem shows how crossed product hosts can be constructed.
It  also isolates a distinguished class of representations:

\begin{defn}
\mlabel{crossrepDef}
Let $\alpha \: G \to \Aut(\cA)$ be a $C^*$-action and $(\cL,\eta)$ be a host algebra for $G$. Then
a covariant $\cL$-representation
$(\pi, U)\in\Rep_\cL(\alpha, \cH)$ is called a
{\it cross representation} for $(\alpha,\cL)$
if any of the equivalent conditions (b)(i)--(iv) of
Theorem~\ref{thm:5.1} hold.
We write $\Rep_\cL^\times(\alpha, \cH)$ for the set of cross representations for  $(\alpha,\cL)$
on $\cH$.

Condition (b)(ii) is also equivalent to
\[
\pi(\cA) U_{\cL}(\cL) \subeq
\br U_{\cL}(\cL) \pi(\cA). \qquad\hbox{which is equivalent to}\qquad
\br\pi(\cA) U_{\cL}(\cL). =
\br U_{\cL}(\cL) \pi(\cA)..
\]
\end{defn}

\begin{rem} \mlabel{rem:3.8}
(i) Note that for a cross representation $(\pi, U),$ we get from these last
conditions that
\[
\al C.=C^*\big(\pi(\cA) U_\cL(\cL)\big)=\br\pi(\cA) U_{\cL}(\cL).\,.
\]
(ii)  From Remark~\ref{rem:factorize} we know that
$U_\cL(\cL)\cB(\cH)$ is a closed right ideal of $\cB(\cH)$, so that
\begin{equation}
  \label{eq:right-ideal}
\cR := \{ L \in \cL \mid \pi(\cA) U_\cL(L) \subeq U_\cL(\cL)\cB(\cH) \}
\end{equation}
is a closed right ideal of $\cL$ and condition (ii)
in Theorem~\ref{thm:5.1}(b)  is equivalent to
$\cR = \cL$. If an element $L \in \cL$ generates $\cL$ as a closed right ideal,
it follows in particular that the statement
$\cR = \cL$ is equivalent to $L \in \cR$.
\end{rem}

The  cross representations satisfy a range of permanence properties
(cf.~\cite[Prop.~5.3]{GrN14}), e.g.\ we have closure with respect to \
the taking of subrepresentations, arbitrary multiples,
and finite direct sums (but not infinite ones).
To express our existence theorem, we need the following:
\begin{defn} \mlabel{def:unicyc}
Cyclic representations of $\cA \rtimes_\alpha G_d$ are obtained from
states through the GNS construction. Let $\fS_\cL$
denote the set of those states
$\omega$ on $\cA\rtimes_{\alpha} G_d$ which thus produce
a covariant $\cL$-representation
$(\pi_\omega,U_\omega)\in{\rm Rep}_\cL(\alpha,\cH_\omega).$
This allows us to define the {\it universal  covariant
$\cL$-representation}
$(\pi_u,U_u)\in{\rm Rep}_\cL(\alpha,\cH_u)$ by
\[
\pi_u:= \bigoplus_{\omega\in\fS_\cL}
\pi_\omega,\quad U_u:= \bigoplus_{\omega\in\fS_\cL} U_\omega
\quad\hbox{on}\quad  \cH_u=\bigoplus_{\omega\in\fS_\cL}\cH_\omega.
\]
Clearly $\cH_u = \{0\}$ if $\fS_\cL = \emptyset.$
If $G$ is locally compact and $\cL=C^*(G)$ then
$(\pi_u,U_u)\in{\rm Rep}_\cL(\alpha,\cH_u)$ is just the
 universal  covariant representation
$(\pi_{co},U_{co})\in{\rm Rep}(\alpha,\cH_{co})$
defined before.
We will use $(\pi_u,U_u)$ below to prove the existence of
crossed product hosts.  We use the notation
 $U_{u,\cL} :=(\eta^*)^{-1}(U_u)\in \Rep(\cL,\cH_u)$
 for the associated representation of~$\cL$.
\end{defn}

\begin{thm} \mlabel{thm:exist} {\rm(Existence Theorem; \cite[Thm.~5.6]{GrN14})}
Let $(\cL,\eta)$ be a host algebra for $G$
and $\alpha \: G \to \Aut(\cA)$ be a $C^*$-action. Then the following are equivalent:
\begin{itemize}
\item[\rm(i)] There exists a full crossed product host
$(\cC, \eta_\cA, \eta_\cL)$ for $(\alpha,\cL)$.
\item[\rm(ii)] The universal covariant $\cL$-representation $(\pi_u, U_u)$
of $(\cA,G,\alpha)$ on $\cH_u$ is a cross representation.
\item[\rm(iii)]
$\;\Rep_\cL(\alpha, \cH)=\Rep_\cL^\times(\alpha, \cH)\;$ for all Hilbert spaces $\cH$.
\end{itemize}
\end{thm}
Below in Subsection~\ref{subsec:8.1} we give an example of a full crossed product host
when the usual crossed product does not exist.
(Unfortunately \cite[Example~5.9]{GrN14} contains an error, hence fails.)

Let $(\cL,\eta)$ be a host algebra for the topological group $G$
and $\alpha \: G \to \Aut(\cA)$ be a $C^*$-action. We assume that we are given a
set of cyclic covariant $\cL$-representations and form
their direct sum, denoted $(\pi^\oplus,U^\oplus)$.
Following  the construction in Theorem~\ref{thm:5.1}(b),
we put $\al C.:=C^*\big(\pi^\oplus(\cA) U^\oplus_\cL(\cL)\big)$,
which produces
 a triple ${(\cC, \eta_\cA, \eta_\cL)}$ such that
\begin{itemize}
\item[\rm(CP1)] $\eta_\cA \: \cA \to M(\cC)$ and
$\eta_\cL \: \cL \to M(\cC)$ are morphisms of $C^*$-algebras.
\item[\rm(CP4)] $\eta_\cA(\cA) \eta_\cL(\cL) \subeq \cC$
and $\cC$ is generated by this set as a $C^*$-algebra.
\end{itemize}
As $(\pi^\oplus,U^\oplus)$ need not be a cross representation,
(CP2) and (CP5) will fail in general.
If (CP2) fails,  the covariance requirement (CP3) does not make sense,
as it uses the multiplier extension $\tilde\eta_\cL \: M(\cL)\to M(\cC)$.
Covariance will have to be expressed differently,
and our first task is to obtain an adequate covariance condition
to replace (CP3).

\begin{defn}\mlabel{def:CP3'}
Assume (CP1) and (CP4) for $(\cC, \eta_\cA, \eta_\cL)$ in the context above.
\begin{itemize}
\item[\rm(i)]
For any Hilbert space $\cH$, a (non-degenerate) representation $\rho\in\Rep(\cC, \cH)$ is called an
{\it $\cL$-representation} if ${\tilde\rho\circ\eta_\cL:\cL\to  \cB(\cH)}$ is a
non-degenerate representation of $\cL$.  We write
$\Rep_\cL(\cC, \cH)$ for the set of
$\cL$-representations of $\cC$ on $\cH$.\\
 For each $\rho\in\Rep_\cL(\cC, \cH)$
there is a unitary $\cL$-representation $U^\rho:G\to \U(\cH)$,
which is uniquely specified by
$U^\rho(g)\cdot(\tilde\rho\circ\eta_\cL)(L)
=(\tilde\rho\circ\eta_\cL)(\eta(g)L)$ for all $g\in G,$ $L\in\cL$.
\item[\rm(ii)]
 Let $\fS_\cL(\cC):=\{\omega\in\fS(\cC)\,\mid\,\rho_\omega\in
\Rep_\cL(\cC, \cH_\omega)\}$
denote the set of those states
whose GNS representations $\rho_\omega$ are  $\cL$-representations of $\cC$.
We define the {\it universal
$\cL$-representation of~$\cC$},  $\rho_u^\cC\in{\rm Rep}_\cL(\cC,\cH_u^\cC)$ by
\[
\rho_u^\cC:= \bigoplus_{\omega\in\fS_\cL(\cC)}
\rho_\omega
\quad\hbox{on}\quad  \cH_u^\cC=\bigoplus_{\omega\in\fS_\cL(\cC)}\cH_\omega.
\]
When there is no danger of confusion, we will omit the superscript $\cC$.
\item[\rm(iii)]
The {\it covariance condition} is
given by assuming that
the universal $\cL$-representation  $\rho_u$ is faithful and satisfies
\begin{equation}
U^{\rho_u}(g)\cdot (\tilde\rho_u\circ\eta_\cA)(A)
\cdot U^{\rho_u}(g)^*=(\tilde\rho_u\circ\eta_\cA)\big(\alpha_g(A)\big)
\quad\mbox{for all}\quad g\in G,  A\in\cA.\tag{CP3'}
\end{equation}
\end{itemize}
\end{defn}

Any triple $(\cC, \eta_\cA, \eta_\cL)$ obtained
as in Theorem~\ref{thm:5.1}(b) from a covariant $\cL$-representation will satisfy condition (CP3').
Moreover, if we replace $\rho_u$ by any other faithful $\cL$-representation, the resulting covariance condition
will be equivalent to (CP3').

In the context of assuming (CP1) and (CP4) for $(\cC, \eta_\cA, \eta_\cL)$,
consider the hereditary $C^*$-subalgebra
of $M(\cC)$ generated by $\eta_\cL(\cL)\subset M(\cC)$.
It is
\[
M(\cC)_\cL:=\eta_\cL(\cL) M(\cC)\eta_\cL(\cL)=\eta_\cL(\cL) M(\cC)\cap M(\cC)\eta_\cL(\cL)
\]
We put $\cC_\cL:=M(\cC)_\cL\cap\cC=\eta_\cL(\cL) \cC\eta_\cL(\cL)$, which is a hereditary subalgebra of $\cC$.
\begin{defn} \mlabel{def:al}
Given a triple $(\cC, \eta_\cA, \eta_\cL)$ satisfying (CP1) and (CP4)
and $\cC_\cL$ as  above, let
\begin{equation}
  \label{eq:AL}
\cA_\cL:=\big\{A\in \cA\mid \eta_\cA(A)\cC_\cL\subseteq\cC_\cL\quad\hbox{and}\quad \eta_\cA(A^*)\cC_\cL\subseteq\cC_\cL\big\}\,.
\end{equation}
This includes the commutant of $\eta_\cL(\cL)$ in $\cA$.
\end{defn}

\begin{prop} \mlabel{propALchar} {\rm(\cite[Prop.~7.9]{GrN14})}
For a $C^*$-action $(\cA, G, \alpha)$
and a triple $(\cC, \eta_\cA, \eta_\cL)$ satisfying
{\rm(CP1)} and {\rm(CP4)}, the following assertions hold:
\begin{itemize}
\item[\rm(i)] $\cA_\cL = \{ A \in \cA \mid
\eta_\cA(A)\eta_\cL(\cL) \subeq \eta_\cL(\cL)\cC\;\;\hbox{and}\;\;
\eta_\cA(A^*)\eta_\cL(\cL) \subeq \eta_\cL(\cL)\cC\}$,
\item[\rm(ii)] For any (hence all) approximate identities $(E_j)_{j \in J}$
of $\cL$, we have\\
$\cA_\cL =\left\{ A \in \cA \,\big|\,
\big(\eta_\cL(E_j) -\1\big)\eta_\cA(B) \eta_\cL(L)
\to 0 \;\;\mbox{ for} \;\;B = A\;\hbox{and}\;A^*,\;\hbox{and for all}\;\;L \in \cL  \right\}$.
\item[\rm(iii)] In addition, assume {\rm(CP3')}. Then $\cA_\cL$ is an  $\alpha_G$-invariant subalgebra of $\cA$.
Moreover $\cA_\cL$ contains all the elements of $\cA$ which are invariant
with respect to  $\alpha_G$.
\end{itemize}
\end{prop}

\begin{cor}\mlabel{corAisAL} {\rm(\cite[Cor.~7.10]{GrN14})}
With $\alpha \: G \to \Aut(\cA)$ and a triple $(\cC, \eta_\cA, \eta_\cL)$ satisfying
{\rm(CP1)} and {\rm(CP4)} as above, we have
\begin{itemize}
\item[\rm(i)]
 $\cA_\cL=\cA$ if and only if $\;\cC_\cL=\cC$
\item[\rm(ii)] If  $(\cC, \eta_\cA, \eta_\cL)$ is
 constructed from the universal covariant $\cL$-representation
of $(\cA, G, \alpha)$, then a full crossed product host exists
if and only if   $\cA_\cL=\cA$.
\end{itemize}
\end{cor}

The preceding corollary implies that:
\begin{cor} \mlabel{cor:3.16}
A covariant
$\cL$-representation   ${(\pi,U)}\in\Rep_\cL(\alpha, \cH)$ is a cross
representation if and only if $\pi(\cA)=\pi(\cA)_\cL.$
\end{cor}

\begin{prop} \mlabel{contSpace} {\rm(\cite[Prop.~8.2]{GrN14})}
Let  $(\cL,\eta)$ be a host algebra for a topological group $G$ such that
the multiplier action $\eta \: G \to \U(M({\cal L}))$
is strictly continuous
and let $(\cA, G, \alpha)$
be a $C^*$-action for which we have a
triple $(\cC, \eta_\cA, \eta_\cL)$ satisfying
{\rm(CP1)}, {\rm(CP3')} and {\rm(CP4)}. We define $\cA_\cL$ as in \eqref{eq:AL}.
Then the following assertions hold:
\begin{itemize}
\item[\rm(i)] The subspace $\eta_\cL(\cL)\cC$ is contained in the closed
right ideal of $\cC$:
\[ \cC_c^L := \{ C \in \cC \mid \lim_{g \to \1} U^{\rho_u}(g)\cdot\rho_u(C)= \rho_u(C)\}.\]
\item[\rm(ii)] For $A \in \cA$ we have
$\eta_\cA(A) \eta_\cL(\cL) \subeq \cC_c^L$
if and only if for each $L \in \cL$ the map
\[ G \to \cC, \quad g \mapsto \eta_\cA(\alpha_g(A))\eta_\cL(L)\]
is continuous.
\item[\rm(iii)] In addition, let $G$ be locally compact and
$\cL = C^*(G)$.
Then $\eta_\cL(\cL) \cC = \cC_c^L$ and
\[ \cA_\cL = \big\{ A \in \cA \mid
\lim_{g \to \1} \eta_\cA(\alpha_g(B)) \eta_\cL(L) = \eta_\cA(B) \eta_\cL(L)
\mbox{ for all } L \in \cL\mbox{ and }B \in\{ A,\, A^*\}\big\}.\]
\end{itemize}
\end{prop}

By part~(iii) we have that $\cA_c\subeq \cA_\cL$ if $G$ is locally compact and $\cL = C^*(G)$.
In fact, from~(iii) we also get another characterization of $\cA_\cL$:
\begin{cor}
Let $G$ be locally compact, $\cL = C^*(G)$, and
$(\cA, G, \alpha)$ be a $C^*$-action. Then
\[
A\in\cA_\cL\quad\hbox{if and only if}\quad \eta_\cA(A) \eta_\cL(\cL)\cup\eta_\cL(\cL)\eta_\cA(A)\subseteq\cA_c\,.
\]
\end{cor}
\begin{prf}
Let $A\in\cA_\cL,$ then by Proposition~\ref{contSpace}(iii), we have
$\lim\limits_{g \to \1} \alpha_g(B) L = BL$
for all $ L \in \cL$ and $B \in\{ A,\, A^*\},$ where the actions $\eta_\cA,\;\eta_\cL$ are
understood in the products. Then for all $ L \in \cL$
\[
\left\|\alpha_g(AL)-AL\right\|\leq\|\alpha_g(A)\|\|\alpha_g(L)-L\|+\left\|\big(\alpha_g(A)-A\big)L\right\|.
\]
As the right hand side approaches zero as $g \to \1$, this implies that $A\cL\subset\cA_c,$
and likewise we get that $\cL A\subset\cA_c.$

Conversely, assume that $A\cL\cup\cL A\subseteq\cA_c,$ hence $\lim\limits_{g \to \1} \alpha_g(B L) = BL$
for all $ L \in \cL$ and $B \in\{ A,\, A^*\}.$ Then
\[
\left\|\alpha_g(A)L-AL\right\|\leq\|\alpha_g(A)\|\|L-\alpha_g(L)\|+\left\|\alpha_g(AL)-AL\right\|.
\]
As the right hand side approaches zero as $g \to \1$, this implies that
$\lim\limits_{g \to \1} \alpha_g(A) L = AL$
and likewise we also get this limit for $A^*$, hence $A\in\cA_\cL.$
\end{prf}
\begin{cor}\mlabel{corContX} {\rm(\cite[Cor.~8.4]{GrN14})}
Let $G$ be locally compact, $\cL = C^*(G)$, and
$(\cA, G, \alpha)$ be a $C^*$-action.  If  $(\cC, \eta_\cA, \eta_\cL)$ is
 constructed from the universal covariant $\cL$-representation
of $(\cA, G, \alpha)$, then the following are equivalent:
\begin{itemize}
\item[\rm(i)]  A full crossed product host exists.
\item[\rm(ii)] $\lim_{g \to \1} \eta_\cA(\alpha_g(A)) \eta_\cL(L) = \eta_\cA(A) \eta_\cL(L)$
for $A\in\cA$ and $L \in \cL$.
\item[\rm(iii)]  The conjugation action of $G$ on $\cC$ is strongly continuous.
\end{itemize}

These conditions imply that the maps $G \to M(\cC),
g\to\eta_\cA(\alpha_g(A))$ are continuous with respect to  the strict topology of $\cC$
for all $A\in\cA$.
If $\cA$ is unital, then {\rm(i)--(iii)} are  equivalent to
\begin{itemize}
  \item[\rm(iv)] For every $A \in \cA$, the map
$G \to M(\cC), g \mapsto \eta_\cA(\alpha_g(A))$ is strictly continuous.
\end{itemize}
\end{cor}

\begin{cor}\mlabel{corRCont} {\rm(\cite[Cor.~8.6]{GrN14})}
For  $G=\R$ and  $\cL=C^*(\R)$, let $(\cA, \R, \alpha)$
be a $C^*$-action for which we have a
triple $(\cC, \eta_\cA, \eta_\cL)$ satisfying
{\rm(CP1)}, {\rm(CP3')} and {\rm(CP4)}.  For a Hilbert space~$\cH$, consider
the injection:
 $$ \eta^*_\times \: \Rep_\cL(\cC,{\cal H}) \to \Rep_\cL(\alpha, \cH),  \quad\hbox{given by}\quad
\eta^*_\times(\rho):=\big(\tilde\rho \circ \eta_\cA, U^\rho\big)=: \big(\pi^\rho, U^\rho\big).$$
Let $\rho\in \Rep_\cL(\cC,{\cal H})$ be faithful, and denote the
spectral measure of $U^\rho$, resp., its infinitesimal generator
by $P_\rho$. Then an element $A \in \cA$ belongs to $\cA_\cL$ if and only if
\[
 \lim_{t \to \infty} P_\rho([-t,t])\pi^\rho(B)P_\rho([-s,s]) = \pi^\rho(B)P_\rho([-s,s])
\quad \mbox{ for } \quad s \in \R_+, B \in\{ A,\, A^*\}.\]
\end{cor}
For the case that $\cA= \cB(\cH)$ in the identical representation, and $P$ is the spectral measure of
a selfadjoint operator $H$ on $\cH$, and the action is $\alpha_t:=\Ad(e^{itH})$, then this
 produces the convenient formula
\[ \cA_\cL = \big\{ A \in  \cB(\cH) \;\big|\;
(\forall s \in \R_+) \lim_{t \to \infty} P[-t,t]BP[-s,s]
= BP[-s,s]\ \ \mbox{ for }\ \  B \in\{ A,\, A^*\}
\big\}, \]
to calculate $\cA_\cL$.

In the case that $G$ is a finite dimensional Lie group, we know that,
for every $\lambda > 0$, the resolvent
\[ R_\lambda=(\lambda\1-\Delta)^{-1} = \int_0^\infty e^{-\lambda t} e^{t\Delta}\, dt
\quad \lambda>0 \]
of the Laplacian $\Delta<0$  on $L^2(G)$ is in $\cL=C^*(G)$.
To see this,  recall from~\cite[Thm.~3.4]{Hu74} that
the positive semigroup $(e^{t\Delta})_{t>0}$ generated by the Laplacian is
represented by convolution with functions $p_t \in L^1(G)$ satisfying
$\|p_t\|\leq 1$ (the heat kernel),
hence $R_\lambda$ is given by convolution with the function $\int_0^\infty e^{-\lambda t} p_t\, dt \in L^1(G)$.

\begin{lem} \mlabel{lem:3.20} For every $\lambda > 0$, the right ideal
$R_\lambda C^*(G)$ is dense.
\end{lem}

\begin{prf} In view of \cite[Thm.~2.9.5]{Dix77}, we have to show that
no state $\omega$ of $C^*(G)$ vanishes on~$R_\lambda$. This follows from the fact that,
for every
continuous unitary representation $(U,\cH)$ of $G$ and $0 \not=\xi \in \cH$, we have
\[ \la \xi, U(R_\lambda)\xi \ra
=  \la\xi, (\lambda\1 - \dd U(\Delta))^{-1} \xi \ra > 0,\]
which
in turn follows from $\dd U(\Delta) \leq 0$, so that $\lambda \1 - \dd U(\Delta)$ is
strictly positive (see Definition~\ref{def:1.1a} for the derived
representation~$\dd U$).
\end{prf}

We may thus identify $R_\lambda$ with an element of
$C^*(G),$ and this is very useful, in that the cross condition only needs to be checked against $R_\lambda$
by the next proposition:

\begin{prop}\mlabel{ResLap}
Given a $C^*$-action $\alpha:G\to \Aut(\cA)$, where $G$ is a finite dimensional Lie group,
fix the host $\cL=C^*(G)$ and
let ${(\pi,U)}\in\Rep_\cL(\alpha,{\cal H}).$
For $\lambda > 0$, let $R_\lambda \in C^*(G)$ be the element representing the resolvent of the Laplacian.
Then ${(\pi,U)}$ is a cross representation for $(\alpha,\cL)$ if and only if
\[ \pi(\cA) U_{\cL}(R_\lambda) \subeq {U_{\cL}(\cL) \cB(\cH)}
\quad \mbox{ for some/any } \quad \lambda>0.\]
\end{prop}

\begin{prf} 
Consider the closed right ideal
$\cR = \{ L \in \cL\,\mid\, \pi(\cA) U_\cL(L) \subeq U_\cL(\cL)\cB(\cH) \}$
of $\cL$ from \eqref{eq:right-ideal}.
Then $(\pi, U)$ is a cross representation if and only if $\cL = \cR$.
By Lemma~\ref{lem:3.20}, this is equivalent to $R_\lambda \in  \cR$ for some
$\lambda >0$, and this completes the proof.
\end{prf}

This generalizes to replacing $R_\lambda$ by any selfadjoint element $E\in\cL$ for which multiplication
of $\cL$ by $E$ is non-degenerate. This in turn
generalizes to any topological group $G$ and any host algebra~$\cL$.

An interesting special case is:
\begin{thm} \mlabel{compactExist} {\rm(\cite[Thm.~6.1]{GrN14})}
Let $(\cA, G, \alpha)$ be a  $C^*$-action and
$(\cL,\eta)$ be a host algebra for $G$.
If ${(\pi,U)}\in\Rep_\cL(\alpha,{\cal H})$ satisfies
$\pi(\cA)U_{\cL}(\cL)\subseteq\cK(\cH)$, then
 $(\pi, U)$ is a cross representation for $(\alpha,\cL)$.
 This holds in particular if $U_{\cL}(\cL)\subseteq\cK(\cH)$.
\end{thm}

\begin{prop} 
{\rm(\cite[Lemma~C.3]{GrN14})}
Let $(U,\cH)$ be a continuous unitary
representation of the locally compact abelian group $G$ and
$U_{C^*(G)} \: C^*(G) \to  \cB(\cH)$ be the
associated integrated representation. Then
\[
U_{C^*(G)}(C^*(G))\subseteq\cK(\cH)
\]
if and only if the spectral measure $P$ of $U$ is a locally finite sum of point measures
with finite-dimensional ranges.
For $G = \R$ and $U_t = e^{itA}$, this condition is equivalent to the compactness of the
resolvent $(A + i \1)^{-1}$.
\end{prop}


Let $\omega$ be an invariant state of $\cA$ for
a  given $C^*$-action $(\cA, G, \alpha)$ and
$(\pi_\omega, U^\omega)$ be the corresponding covariant pair, where
$U^\omega: G\to\U(\cH_\omega)$ is the GNS unitary group determined by
\[ U^\omega_g \pi_\omega(A)\Omega_\omega
= \pi_\omega(\alpha_g(A))\Omega_\omega \quad \mbox{ for } \quad
g \in G, A \in \cA.\]
For any topological group $G$, we can now ask for conditions
on $\omega$ that ensure the continuity of the unitary representation
$U^\omega$, and if $G$ is locally compact, we can further try to see
when $(\pi_\omega, U^\omega)$ is a cross representation with respect to
$\cL = C^*(G)$. To formulate these characterizations, we write
\[
(\cA^*)_c := \{ \omega \in  \cA^*\mid {\lim\limits_{g\to \1}\|\alpha_g^*\omega-\omega\|}=0\}
\]
for the closed subspace of $\cA^*$ consisting of the $\alpha$-continuous
elements.

\begin{prop}
For a $C^*$-action $\alpha:G\to \Aut(\cA)$, where $G$ is locally compact,
and a $G$-invariant state $\omega$ of $\cA$, the following assertions hold:
\begin{itemize}
\item[\rm(a)] $(\pi_\omega,U^\omega)$ is a covariant representation,
i.e., $U^\omega$ is continuous, if and only if
$\cA\omega\cA\subseteq(\cA^*)_c$.
\item[\rm(b)] If {\rm(a)} is satisfied, then
$(\pi_\omega,U^\omega)$ is a covariant cross representation with respect to
the host $\cL=C^*(G)$ if and only if,
for all $A\in\cA$ and $f\in L^1(G)$, we have that
\[\lim_{g \to \1} \int_G \omega\big(B^*(\alpha_g(A)-A)\alpha_h(C)\big)\,f(h)\,dh = 0, \]
uniformly with respect to
$B$ and $C$ in the set $\{X \in \cA \,\mid\, \omega(X^*X) \leq 1\}$.
\end{itemize}
\end{prop}

\begin{prf}  (a) follows from \cite[Prop.~2.26(iii)]{BGN17}.

(b) Consider the version of the cross condition in Proposition~\ref{contSpace}(iii):
\[\lim_{g \to \1} \pi_\omega(\alpha_g(A)-A) U^\omega_\cL(L) = 0
\mbox{ for all } L \in \cL\mbox{ and }A \in\cA.\]
Here it is enough to let $L$ range over the dense subspace $L^1(G)\subset C^*(G)=\cL$.
Note that
\[
 \|\pi_\omega(\alpha_g(A)-A) U^\omega_\cL(L)\|
=\sup\Big\{\frac{|\la v,(\alpha_g(A)-A) U^\omega_\cL(L)w\ra |}{\|v\|\|w\|}\,\Big|\,\;
 v,\; w\in\cH_\omega\backslash\{0\}\Big\}.
\]
As $U^\omega_g\Omega_\omega=\Omega_\omega$ for all $g$,
we  have $U^\omega_g\pi_\omega(C)\Omega_\omega=\pi_\omega(\alpha_g(C))\Omega_\omega$,
 hence
\[ U^\omega_\cL(f)\pi_\omega(C)\Omega_\omega=\int_G f(h)\pi_\omega(\alpha_h(C))\Omega_\omega\,dh=:\pi_\omega(\alpha_f(C))\Omega_\omega
\quad \mbox{  for  } \quad f\in L^1(G),C\in\cA \]
by an abuse of notation (note that $\pi_\omega(\alpha_f(C))\in\pi_\omega(\cA)''$).
By setting
 $v=\pi_\omega(B)\Omega_\omega$ and  $w=\pi_\omega(C)\Omega_\omega$ for
 $B,C\in\cA$ and $L= f\in L^1(G)\subset\cL$, we get
 \[
 \la v,(\alpha_g(A)-A) U^\omega_\cL(f)w\ra=\omega\big(B^*(\alpha_g(A)-A)\alpha_f(C)\big).
 \]
Thus,
\begin{eqnarray*}
 &&\!\!\!\!\|\pi_\omega(\alpha_g(A)-A) U^\omega_\cL(f)\|   \\[1mm]
 &&=\sup\bigg\{
 \frac{|\omega\big(B^*(\alpha_g(A)-A)\alpha_f(C)\big)|}{\big[\omega(B^*B) \omega(C^*C)\big]^{1/2}}\;\bigg|\;
 B,\,C\in\cA, \  0 < \omega(B^*B), \omega(C^*C) \bigg\},
 \end{eqnarray*}
where we use the short-hand notation
\[\omega\big(B^*(\alpha_g(A)-A)\alpha_f(C)\big):=\int_G \omega\big(B^*(\alpha_g(A)-A)\alpha_h(C)\big)\,f(h)\,dh\,.\]
We conclude that the covariant pair  ${(\pi_\omega,U^\omega)}$ is a
 cross representation if and only if
\[
0=\lim_{g \to \1}\;\sup\bigg\{
 \frac{|\omega\big(B^*(\alpha_g(A)-A)\alpha_f(C)\big)|}{\big[\omega(B^*B) \omega(C^*C)\big]^{1/2}}\;\bigg|\;
 B,\,C\in\cA,\ 0 < \omega(B^*B),  \omega(C^*C) \leq 1 \bigg\}\]
holds for all $A\in\cA$ and  $f\in L^1(G)$, and this is condition~(b).
\end{prf}

Below we will also need the following variant of Corollary~\ref{corContX}(iii):

 \begin{prop}\mlabel{contC}
Given a $C^*$-action $\alpha:G\to \Aut(\cA)$ where $G$ is locally compact,
fix the host $\cL=C^*(G),$
let ${(\pi,U)}\in\Rep_\cL(\alpha,{\cal H})$ and put $\al C.:=C^*\big(\pi(\cA) U_\cL(\cL)\big)$
with morphisms
 ${\eta_\cA \: \cA \to M(\cC)}$ and
$\eta_\cL \: \cL \to M(\cC)$ as in {\rm Theorem~\ref{thm:5.1}}.
Then ${(\pi,U)}$ is a cross representation if and only if
the conjugation action $\Ad U$ of $G$ on $\cC$
 is strongly continuous.
\end{prop}

\begin{prf} If the  action $\Ad U$ of $G$ on $\cC$
 is strongly continuous, then as $\Ad U$ of $G$ on $U_{\cL}(\cL)$
 is also strongly continuous,
 it follows from
 \begin{align}\label{eq:triang}
&\pi(\alpha_g(A))U_{\cL}(L)  - \pi(A)U_{\cL}(L) \notag\\
&= \Big(U_g(\pi(A)U_\cL(L))U_g^{-1}  - \pi(A)U_{\cL}(L)\Big)
+ U_g\big(\pi(A)U_\cL(\alpha^\cL_{g^{-1}}(L)- L)\big)U_g^{-1}
 \end{align}
that
\begin{equation}
  \label{eq:limcond}
\lim_{g \to \1} \pi(\alpha_g(A))U_{\cL}(L) =
\pi(A)U_{\cL}(L)\quad \mbox{ for all } \quad A\in\cA, L\in\cL.
\end{equation}
By Proposition~\ref{contSpace}(iii), this implies that
 $\pi(\cA)=\pi(\cA)_{\cL},$ i.e.  ${(\pi,U)}$ is a cross representation.

If, conversely, ${(\pi,U)}$ is a cross representation, then
\eqref{eq:limcond} follows from Proposition~\ref{contSpace}(iii). As
 $\al C.:=C^*\big(\pi(\cA) U_\cL(\cL)\big),$ it suffices to show that
 \[\lim_{g \to \1} U_g\pi(A)U_{\cL}(L)U_g^* =\pi(A)U_{\cL}(L) \quad \mbox{ for all } \quad
A\in\cA, L\in\cL.\]
This follows from the identities \eqref{eq:triang},
\eqref{eq:limcond}, and the strong continuity of $\Ad U$ on $U_{\cL}(\cL)$.
\end{prf}

\subsection{Perturbations} 
\label{CrossPertb}

In general if one does not have the usual case, it is hard to find cross representations (hence crossed product hosts).
For a small class of situations, we gave
methods of finding and constructing cross representations
 in \cite{GrN14}. Here we want to continue with that line of enquiry, and in particular investigate whether cross representations are stable under perturbations.
This will allow us to extend the classes of $C^*$-actions known to have cross representations.
Perturbations of one-parameter groups is a large area of study, so there is much theory available to investigate the question.
In Sections~\ref{RepPosCros} and \ref{SpecCondCovRep} below we will consider how spectral conditions relate to cross representations.

Consider the one-parameter case, so  $U:\R\to\U(\cH), t \mapsto e^{itH},$
is a strong operator continuous unitary one-parameter group
such that $\alpha_t:=\Ad U_t$ defines  an action
 $\alpha:\R\to \Aut(\cM)$  on  a von Neumann algebra
$\cM\subseteq \cB(\cH)$.
Observe that if $\cL=C^*(\R)$, then
\[ \cM_c\subset\cM_{\cL}
=\{A\in\cM\,\mid\, \eta_{\cM}(B)\eta_{\cL}(\cL)\subseteq\eta_{\cL}(\cL)\cB(\cH)
\quad \mbox{ for } \quad B \in \{A,A^*\}\}.\]
If $\cM=\pi(\cA)''$ and $\alpha$ is not strongly continuous on $\cA$, then $\pi(\cA)$ need not be in $\cM_c$, though
it can still be in~$\cM_{\cL}$ (\cite[Ex.~5.11, Ex.~9.1]{GrN14}).

We first consider perturbations.
\begin{prop}
\label{crossperturb}
Given a concrete $C^*$-subalgebra $\cA\subseteq \cB(\cH)$, let $H$ and $B$ be selfadjoint operators such that $H+B$ is essentially selfadjoint,
let $U^{(0)}_t:=e^{itH}$ and  $U_t:={\exp(it(\oline{H+B}))}$, and assume that $\cA$ is preserved by both
$\alpha^{(0)}_t:=\Ad U^{(0)}_t$ and $\alpha_t:=\Ad U_t$.
Assume that either:
\begin{itemize}
\item[\rm(i)]  $B(i\1-H)^{-1}$ and $B(i\1-H-B)^{-1}$ are both bounded (e.g.\ if $B$ is bounded), or
\item[\rm(ii)] $H$ and $B$ are positive and
both $B(\1+H)^{-1}$ and $B(\1+H+B)^{-1}$ are bounded.
\end{itemize}
Then
\begin{eqnarray*}
\cA^{(0)}_\cL &:=& \big\{ A \in \cA \mid
BU^{(0)}_\cL(\cL) \subeq U^{(0)}_\cL(\cL) \cB(\cH)\quad\hbox{for}\quad B\in\{A,A^*\}\big\}\\[1mm]
= \cA_\cL &:=& \big\{ A \in \cA \mid
BU_\cL(\cL) \subeq U_\cL(\cL) \cB(\cH)\quad\hbox{for}\quad B\in\{A,A^*\}\big\}.
\end{eqnarray*}
Hence  $(\cA,U^{(0)})$ is  a cross representation if and only if $(\cA,U)$ is  a cross representation.
\end{prop}

\begin{prf}
(i) By the second resolvent identity we have
\[ (i\1-H-B)^{-1}-(i\1-H)^{-1}=(i\1-H-B)^{-1}B(i\1-H)^{-1}=(i\1-H)^{-1}B(i\1-H-B)^{-1},\]
so
\[ (i\1-H-B)^{-1}=(i\1-H)^{-1}+(i\1-H)^{-1}\big[B(i\1-H-B)^{-1}\big]\in (i\1-H)^{-1} \cB(\cH),\]
and
\[ (i\1-H)^{-1}=(i\1-H-B)^{-1}-(i\1-H-B)^{-1}\big[B(i\1-H)^{-1}\big]\in(i\1-H-B)^{-1} \cB(\cH).\]
Thus $(i\1-H)^{-1} \cB(\cH)=(i\1-H-B)^{-1} \cB(\cH),$ and likewise
\[ \cB(\cH)(i\1-H)^{-1}= \cB(\cH)(i\1-H-B)^{-1}.\]
If $A\in\cA^{(0)}_\cL$, then $A(i\1-H)^{-1}\in U^{(0)}_\cL(\cL) \cB(\cH)=
C^*((i\1-H)^{-1}) \cB(\cH)=\overline{(i\1-H)^{-1} \cB(\cH)}$ by definition, hence by the preceding
\[
A(i\1-H-B)^{-1}\in A(i\1-H)^{-1} \cB(\cH)\subseteq \overline{(i\1-H)^{-1} \cB(\cH)}=\overline{(i\1-H-B)^{-1} \cB(\cH)}=U_\cL(\cL) \cB(\cH).
\]
As this also holds for $A^*$, we have $A\in\cA_\cL$, i.e.\ $\cA^{(0)}_\cL\subseteq \cA_\cL$.
Likewise we get the converse inclusion, hence  $\cA^{(0)}_\cL= \cA_\cL$.

(ii) In the case that $H$ and $B$ are positive, then $U^{(0)}_\cL(\cL)=C^*((\1+H)^{-1})$ and
$U_\cL(\cL)={C^*((\1+H+B)^{-1})}$, so the rest of the proof follows by transcribing the one above for the replacement of
the resolvents $(i\1-H)^{-1},\;(i\1-H-B)^{-1}$ by $(\1+H)^{-1},\;(\1+H+B)^{-1}$ respectively.
\end{prf}
Thus cross representations are stable with respect to \ bounded perturbation of the generator of the dynamics.
A more interesting example is the following.

\begin{ex}
\label{NoGap}
We return to \cite[Ex.~9.3]{GrN14} to resolve a question which was left open.
There we proved that the Fock representation for a bosonic quantum field
is a cross representation for the dynamics induced by a second quantized
 one-particle Hamiltonian for which zero is isolated in its spectrum. We now prove that the
hypothesis that zero is isolated in its spectrum is unnecessary,
hence the Fock representation is
a cross representation for the second quantized dynamics produced by any (positive) Hamiltonian.

We start by recalling the notation.
Let $(\cH,\sigma)$ consist of a non-zero complex Hilbert space
 $\cH$ and the symplectic form
$\sigma:\cH\times\cH\to\R, \sigma(x,y):={\rm Im}{\langle x,y\rangle}$.
Note that $\U(\cH)\subset \Sp(\cH,\sigma)$.
The Weyl algebra $\cA=\ccr \cH,\sigma.$ carries an  action
$\alpha:\Sp(\cH,\sigma)\to\Aut(\cA)$ determined by $\alpha\s T.(\delta_x):=\delta_{T(x)}$.
The  Fock representation $\pi_F:\cA\to \cB(\cF(\cH))$ is given as follows.
 The bosonic Fock space is
\[
\cF(\cH):=\bigoplus_{n=0}^\infty\otimes_s^n\cH\,,\quad
\otimes_s^n\cH\equiv\hbox{symmetrized Hilbert tensor product of
$n$ copies of $\cH$}
\]
with the convention $\otimes_s^0\cH:=\C.$ The finite particle space
$\cF_0(\cH):={\Spann\{\otimes_s^n\cH\mid n=0,1,\cdots\}}$ is dense in $\cF(\cH)$.
For each $f\in\cH$ we define  on $\cF_0(\cH)$ a (closable)
 creation operator $a^*(f)$ by
\[
a^*(f)\, \big(v_1\otimes_s\cdots\otimes_s v_n\big):= \sqrt{n+1}\,
S\big(f\otimes v_1\otimes\cdots\otimes v_n\big)
=:\sqrt{n+1}\,f\otimes_s v_1\otimes_s\cdots\otimes_s v_n, \]
where $S$ is the symmetrizing operator. Define on $\cF_0(\cH)$ an essentially selfadjoint
operator by $\phi(f):=\big(a^*(f)+a(f)\big)/\sqrt{2}$
 where $a(f)$ is the adjoint of $a^*(f)$.
The Fock representation $\pi_F:\cA\to
\cB\big(\cF(\cH)\big)$ is then defined by $\pi_F(\delta_f)=\exp(i\overline{\phi(f)})$,
for all $f\in\cH,$ and it is irreducible.

Given a strong operator continuous one-parameter group
$U_t=\exp(itH)$ in $\U(\cH)$, where $H$ is selfadjoint,
define a unitary group
$\Gamma(U_t)$ in $\U\big(\cF(\cH)\big)$ by
\[ \Gamma(U_t)\big(v_1\otimes_s\cdots\otimes_s v_n\big):=
U_tv_1\otimes_s\cdots\otimes_s U_tv_n \]
which is strong operator continuous and whose
generator is given on $\cF_0(\cH)$ by
\[
\dd\Gamma(H)\big(v_1\otimes_s\cdots\otimes_s v_n\big)=
Hv_1\otimes_sv_2\otimes_s\cdots\otimes_s v_n+\cdots+v_1\otimes_s\cdots\otimes_s v_{n-1}\otimes_s Hv_n
\]
for  $v_j \in \cD(H).$
We then have the covariance $\pi_F(\alpha\s U_t.(A))=\Gamma(U_t)\pi_F(A)\Gamma(U_t)^*$.
If $P_n$ denotes the projection onto the $n$-particle space $\otimes_s^n\cH$,
 then $P_n$ commutes with
$\dd\Gamma(H)$ and $\Gamma(U)$ by construction.

In this example we show:
\begin{prop} \mlabel{prop:3.29}
If $H\geq 0$,  then  $(\pi_F,\Gamma(U))$ is  a
cross representation for $(\alpha,C^*(\R))$, where
$U_t = e^{itH}$ for $t \in \R$ and $\alpha \: \R \to \Aut(\cA)$
is given by $\alpha_t(\delta_x) = \delta_{U_t x}$ for $x \in \cH$.
\end{prop}

We already know that this is the case when
either $0\not\in\sigma(H),$ or
zero is isolated in the spectrum $\sigma(H).$
For a fixed operator $H\geq 0$ on $\cH$, decompose  $H=H_1+H_\infty$,
where $H_1:=\chi_{[0,1]}(H)\cdot H$
and $H_\infty:=H-H_1$. Then, for the spectra we have the
inclusions $\sigma(H_1)\subseteq [0,1]$, $\sigma(H_\infty)\subseteq \{0\}\cup[1,\infty)$,
and $H,\,H_1$ and  $H_\infty$ strongly commute. Now strong commutativity is equivalent to
the commutativity of the unitary one-parameter groups generated by these
selfadjoint operators (cf. \cite[Thm.~VIII.13]{RS80}). Hence, by the definition of
their second quantized unitary groups, these also commute, hence their generators
strongly commute, and these are $\dd\Gamma(H),\,\dd\Gamma(H_1)$ and  $\dd\Gamma(H_\infty).$
Define $U^{(\infty)}_t:=\exp(itH_\infty)$,
$\alpha_t:=\Ad \Gamma(U_t)$ and $\alpha^{(\infty)}_t:=\Ad \Gamma(U^{(\infty)}_t)$ on $\cM= \cB(\cF(\cH)).$ Let
\begin{eqnarray*}
\cM_\cL &:=& \big\{ M \in \cM \mid
B\cdot\Gamma(U)_\cL(\cL) \subeq \Gamma(U)_\cL(\cL)\cdot \cB(\cF(\cH))\\
&&\qquad\qquad\qquad\hbox{for}\quad B\in\{M,M^*\}\big\},\\[1mm]
\cM^{(\infty)}_\cL &:=& \big\{ M \in \cM \mid
B\cdot\Gamma(U^{(\infty)})_\cL(\cL) \subeq \Gamma(U^{(\infty)})_\cL(\cL)\cdot \cB(\cF(\cH))\\
&&\qquad\qquad\qquad\hbox{for}\quad B\in\{M,M^*\}\big\}.
\end{eqnarray*}
As $0$ is either isolated or not contained
in $\sigma(H_\infty)$, we know  by \cite[Ex.~9.3]{GrN14}  that $\pi_F(\cA)\subseteq\cM^{(\infty)}_\cL$, i.e.\ that
 $(\pi_F, \Gamma(U^{(\infty)}))$ is a cross representation for $(\alpha^{(\infty)},\cL)$.
We want to apply Proposition~\ref{crossperturb}(ii) to conclude that $\cM_\cL=\cM^{(\infty)}_\cL$, which will imply that
$(\pi_F, \Gamma(U))$ is a cross representation for $(\alpha,\cL)$ as well.
As $\dd\Gamma(H_\infty)$ and $\dd\Gamma(H_1)$ are positive we only need to show that
the two operators
\[
\dd\Gamma(H_1)(\1+\dd\Gamma(H_\infty))^{-1}\quad\hbox{and}\quad \dd\Gamma(H_1)(\1+\dd\Gamma(H) )^{-1}\,
\]
are bounded on $\cF(\cH)$. We start with the second one.
As
\[ 0\leq\dd\Gamma(H_1)\leq \dd\Gamma(H_1) + \dd\Gamma(H_\infty)=\dd\Gamma(H)\leq\1+\dd\Gamma(H), \]
 we have
 \[
 0\leq(\1+\dd\Gamma(H))^{-1/2}\dd\Gamma(H_1)(\1+\dd\Gamma(H))^{-1/2}=\dd\Gamma(H_1)(\1+\dd\Gamma(H))^{-1}\leq\1,  \]
  using strong commutativity, hence $\dd\Gamma(H_1)(\1+\dd\Gamma(H) )^{-1}\in \cB(\cF(\cH))$.

For the first operator, recall that  $P_n$ commutes with
$\dd\Gamma(H_1)$ and $\dd\Gamma(H_\infty)$ and $\1=\sum\limits_{n=0}^\infty P_n$ gives an orthogonal decomposition of $\cF(\cH)$.
Hence it suffices to show that the sequence ${\big\|\dd\Gamma(H_1)(\1+\dd\Gamma(H_\infty))^{-1}P_n\big\|}$ is uniformly bounded with respect to~$n$.
Recall from the \cite[Cor., p.~301]{RS80}
that, for any selfadjoint operator $B$ on $\cH$, we have
\[ \sigma(\dd\Gamma(B)P_n)
=\overline{\Big\{\sum\limits_{k=1}^n\lambda_k\mid\lambda_k\in\sigma(B)\Big\}}.\]
Thus $\|\dd\Gamma(H_1)P_n\|\leq n$. On the other hand, by $\sigma(H_\infty)\subseteq [1,\infty)$ we have\break
$\sigma((\1+\dd\Gamma(H_\infty))P_n)\subseteq {[n+1,\infty)}$, hence $\sigma((\1+\dd\Gamma(H_\infty))^{-1}P_n)\subseteq {[0,1/(n+1)]}$
by the Spectral Mapping Theorem,
so $\big\|(\1+\dd\Gamma(H_\infty))^{-1}P_n\big\|\leq 1/(n+1)$. Combining these:
\[
\big\|\dd\Gamma(H_1)(\1+\dd\Gamma(H_\infty))^{-1}P_n\big\|\leq \frac{n}{n+1} < 1
\]
which gives the desired uniform bound. Thus
\[
\dd\Gamma(H_1)(\1+\dd\Gamma(H_\infty))^{-1}\in \cB(\cF(\cH))
\] and so, via
Proposition~\ref{crossperturb}(ii), we get that
$(\pi_F, \Gamma(U))$ is a cross representation for $(\alpha,\cL)$.
\end{ex}


\begin{ex}
\label{h1h2}
This example continues \cite[Rem.~C.4]{GrN14}.
Let $h_1(\lambda)=\lfloor\lambda\rfloor$
and $h_2(\lambda):=\lambda - h_1(\lambda)$. Then we decompose
 $H=h_1(H)+h_2(H)$ and obtain  $\alpha_t=\alpha^{(1)}_t\circ\alpha_t^{(2)}$, where
$\alpha^{(j)}_t=\Ad U^{(j)}_t\in\Aut( \cB(\cH))$, $U^{(j)}_t:=\exp(ith_j(H))$ and
$U_t:=e^{itH}$.
As $\|h_2\|_\infty=1$, the operator $H$ is a bounded perturbation of $h_1(H)$, and so by
Proposition~\ref{crossperturb} we conclude that
\begin{eqnarray*}
 \cB(\cH)^{(1)}_\cL &:=& \big\{ A \in  \cB(\cH) \mid
BU^{(1)}_\cL(\cL) \subeq U^{(1)}_\cL(\cL) \cB(\cH)\quad\hbox{for}\quad B\in\{A,A^*\}\big\}\\[1mm]
=  \cB(\cH)_\cL &:=& \big\{ A \in  \cB(\cH) \mid
BU_\cL(\cL) \subeq U_\cL(\cL) \cB(\cH)\quad\hbox{for}\quad B\in\{A,A^*\}\big\}.
\end{eqnarray*}
Thus $\cA\subseteq  \cB(\cH)$ with $U_t$ is a cross
representation for $\alpha$ if and only if $\cA\subseteq  \cB(\cH)^{(1)}_\cL $,
and $\cA$ is preserved by $\alpha_t$ (but not necessarily $\alpha^{(1)}_t$).
This is very convenient, as $\alpha^{(1)}_{2\pi}=\id,$
so that  it is actually a representation of the circle group $\T \cong \R/2\pi \Z$.
If $P$ is the spectral measure of $H$ and $P_n := P[n,n+1)$, then
$U^{(1)}_{\cL}(\cL)=C^*\big(\{P_n \,\mid\,n\in\Z \}  \big)$. Thus
\begin{equation}
\label{CrossDiscrete}
  \cB(\cH)^{(1)}_\cL  = \big\{ A \in  \cB(\cH) \;\big|\;
(\forall k \in \Z) \lim_{n \to \pm\infty} P_n B P_k= 0,\quad B \in\{ A,\, A^*\}
\big\},
\end{equation}
which can in fact already be proven from Corollary~\ref{corRCont}.
Clearly, if $H\geq 0$ then  the
limit $n\to-\infty$ is omitted from the condition. This leads to the following matrix picture. We can write
 $A = (A_{jk})_{j,k\in \Z}$ as a matrix with
$A_{jk} =P_jAP_k$, and keep in mind that the convergence
$A=\sum\limits_{j\in\Z}\sum\limits_{k\in\Z}A_{jk}=\sum\limits_{j\in\Z}\sum\limits_{k\in\Z}P_jAP_k$ is in general
with respect to  the strong operator topology.
If we now form the matrix $M_A=(\|A_{jk}\|)_{j,k\in \Z}$, then the condition above states  that
$A\in  \cB(\cH)^{(1)}_\cL $ if and only if $\lim\limits_{j\to\pm\infty}\|A_{jk}\|=0=\lim\limits_{k\to\pm\infty}\|A_{jk}\|$
i.e., the real matrix $M_A$ has $c_0\hbox{--rows}$ and columns.
The Arveson spectral spaces can also be expressed in terms of the properties of this matrix,
as in \cite[Example~4.5]{BGN17}.
\end{ex}

In Proposition~\ref{crossperturb} above
 we saw that cross representations are stable with respect to bounded perturbations.
A natural question is then whether the crossed product hosts for the original and the perturbed actions
are the same.

 \begin{prop}
\label{PerturbCPH}
Given a concrete $C^*$-algebra $\cA\subseteq \cB(\cH)$, let $H$ and $B$ be selfadjoint operators such that $H+sB$ is essentially selfadjoint
for all $s\in(-\varepsilon,\varepsilon)$ and some $\varepsilon>0$.
Let $U^{(s)}_t:=\exp(it(\oline{H+sB}))$  and assume that $\cA$ is preserved by both
$\alpha^{(0)}_t:=\Ad U^{(0)}_t$ and $\alpha^{(s')}_t:=\Ad U^{(s')}_t$ for a fixed $s'\in(-\varepsilon,\varepsilon)$.
Assume that ${(\cA,U^{(0)})}$ (resp. ${(\cA,U^{(s')})}$)   is a cross representation for  $\alpha^{(0)}$
(resp.  $\alpha^{(s')}$)
with respect to $\cL=C^*(\R)$, so that we obtain the respective crossed product hosts
$\al C.^{(0)}:=C^*\big(\cA U^{(0)}_\cL(\cL)\big)$ and $\al C.^{(s')}:=C^*\big(\cA U^{(s')}_\cL(\cL)\big)$.
Assume that 
$B(i\lambda\un-H)^{-1}\in\al C.^{(0)}$   for some $\lambda\in\R\backslash \{0\}$
(e.g.~if $B\in\cA$).
\begin{itemize}
\item[\rm(i)]
Then, for all $s\in\R$ such that $\big\|sB(i\lambda\un-H)^{-1}\big\|<1$, we have:
\[
(i\lambda\un-H-sB)^{-1}=(i\lambda\un-H)^{-1}\sum_{n=0}^\infty\big(sB(i\lambda\un-H)^{-1}  \big)^n\in\al C.^{(0)},
\]
 and the series converges in norm.
\item[\rm(ii)] If $\big\|s'B(i\lambda\un-H)^{-1}\big\|<1$, then $\al C.^{(0)}\supseteq\al C.^{(s')}$.
\item[\rm(iii)] If we also have $B(i\lambda\un-H-s'B)^{-1}\in\al C.^{(s')}$   for some $\lambda\in\R\backslash \{0\}$ (e.g.~if $B\in\cA$) and
$\big\|s'B(i\lambda\un-H-s'B)^{-1}\big\|<1$ then $\al C.^{(0)}=\al C.^{(s')}$.
\end{itemize}
\end{prop}

\begin{prf}
(i) By applying \cite[Thm.~5.11]{We80} to $T:=i\lambda\un-H$ and $S:=-sB$, using the given hypotheses, we immediately obtain
the norm-convergent series
\[
(i\lambda\un-H-sB)^{-1}=(i\lambda\un-H)^{-1}\sum_{n=0}^\infty\big(sB(i\lambda\un-H)^{-1}  \big)^n\,.
\]
Since by hypothesis all terms of the series are in $\al C.^{(0)}$, by norm convergence so is the limit,
and hence $(i\lambda\un-H-sB)^{-1}\in\al C.^{(0)},$ as $(i\lambda\un-H)^{-1}\in  U^{(0)}_\cL(\cL)$. This proves (i).

(ii)
By part (i) we have $(i\lambda\un-H-s'B)^{-1}\in\al C.^{(0)}$ and hence
\[
\cA U^{(s')}_\cL(\cL) =    \cA\,C^*\big((i\lambda\un-H-s'B)^{-1}\big)\subset \al C.^{(0)},
\]
from which it follows that $\al C.^{(s')}=C^*\big(\cA U^{(s')}_\cL(\cL)\big)\subseteq \al C.^{(0)}.$

(iii) The hypotheses allow us to interchange in the inclusion in (ii) the operators  $H$ with $H+s'B$ to obtain the
reverse inclusion, hence equality.
\end{prf}
Given the hypothesis that $B(i\lambda\un-H)^{-1}\in\al C.^{(0)}$
 for some $\lambda\in\R\backslash \{0\}$, then by the first resolvent relation,
we also have for any $\mu\in\R\backslash \{0\}$ that
\[
B(i\mu\un-H)^{-1}=B\Big[ (i\lambda\un-H)^{-1}+i(\lambda-\mu)(i\lambda\un-H)^{-1}(i\mu\un-H)^{-1}  \Big]\in\al C.^{(0)}\,.
\]
By taking adjoints, we also get that  $(i\mu\un-H)^{-1}B\in\al C.^{(0)}$ for all $\mu\in\R\backslash \{0\}$.


In Example~\ref{h1h2} above, we had to decompose $H$
into more convenient parts. As this is a technique we will often use, we
prove the general lemma below to address this situation.
We have $H=h_1(H)+h_2(H)$, where $h_i\: \R \to \R$
are real Borel functions satisfying $h_1 + h_2 = \id_\R$.
Then $\alpha_t=\alpha^{(1)}_t\circ\alpha_t^{(2)}$ where
$\alpha^{(j)}_t=\Ad U^{(j)}_t\in\Aut( \cB(\cH))$ and $U^{(j)}_t:=\exp(ith_j(H)).$

\begin{lem}  \mlabel{Hdecomp}
With the notation from above, fix the host as $\cL=C^*(\R)$, and assume
for a concrete $C^*$-algebra that $\cA\subseteq  \cB(\cH)_{\cL}$ with respect to
both actions  $(\alpha^{(j)})_{j=1,\,2}$. If
\[
\{f(H)\mid f\in C_0(\sigma(H))\}\subseteq C^*\big(\{f_1(h_1(H))f_2(h_2(H))\;\mid\;
f_j\in C_0(\sigma(h_j(H))),\;j=1,\,2  \}\big),
\]
then ${(\cA,U)}$ is a cross representation for $\alpha:\R\to \Aut(\cA)$.
\end{lem}

\begin{prf}
From the hypotheses we obtain for all $A\in\cA$ and each
$f_j\in C_0(\sigma(h_j(H)))$ that
\[
Af_j(h_j(H))=g_j(h_j(H))B_j
\]
for some $g_j\in C_0(\sigma(h_j(H)))$ and $B_j\in \cB(\cH)$. Thus
\[
Af_1(h_1(H))f_2(h_2(H))=g_1(h_1(H))B_1f_2(h_2(H))=g_2(h_2(H))B_2f_1(h_1(H)), \]
hence $\lim\limits_{t\to 0} (U^{(j)}_t-\1)Af_1(h_1(H))f_2(h_2(H))=0$ for both $j=1,\,2.$ Using
\begin{align*}
\|(U_t-\1)D\|
&=\|(U^{(1)}_tU^{(2)}_t-\1)D\|
\leq\|U^{(1)}_t(U^{(2)}_t-\1)D\|+\|(U^{(1)}_t-\1)D\|\\
&=\|(U^{(2)}_t-\1)D\|+\|(U^{(1)}_t-\1)D\|
\end{align*}
for any $D\in \cB(\cH)$, we conclude that
\[ \lim\limits_{t\to 0}\big(U_t-\1)Af_1(h_1(H))f_2(h_2(H))=0, \quad \mbox{ hence }\quad
Af_1(h_1(H))f_2(h_2(H))\in U_{\cL}(\cL) \cB(\cH).\]
This also holds for
 $A^*\in\cA$, so by the hypothesis that the set
\[ \{f_1(h_1(H))f_2(h_2(H))\mid f_j\in C_0(\sigma(h_j(H)))\} \]
generates a $C^*$-algebra containing $U_{\cL}(\cL)$,
 we conclude that $\cA\subseteq  \cB(\cH)_{\cL}$ with respect to\
the action  $\alpha$, i.e.\ the given representation is a cross representation.
\end{prf}

\begin{rem}
We can use Proposition~\ref{crossperturb} to give a partial answer to the following
natural question.
 If $(\cA,U)$ and $(\cA,V)$ with $\cA\subseteq \cB(\cH)\supset U_{\R}\cup V_{\R}$
 are both covariant $\cL\hbox{-representations}$ for the same $\alpha$
and one is a cross representation, is the other one also
a cross representation?

It is true if $\cA= \cB(\cH)$ because
in this case there exists a $\lambda \in \R$ with
$U_t = e^{i\lambda t} V_t$ for all $t \in \R$. More concretely,
we recall from \cite[Ex.~5.11]{GrN14} that $U_t = e^{itH}$ defines a cross
representation for $\cA = \cB(\cH)$ if and only if $(i\1 - H)^{-1}$ is a compact
operator. By the resolvent formula, this property is stable under bounded
perturbations.

If we assume that $U_t:=e^{itH}$ and  $V_t:={\exp(it(\oline{H+B}))}$,
where $H$ and $B$ are selfadjoint operators such that $H+B$ is essentially selfadjoint,
then  Proposition~\ref{crossperturb} applies, and we get that
boundedness of  $B(i\1-H)^{-1}$ and $B(i\1-H-B)^{-1}$ implies that  $(\cA,U)$ is a cross representation if and only if $(\cA,V)$ is a cross representation.

For positive implementing groups $U_t$ and $V_t$, Theorem~\ref{CrossRepsSame} below further
analyzes the question.
\end{rem}


\subsection{Inner cross representations of $W^*$-algebras}
\label{InnerW}

Here we want to examine  cross representations for $W^*$dynamical systems.
As remarked before, the cross property of a covariant representation does not in general extend to
the $W^*$-dynamical system generated by the covariant representation. Recall the following example:

\begin{ex}
\label{CrossNoExt}
Let $\cA=\cK(\cH)$ and consider its identical representation $\pi$.
Further,   let  $\alpha_t=\Ad V_t$ for $V_t=e^{itH}$ for unbounded $H=H^*$,
 where $(i\1-H)^{-1}\not\in\cK(\cH)$. Then
 $(\pi,V)$ is a cross representation  of $\alpha:\R\to \Aut(\cA)$.
For $\cM=\pi(\cA)''= \cB(\cH)$,  the pair   $(\cM,V)$
is not a cross representation for $\beta_t:={\rm Ad}\,V_t$ on $\cM$
(cf.  \cite[Ex~5.11]{GrN14}), even though $(\pi,V)$ is so for~$\alpha$.
\end{ex}
In the converse direction the cross property is conserved, i.e. if
$(\cM,V)$ is a cross representation, then so is its restriction to
the ${\rm Ad}\,V\hbox{--invariant}$ subalgebra $\pi(\cA)$, where $\cM=\pi(\cA)''$
(cf. \cite[Prop.~5.3(iv)]{{GrN14}}).

It is well-known that the requirement for the algebra $\cM$ on which we have an action
$\alpha:\R\to\Aut(\cM)$,
to be a $W^*$-algebra, places strong restrictions on it and on its covariant representations.
For instance, if the action is strongly continuous, it must be
uniformly continuous, hence it is inner and implemented by a continuous unitary one-parameter group in
the $W^*$-algebra (cf. \cite[Ex.~XI.3.6]{Ta03}),
 or if a covariant representation is positive, then the action is inner in that representation,
 and implemented by a $W^*$-continuous unitary one-parameter group in
the $W^*$-algebra
 (cf.~Borchers--Arveson Theorem \cite[Thm.~3.2.46]{BR02}).
We will see likewise below in Theorem~\ref{W-CrossRepsChar} that the cross condition
places strong restrictions on an inner $W^*$-dynamical system. We will only consider the
case where the $W^*$-dynamical system is inner.

\begin{defn} Let $\cM$ be a $W^*$-algebra and
$(V_t)_{t \in \R}$ be a weakly continuous unitary one-parameter group
of $\U(\cM)$. Denote the inner action ${\rm Ad}\,V_t$ of $\R$ on $\cM$
by $\alpha_t$.

(a) We call a normal representation
$(\pi, \cH)$ of $\cM$ a {\it $V$-cross representation} if
the pair $(\pi, U)$ with $U_t := \pi(V_t)$ is a cross representation
for $(\cA,\cL, \alpha) = (\cM, C^*(\R), \alpha)$   in the sense of
Subsection~\ref{Cross-review}.

(b) We call a normal representation
$(\pi, \cH)$ of $\cM$ {\it $V$-bounded} if
$\pi \circ V$ is a norm-continuous one-parameter group.

(c)  We call a normal representation
$(\pi, \cH)$ of $\cM$ of {\it bounded type} if
it is an orthogonal direct sum of $V$-bounded representations.
\end{defn}

\begin{rem} Every $V$-bounded representation $(\pi,\cH)$
is a $V$-cross representation by continuity of the action $t\mapsto{\rm Ad}(V_t)$.
As $(V_t)_{t \in \R}\subset\U(\cM)$, the spectral projection map $P$ of $V$ has range in $\cM$,
hence a normal representation
$(\pi, \cH)$ of $\cM$ is  $V$-bounded if and only if there is some $C_\pi>0$ such that $\pi(P(-\infty,-C_\pi))
=0 =\pi(P(C_\pi,\infty))$. As $\pi(\cM)\cong Z_\pi\cM$ for some central projection $Z_\pi\in Z(\cM)$,
$\pi$ is $V$-bounded if and only if $Z_\pi P(-\infty,-C_\pi)=0=Z_\pi P(C_\pi,\infty)$. It is clear that a direct sum of
$V$-bounded representations need not be $V$-bounded, unless it is a finite direct sum.
Thus representations of bounded type need not be  $V$-bounded.
\end{rem}
\begin{prop}
\label{NormDecomp}
Given $(\cM,V)$, then any normal representation $(\pi,\cH)$ has a decomposition
\[ \pi=\pi_N\oplus\pi_T  \]
where $\pi_N$ is of bounded type, 
and where $\pi_T$ contains no
$V$-bounded subrepresentation  other than the trivial one.
Furthermore, $\pi_N(\cM)$ is an $\ell^\infty$-direct sum of ideals $(\cM^j)_{j \in \N}$ such that
$(\cM^j, V^j)$ has a faithful  normal representation for which
 $V^j_t =e^{-it H_j}$ for a bounded
operator $H_j \in \cM_j$. Here $V^j_t$ is the projection of $\pi_N(V_t)$ onto the
$j$-th component.
\end{prop}

\begin{prf} Let $\cH_N \subeq \cH$ denote the $\pi(\cM)$-invariant subspace
generated by all vectors $\xi \in \cH$ generating a $V$-bounded cyclic
subrepresentation. Then the invariant subspace $\cH_T := \cH_N^\bot$
contains no non-zero $V$-bounded subrepresentation.
An application of Zorn's Lemma now shows that the representation
$(\pi_N,\cH_N)$ is a direct sum of $V$-bounded cyclic subrepresentations.

For the last claim, we consider for every $j \in \N$ the maximal
$\pi(\cM)$-invariant subspace of $\cH$ on which the the infinitesimal generator
of $V$ has norm $\leq j$. Then $\bigcup_j \cH_j$ is dense in $\cH_N$, so that
$\cH_N$ is the Hilbert space direct sum of the subspaces $\cK_j := \cH_{j+1} \cap \cH_j^\bot$.
It follows from their definition, that the subspaces $\cH_j$ are also invariant
under the commutant $\pi(\cM)'$. This shows that the projections
$Z_j$ onto $\cH_j$ are central in $\pi(\cM)$, which leads to the $\ell^\infty$-direct sum
decomposition $\pi_N(\cM) \cong \oplus_{j = 1}^\infty Z_j\pi_N( \cM)$.
\end{prf}

 Using this proposition, we can state the main result of this subsection, in which we
 characterize the cross representations for this system. We will call a spectral projection
 $P(E)$ of a one-parameter group $V:\R\to \U(\cH)$ {\it finite} if
$E \subset \R$ is a bounded measurable set.

\begin{thm} \mlabel{W-CrossRepsChar}
For the pair $(\cM,V)$, a normal representation
\[ \pi=\pi_N\oplus\pi_T \quad\hbox{on}\quad\cH=\cH_N\oplus\cH_T
\]
decomposed as in {\rm Proposition~\ref{NormDecomp}}
defines a cross-representation $(\pi(\cM),\pi \circ V)$ for $\alpha$ if and only if
  both of the following conditions are satisfied:
\begin{itemize}
\item[\rm(i)] For every finite spectral projection
$P$ of $\pi_N(V)$, there exists a finite spectral projection $Q$ of $\pi_N(V)$ with
\begin{equation}
  \label{eq:scon}
 \pi_N(\cM) P \subeq Q \pi_N(\cM).
\end{equation}
\item[\rm(ii)]  $\pi_T(\cM)$ is an orthogonal, at most countable, direct sum of ideals
$\cM_j \cong  \cB(\cH_j)$, $j \in J$,
and the infinitesimal generator $H_j$ of the one-parameter group
$(V^j_t)_{t \in \R}$ in $\cM_j$, has compact resolvent.
Here $V^j_t$ is the projection of $\pi_T(V_t)$ onto the $j$-th component.
Moreover, every  finite spectral projection of $\pi_T(V)$ generates an ideal in $\pi_T(\cM)$
consisting of a finite direct sum of the ideals $\cM_j$.
\end{itemize}
  \end{thm}

This theorem is based on the special case $(\cB(\cH),V)$, for which we have seen in
\cite[Ex.~5.11]{GrN14} that the cross condition is equivalent to
the compactness of $(i\1 - H)^{-1}$ if $V_t = e^{itH}$ for $t \in \R$.
We also recall that the corresponding conjugation action of $\R$ on $\cB(\cH)$
is strongly continuous if and only if $H$ is bounded
\cite[Prop.~5.10]{GrN14} and \cite[Ex.~XI.3.6]{Ta03}.

  \begin{prf} {\it(first part)}
 By  \cite[Prop.~5.3]{{GrN14}}, the set of cross representations is closed
under subrepresentations and finite direct sums.
Thus $(\pi(\cM), \pi \circ V)$ is a cross representation if
and only if both $(\pi_N(\cM), \pi_N \circ V)$
and $(\pi_T(\cM), \pi_T \circ V)$ are cross representations.
Criteria (i) and (ii) are the conditions for
$(\pi_N(\cM), \pi_N \circ V)$
and $(\pi_T(\cM), \pi_T \circ V)$, resp., to be cross representations,
which we now prove.

 (i) Clearly \eqref{eq:scon} implies the cross condition for $\pi_N$.
Suppose that \eqref{eq:scon} is not satisfied. Then there exists a
$P$ such that the left multiplication action of $\pi_N(V_t)$ on $\pi_N(\cM) P$ is not norm continuous.
We thus find pairwise different elements
$j_n \in \N$ and bounded subsets $E_n \subeq \R$ with
$E_n \cap [-n,n] = \eset$ and $P(E_n)_{j_n} \cM_{j_n} P_{j_n} \not= \{0\}$,
where $P_j$ denotes the component of $P$ in $\cM_j$ in the decomposition
of Proposition~\ref{NormDecomp}.
Here we use that, for any finite subset $F = \{ j_1, \ldots, j_n\}\subeq \N$,
the right multiplication action on $\sum_{j \in F} \cM_j$ is norm continuous.

Pick $M_n \in P(E_n)_{j_n} \cM_{j_n} P_{j_n}$ with $\|M_n\| = 1$, so that the
sum $M := \sum_n M_n \in\pi_N( \cM)$ converges in the weak topology. Then
$M$ does not satisfy the cross condition because
\[ {\|M - P[-n,n]M\|} \not\to 0 \quad \mbox{ for } \quad n \to \infty.\]
 This proves (i)

 To prove (ii), we need the following definition and lemmas.
  \end{prf}
\begin{defn} If $\cM$ is a $W^*$-algebra, then the weakly closed
left ideals of $\cM$ are of the form $\cM P$, where
$P$ is a hermitian projection (\cite[Prop.~1.10.1]{Sa71}). We say that a normal
state $\omega \in \fS_n(\cM)$ is {\it supported by $\cM P$} if
\[ \cM(\1-P)  \subeq \ker \omega, \quad \mbox{ or, equivalently, } \quad
 \1 - P \subeq (\cM\omega)^\bot.\]
\end{defn}

\begin{lem} \mlabel{lem:4.3} Let $\omega \in \cM_*$ be a state supported on a
$\sigma(\cM,\cM_*)\hbox{--closed}$ left 
ideal $\cJ \subeq \cM$ on which the left multiplication action
of $\R$ given by $(M,t) \mapsto  V_tM$ is norm continuous.
Then $(\pi_\omega, \cH_\omega)$ is a $V$-bounded representation.
\end{lem}

\begin{prf} The assumption that $\R \times \cJ \to \cJ,
(t,M) \mapsto  V_tM$ is norm continuous, implies that there exist
real numbers $a < b$ such that
$\cJ =  P[a,b]\cJ$, where $P[a,b]$ is the corresponding
spectral projection of $V$ for which
$V_t = \int_\R e^{itx}\, dP(x)$.

Writing $\cJ =  \cM Q $ for a projection $Q$ (\cite[Prop.~1.10.1]{Sa71}), we have
 $\cM = \cJ \oplus  \cM(\1 - Q)$.
As $\omega$ is supported on $\cJ$ we have
 $\cM(\1 - Q)\subseteq N_\omega\subset{\rm Ker}(\omega)$ where $N_\omega$ is the left kernel of $\omega$.
 Thus $\xi_\omega(\cM)=\xi_\omega(\cJ)=\cM/N_\omega\subset\cH_\omega$ is dense, where
 $\xi_\omega:\cM\to\cH_\omega, M \mapsto \pi_\omega(M)\Omega_\omega$ denotes the GNS map.
As
\[
\pi_\omega(P[a,b])\xi_\omega(M)
= \xi_\omega( P[a,b]M) = \xi_\omega(M)\qquad\hbox{for all $M\in\cJ$,}
\]
 we obtain
 that $\pi_\omega(P[a,b]) = \1_{\cH_\omega}$, i.e.,
that $\pi$ is $V$-bounded.
\end{prf}

\begin{lem} \mlabel{lem:4.4}  Let $\cJ \subeq \cM$ be a non-zero weakly closed left ideal.
Then there exists a normal state $\omega \in \fS_n(\cM)$ supported by $\cJ$.
\end{lem}

\begin{prf} If $0 \not= M \in \cJ$, then
$M^*M \in \cJ$, and since this element is hermitian,
$0 \not= M^*M \in \cJ \cap \cJ^*$.
Let $P \in \cM$ be a projection with
$\cJ = \cM P$ (\cite[Prop.~1.10.1]{Sa71}). Then
$\cJ \cap \cJ^* = P \cM P$. If $\nu \in \fS_n(\cJ \cap \cJ^*)$,
then $\omega(M) := \nu(PMP)$ is a normal state of $\cM$ supported by $\cJ$.
\end{prf}

\begin{lem} \mlabel{lem:4.6} Let $\cM$ be a $W^*$-algebra.
\begin{itemize}
\item[\rm(i)] If $Q \in \cM$ is a minimal projection, then the weakly closed
(two-sided) ideal of $\cM$ generated by $Q$ is isomorphic to some $\cB(\cH)$.
\item[\rm(ii)] If $\cM$ is generated, as a two-sided weakly closed ideal of $\cM$, by
minimal projections, then
\[ \cM \cong \bigoplus_{j \in J}  \cB(\cH_j) \quad \mbox{ for Hilbert spaces }
\quad \cH_j.\]
\end{itemize}
\end{lem}

\begin{prf} (cf.\ \cite[p.~354]{Bla06}) (i)  Let $\cJ \subeq \cM$ be the weakly
closed ideal generated by a minimal projection $Q \in \cM$.
Let $0 \not = Z \in Z(\cJ)$ be a projection. Then
$Q = QZ + Q(\1- Z)$, so that $QZ \not=0$ implies $Q(\1 - Z)= 0$ by minimality.
Thus $Q = QZ$ so $Z$ is the identity of $\cJ$.
We conclude for  the center $Z(\cJ) =\C Z$, so that $\cJ$ is factor.
Since $\cJ$ contains a minimal projection,
 \cite[Cor.~I.8.3]{Dix81} implies that $\cJ \cong  \cB(\cH)$
for some Hilbert space~$\cH$ (it also follows from  \cite[Cor.~5.5.8]{Pe89}).

(ii) Let $\cM$ be generated, as a two-sided weakly closed ideal of $\cM$,
 by minimal projections.
Let $\{P_j\,\mid\,j\in J\}$ be a set
of minimal projections which generates $\cM$
as a two-sided weakly closed ideal of $\cM$. Let $\cJ_j \subeq \cM$ be the weakly
closed ideal generated by  $P_j \in \cM$. Then by \cite[Prop.~1.10.5]{Sa71},
there is a central projection $Z_j \in Z(\cM)$  with
$\cJ_j = \cM Z_j$. Clearly $Z_j$ is minimal, so that we either have
$Z_j Z_k = Z_j$ or $Z_j Z_k = 0$.
We conclude that there is a subset $J'\subseteq J$ such that
$ \cM \cong \bigoplus\limits_{j \in J'} \cJ_j$. The assertion
then follows from part (i).
\end{prf}

\begin{lem} \mlabel{lem:4.9}
Assume that $(\pi_T(\cM),\pi_T(V))$ satisfies the cross condition.
Let $P = P(E)$ be a spectral projection of $\pi_T(V)$, where
$E \subeq \R$ is a bounded measurable subset.
Then $P$ is a finite orthogonal sum of minimal projections in $\pi_T(\cM)$.
\end{lem}

\begin{prf} If $P$ is not a finite orthogonal sum of minimal projections in $\pi_T(\cM)$, then it can be written
as $P = \sum_{j = 1}^\infty P_j$, where the $P_j\in \pi_T(\cM)$ are non-zero pairwise orthogonal
projections.
Recalling that  $\pi_T$ contains no
$V$-bounded subrepresentation  other than the trivial one,
Lemmas~\ref{lem:4.3} and \ref{lem:4.4} imply that the left multiplication
actions $(t,M) \mapsto \pi_T(V_t)M$ of $\R$ on the left ideals $\pi_T(\cM) P_j$ are not norm continuous.
Hence we can inductively choose $a_j \in \R$, $j \in \N$,
such that the spectral projections $Q_j = P[a_j, a_j+1]$ satisfy
$Q_j \pi_T(\cM) P_j \not=\{0\}$ and $|a_{j+1}| > |a_j| + 1$.

Let $M_j \in Q_j \pi_T(\cM)  P_j$ be an element with $\|M_j\| = 1$.
Since the sequences of projections $(P_j)_{j \in \N}$ and $(Q_j)_{j\in \N}$
are mutually disjoint, the series
$M := \sum_{j = 1}^\infty M_j$ converges weakly in $\pi_T(\cM)$. In fact,
if 
$v \in \cH_T$, then
\[ \sum_{j = 1}^\infty \|M_jv\|^2
= \sum_{j = 1}^\infty \|M_j P_j v\|^2
\leq \sum_{j = 1}^\infty \|P_j v\|^2
\leq \|v\|^2.\]
Since the vectors $M_jv$ are also mutually orthogonal,
$Mv := \sum_j M_j v$ defines a bounded operator on $\cH_T$ with $\|M\| \leq 1$.
By construction, the element $M \in P\pi_T(\cM)$ does not satisfy the cross condition
\[ \lim_{n \to \infty} \|  M P-  P[-n,n]MP\| = 0. \]
This contradicts our assumption that $(\pi_T(\cM),\pi_T(V))$ satisfies the cross condition.
Thus $P$ is a finite orthogonal sum of minimal projections.
\end{prf}


\begin{prf} {\it (of Theorem~\ref{W-CrossRepsChar}(ii))}:\\
Assume that $(\pi_T(\cM),\pi_T(V))$ satisfies the cross condition.
 Let $\cJ \trile \pi_T(\cM)$ be the weakly closed ideal generated by the
finite spectral projections of $V$. Lemma~\ref{lem:4.9} further
implies that $\pi_T(\cM)$ is generated by minimal projections.
In view of
$\1 = P(\R) = \lim_{n \to \infty} P[-n,n]$, we then have $\1 \in \cJ$, i.e.,
$\pi_T(\cM) = \cJ$, so that Lemma~\ref{lem:4.6} implies that
\begin{equation}
\label{Mdecomp}
 \pi_T(\cM)  = \bigoplus_{j \in J} \cM_j\cong \bigoplus_{j \in J}  \cB(\cH_j)
 \end{equation}
for Hilbert spaces $\cH_j$. Now the corresponding unitary one-parameter group
$V^j_t \in \U(\cH_j)$ satisfies the cross condition for $\cB(\cH_j)$,
so that the compactness of the resolvent of $H_j$ follows from
\cite[Ex.~5.11]{GrN14}.

As each finite spectral projection $P$ of $\pi_T(V)$ is a finite orthogonal sum of minimal projections in $\pi_T(\cM)$,
we can build up a compatible set of minimal projections for the projections $ P[-n,n]$ such that the
set of  minimal projections of  $ P[-n,n]$ is contained in that of $ P[-n-1,n+1]$.
This can be done by starting from  $ P[-1,1]$ and proceeding inductively by letting the set of minimal projections of $ P[-n-1,n+1]$
be the union of that of $ P[-n,n]$ and of ${P\big([-n-1,-n)\cup(n,n+1]   \big)}$. By $\1 = P(\R) = \lim_{n \to \infty} P[-n,n]$,
this produces a countable resolution of the identity
\[
\1=\sum_{k = 1}^\infty P_k
\]
into minimal projections. If we use this resolution to obtain the decomposition (\ref{Mdecomp}) above,
then each summand $\cM_j\cong  \cB(\cH_j)$ is the weak operator closure of $\cM P_k\cM$ for some $k\in\N$, and clearly
$J$ is countable. Consider $ P[-n,n]=\sum\limits_{k\in K_n}P_k$ where $K_n$ is finite.
It follows that
\[
P[-n,n]\in \bigoplus_{j \in J_n} \cM_j,\quad\hbox{where}\quad J_n:=\{j\in J\,\mid\,P_k\in \cM_j\;\hbox{for some}\;
k\in K_n\}.\]
As the only weakly closed ideals of $\bigoplus_{j \in J_n} \cM_j$ are direct sums of $\cM_j$, and $P[-n,n]$
is not in any of these except the full sum, it follows that the weakly closed ideal generated by
$P[-n,n]$ is all of $\bigoplus_{j \in J_n} \cM_j$. It is obvious that $J_n$ is finite.
If $E\subset [-n,n]$ is measurable, then $P(E)=P(E)P[-n,n]$ hence the weakly closed ideal generated by
$P(E)$ is contained in the one generated by $P[-n,n]$, which is therefore a finite direct sum of the $\cM_j\hbox{'s}$.
 This proves
Theorem~\ref{W-CrossRepsChar}(ii) in one direction.

In the opposite direction, assume that
$\pi_T(\cM)=\bigoplus_{j \in J} \cM_j\cong \bigoplus_{j \in J}  \cB(\cH_j) $  for a
subset~$J\subeq \N$,
where the infinitesimal generator $H_j$ of the one-parameter group
$(V^j_t)_{t \in \R}$ in $\cM_j$ has compact resolvent
and every finite spectral projection of $\pi_T(V)$ generates an ideal in $\pi_T(\cM)$
consisting of a finite direct sum of the ideals $\cM_j$. Let $M\in \pi_T(\cM)$, and let $P$ be a
finite spectral projection of $\pi_T(V)$. Then $MP\in\bigoplus_{j \in J_0} \cM_j$,
 where $J_0$ is finite.
As $\cM_j = \cM Z_j$ for a central projection $Z_j$ and $Z_j\pi_T(V_t)=V^j_t$, we have
$Z_jP(E)=P_j(E)$ for any measurable set $E\subseteq\R$ and where $P(E)$ (resp. $P_j(E)$) is the spectral measure
of $\pi_T(V)$ (resp. $V_j$). As $H_j$  has compact resolvent, $(\cM_j, V_j)$ is a cross representation. This implies that, for any bounded measurable subset
$E \subeq \R$, we have
\[
P[-n,n]MP(E) = \sum_{j\in J_0}P_j([-n,n])MZ_j P_j(E) \to \sum_{j\in J_0}MZ_j P_j(E) = MP(E)
\]
as $n\to \infty$. Thus
$(\pi_T(\cM), \pi_T \circ V)$ is a cross representation.
\end{prf}

\begin{rem}
\label{W-CrossRepDecomp}
(i) We conclude from Theorem~\ref{W-CrossRepsChar} that,
for an inner $W^*$-dynamical system $(\cM,\Ad V),$ a
cross representation ${(\pi(\cM),\pi(V))}$ has a decomposition
\[ \pi=\pi_N\oplus\pi_T,\qquad  \pi_N(\cM)=\bigoplus_{k=1}^\infty \cM_k\quad\hbox{and}\quad
 \pi_T(\cM)\cong\bigoplus_{j=1}^\infty  \cB(\cH_j),
\]
where the projection onto a component $\cM_k$ of $\pi_N$ is a $V$-bounded representation,
and the projection onto a component $\cB(\cH_j)$  of $\pi_T$ converts $V$ to a one-parameter group
where the generator has compact resolvent.

(ii) If, for a concretely given von Neumann algebra $\cM\subset  \cB(\cH)$
preserved by the adjoint action $\Ad V$ defined by a unitary
one-parameter group $(V_t)_{t \in \R}$,
we have that $V_t\not\in\cM$ for some~$t$, we can of course extend the system to the
new von Neumann algebra $(V_{\R}\cup\cM)''$ to which we can apply
Theorem~\ref{W-CrossRepsChar},
and then restrict the cross representation obtained to $\cM$. Of course this leaves the
possibility that some cross representations of $\cM$ are not such
restrictions of cross representations of $(V_{\R}\cup\cM)''$.
\end{rem}

\section{Positive cross representations -- the one-parameter case}
\mlabel{RepPosCros}

In this section, we consider covariant representations of one-parameter actions
where the generator of the implementing unitary group is positive, and examine
when they are cross representations.
This will be generalized in Sections~\ref{SpecCondRev} and \ref{SpecCondCovRep} to
higher dimensions and to non-abelian groups, and in Theorem~\ref{CrossSpecParam2} we will
see that the one-parameter case fully determines the cross property for the other generalizations.

Let  $U:\R\to\U(\cH)$  be a strong operator continuous unitary one-parameter group
such that $\alpha_t:=\Ad U_t$ defines  an action
 $\alpha:\R\to \Aut(\cA)$  on  a given concrete $C^*$-algebra
$\cA\subseteq \cB(\cH)$. Let $H=H^*$ be its generator, so that $U_t=e^{itH}$. Then
$U(C^*(\R))=\{f(H)\mid f\in C_0(\sigma(H))\}$.
As we seek a crossed product for covariant representations
 $(\pi,U)\in{\rm Rep}(\alpha,\cH)$ which are positive, i.e.\ where $H \geq 0$,
we need to choose a host algebra $\cL$ that
 will produce such unitary representations of $\R$. As $U(f)=\hat{f}(H)$ for $f\in L^1(\R)$,
  it is clear that we should take for our host the quotient algebra
\[ \cL:=C^*_+(\R):= C_0([0,\infty)) \quad \mbox{  of } \quad C_0(\R) \cong C^*(\R), \]
where
the quotient map corresponds to restricting to the closed positive half-line $[0,\infty)$.
If  $U:\R\to\U(\cH)$ is a unitary one-parameter group, then the integrated representation
of $C^*(\R) \cong C_0(\R)$ factors through the quotient
$C_+^*(\R)$ if and only if $H \geq 0$, resp., $\sigma(H) \subeq [0,\infty)$.
Hence we may use the conditions in Corollary~\ref{corRCont}
to check that the covariant representation $(\pi,U)$ is a cross representation.

  Thus, the problem is just that of finding
  cross representations with respect to this host, and if a representation $(\pi,U)\in{\rm Rep}(\alpha,\cH)$
  is not a cross representation, then
we reduce the algebra $\pi(\cA)$ to  $\pi(\cA)_\cL.$ There are natural questions of
  how to obtain such cross representations from positive representations, or from cross representations
  with respect to\ the larger host algebra~$C^*(\R)$.

\begin{rem}
As $U(C^*(\R))$ is generated by a single element $L$, such as $(i\1-H)^{-1}$,
we only need to check that
 $\pi(A)L\in U_{\cL}(\cL)  \cB(\cH)$ (and similar for $A^*$) to conclude that $A\in \cM_{\cL}$.
 If ${(\pi,U)}$ is positive, we can take $L=\exp(-H)>0$. If $E_\lambda$ is an approximate identity in $\pi(C^*(\R))$,
 then we need to check that $\lim\limits_{\lambda\to\infty} E_\lambda\pi(A)L =\pi(A)L$. One natural approximate identity
 is $E_n=L^{1/n}=\exp(-H/n)$ (see \cite[Prop.~II.4.2.1]{Bla06} for more). So the condition to check is:
 \[
\lim_{n\to\infty} \exp(-H/n)\pi(A)\exp(-H) =\pi(A)\exp(-H).
 \]
 There is also the complex semigroup, i.e.\ $\{\exp(zH)|{\rm Re}(z)<0\}$ which
generates  the $C^*$-algebra $U_{\cL}(\cL)$, so it can also be used in this condition.
\end{rem}

  First we want to connect with the Borchers--Arveson Theorem  (cf.~\cite[Thm.~3.2.46]{BR02}
and Theorem~\ref{BA-thm} below),
so we begin with a positive covariant representation.
  The positivity also allows for comparison of generators of different
  implementing unitary groups, which we exploit below:

 \begin{lem}  \mlabel{lem:5.2}
Let $\cH$ be an infinite dimensional Hilbert space and
$0\leq A\leq B$ be selfadjoint operators commuting in the strong sense. Then
\[
(\1+B)^{-1}\in C^*((\1+A)^{-1}(\1+B)^{-1})\big((\1+B)^{-1}\cup(\1+A)^{-1}  \big)''\subset C^*((\1+A)^{-1}) \cB(\cH).
\]
\end{lem}

\begin{prf}
As $0\leq\1+ A\leq\1+ B$, strong commutativity implies
\begin{align*}
0\leq(\1 + A)(\1 + B)^{-1} = (\1+B)^{-1/2}(\1+A)(\1+B)^{-1/2}
\leq(\1+B)^{-1/2}(\1+B)(\1+B)^{-1/2} = \1,
\end{align*}
hence
$(\1+A)(\1+B)^{-1}\in \cB(\cH)$. Thus
\begin{eqnarray*}
(\1+B)^{-1}&=&(\1+A)^{-1}(\1+A)(\1+B)^{-1}\in(\1+A)^{-1} \cB(\cH)\qquad\hbox{and hence}\\[1mm]
(\1+B)^{-1}&\in &(\1+A)^{-1} \cB(\cH)\cap\big((\1+B)^{-1}\cup(\1+A)^{-1}  \big)''.
\end{eqnarray*}
For $R = R^*\in  \cB(\cH)$, the
intersection of any right ideal $C^*(R) \cB(\cH) = \oline{R  \cB(\cH)}$ with a
$C^*$-subalgebra $\cN\subset \cB(\cH)$ containing $R$
is just $C^*(R)\cN$ (an easy consequence of \cite[Thm~3.10.7]{Pe89}), hence
\[
C^*((\1+A)^{-1}) \cB(\cH)\cap\big((\1+B)^{-1}\cup(\1+A)^{-1}  \big)''=C^*( (\1+A)^{-1})\big((\1+B)^{-1}\cup(\1+A)^{-1}  \big)''.
\]
As $\big((\1+B)^{-1}\cup(\1+A)^{-1}  \big)''$ is commutative, it is clear that the previous right ideal is a two sided ideal.
Since it is also obvious that $(\1+B)^{-1}\in C^*( (\1+B)^{-1})\big((\1+B)^{-1}\cup(\1+A)^{-1}  \big)''$ and the intersection of ideals
is their product, we have
\[
(\1+B)^{-1}\in C^*((\1+A)^{-1}(\1+B)^{-1})\big((\1+B)^{-1}\cup(\1+A)^{-1}  \big)''
\]
 as claimed.
 \end{prf}

We first consider the case of positive $W^*$-dynamical systems
(cf. Theorem~\ref{W-CrossRepsChar},
which applies here). See Theorem~\ref{CrossRepsSame}(iii) below for the more
general $C^*$-case. In the following proposition we shall use the
Borchers--Arveson Theorem  (cf.~\cite[Thm.~3.2.46]{BR02}
and Theorem~\ref{BA-thm} below).

 \begin{prop} \mlabel{UVC}
Let $(\cM,\R, \alpha)$ be a $W^*$-dynamical system
on a von Neumann algebra
$\cM\subseteq \cB(\cH)$ and let $(U_t)_{t \in \R}$
be a positive strong operator continuous unitary
one-parameter group on~$\cH$ such that $\alpha_t=\Ad U_t$ on $\cM.$
Let $V:\R\to\cM$ be the positive one-parameter unitary  group provided by the
Borchers--Arveson Theorem. For  $\cL=C^*(\R)$ or $C^*_+(\R)$, we then have:
\begin{itemize}
\item[\rm(i)] $\cC_V:=C^*\big(V_{\cL}(\cL)\cM\big)\subseteq\cM$,
and $(\cM,V)$ is a cross representation
if and only if $\cC_V$ is a non-zero  ideal of $\cM$ contained in $\cM_c$.
\item[\rm(ii)]
$U_{\cL}(\cL)\subseteq  V_{\cL}(\cL) \cB(\cH)$ and
\[
\cC_U:=C^*\big(U_{\cL}(\cL)\cM\big)
\subseteq C^*\big(V_{\cL}(\cL)(\cM\cup U_{\cL}(\cL))''\big)\supseteq\cC_V.
\]
\end{itemize}
\end{prop}

\begin{prf}
(i) As $V_\R\subset\cM$, we have $V_{\cL}(\cL)\subset\cM$ and
hence $V_{\cL}(\cL)\cM\subseteq\cM\supseteq \cC_V$. Now
$(\cM,V)$ is a cross representation if and only if $\cM_{\cL}=\cM$ and this implies that $\cC_V\subset\cM$ is a two-sided ideal of $\cM$
by Proposition~\ref{propALchar}(i) and Corollary~\ref{corAisAL}(i).
That it is in $\cM_c$ follows from Corollary~\ref{corAisAL}(i) as in this case
$\cM V_{\cL}(\cL)\subseteq V_{\cL}(\cL) \cB(\cH)$ hence $\cC_V\subseteq V_{\cL}(\cL) \cB(\cH)$ and $t\mapsto V_tA$ for $A\in V_{\cL}(\cL) \cB(\cH)$
is norm continuous. As it is a $C^*$-algebra, also right multiplication of $\cC_V$ with $U_t$  produces functions norm continuous in $t$, so
$\cC_V\subseteq \cM_c$.

Conversely, let $\cC_V$ be
a non-zero ideal of $\cM$ contained in $\cM_c$.
Then $\cM V_{\cL}(\cL)\subset\cC_V$ and $\alpha$ acts strongly continuously on   $\cC_V$ . Thus, for all $A\in \cM$ and
$L\in V_{\cL}(\cL)$, we have that
\[
\lim_{t\to 0}\big\|(\alpha_t(A)-A)L\big\|\leq\lim_{t\to 0}\|\alpha_t(AL)-AL\|+\lim_{t\to 0}\|\alpha_t(A)(\alpha_t(L)-L)\|=0
\]
as  $\alpha_t=\Ad V_t$ also acts strongly continuously on $V_{\cL}(\cL)$.
Thus $\lim\limits_{t\to 0} \alpha_t(A)L=AL$ for all $L\in V_{\cL}(\cL)$, and likewise for $A^*$, hence
by Proposition~\ref{contSpace} we have that $A\in\cM_{\cL}$. As this holds for all
$A\in \cM$, it follows that $(\cM,V)$ is a cross representation.

(ii) Let $U_t=e^{itB}$ and $V_t=e^{itA}$ for positive operators $A,\,B$ on $\cH$. As $\alpha_t=\Ad U_t=\Ad V_t$ on $\cM,$
and $V_t\in\cM$ we have that $(\Ad  U_t)(V_s)=\alpha_t(V_s)=(\Ad  V_t)(V_s) =V_s$ hence
$U_t V_s = V_s U_t$ for all $t,\,s\in\R$.
Therefore $A$ and $B$ commute in the strong sense, i.e.\ their spectral projections
commute (cf.~\cite[Thm.~VIII.13]{RS80}).
Moreover as $V:\R\to\cM$ is the positive  unitary  group constructed
in the Borchers--Arveson Theorem,
we have that $B\geq A\geq 0$ (cf.~Equation (3.2) and the
subsequent proposition on \cite[p.235]{Arv74}). Thus, by Lemma~\ref{lem:5.2},
we have
\[
U_{\cL}(\cL)=C^*((\1+B)^{-1})\subset C^*((\1+A)^{-1}) \cB(\cH)=V_{\cL}(\cL) \cB(\cH).
\]
By the same lemma, we also have $U_{\cL}(\cL)\subseteq C^*(V_{\cL}(\cL)U_{\cL}(\cL))\big(U_{\cL}(\cL)\cup V_{\cL}(\cL)  \big)''$ so
 \[
U_{\cL}(\cL)\cM \subseteq C^*(V_{\cL}(\cL)U_{\cL}(\cL))\big(U_{\cL}(\cL)\cup V_{\cL}(\cL)  \big)''\cM
\subseteq V_{\cL}(\cL)(\cM\cup U_{\cL}(\cL))''
 \]
 by $V_{\cL}(\cL)\subset\cM\ni\1$, and commutativity of $\big(U_{\cL}(\cL)\cup V_{\cL}(\cL)  \big)''$.
 \end{prf}
In general we do not have $\cC_U\subseteq\cC_V$ because $\cC_V\subseteq\cM$, and if $U$ is outer then
 $\cC_U\supset U_{\cL}(\cL)$ will not be in $\cM$.
 By part (i), for nontrivial actions, the
 Borchers--Arveson positive representation is not a cross representation for simple
 von Neumann algebras (e.g.\ finite factors or countably decomposable type III factors
\cite[Cor.~6.8.4, 6.8.5]{KR86}). This is not as serious as it looks, as in general $\cM=\pi(\cA)''$ and $(\pi,V)$ can be a cross representation,
 even when $(\cM,V)$ is not (cf. the trivial Example~\ref{CrossNoExt}).

 \begin{thm} \mlabel{CrossRepsSame}
Let $(\cM, \R,\alpha)$ be a $W^*$-dynamical system on a von Neumann algebra
$\cM\subseteq \cB(\cH)$ and let $U:\R\to\U(\cH)$
be a positive  unitary one-parameter group
such that $\alpha_t=\Ad U_t$ on $\cM.$
Let $V:\R\to\cM$ be the positive one-parameter unitary  group provided by the
Borchers--Arveson Theorem.
Then, for $\cL=C^*(\R)$ or $C^*_+(\R)$, the following assertions hold:
\begin{itemize}
\item[\rm(i)] We have
\begin{eqnarray*}
\cM^V_\cL &:=& \big\{ M \in \cM \mid
BV_\cL(\cL) \subeq V_\cL(\cL) \cB(\cH)\quad\hbox{for}\quad B\in\{M,M^*\}\big\}\\[1mm]
&\subseteq& \cM^U_\cL := \big\{ M \in \cM \mid
BU_\cL(\cL) \subeq U_\cL(\cL) \cB(\cH)\quad\hbox{for}\quad B\in\{M,M^*\}\big\}.
\end{eqnarray*}
\item[\rm(ii)]
If $(\cM,V)$ is a cross representation, then $(\cM,U)$ is a cross representation.
\item[\rm(iii)]
Let  $(\cA, G, \alpha)$
be a $C^*$-action and let $\pi \: \cA \to  \cB(\cH)$ be a
representation such that
 $\pi(\cA)\subseteq\cM$, and $\pi\circ\alpha_t=\Ad U_t\circ\pi=\Ad V_t\circ\pi$,
 with $U,\,V,\,\cM$ as above. If
 $(\pi,V)$ is a cross representation, then $(\pi,U)$ is a cross representation.
\end{itemize}
\end{thm}

\begin{prf}
(i) 
We have the situation in Proposition~\ref{UVC}, hence we may take
 $U_t=e^{itB}$ and $V_t=e^{itA}$ for positive operators $A,\,B$ on $\cH$,
where $A$ and $B$ commute in the strong sense, and  $B\geq A\geq 0$. Moreover $C:=B-A\geq 0$ is the generator of
the homomorphism $t\mapsto W_t:=U_tV_{-t}$. As the conjugation action of $W_t$ on $\cM$ is trivial, we have $W_t\in\cM'$.

Assume for an $M \in \cM$ that
$MV_\cL(\cL) \subeq V_\cL(\cL) \cB(\cH)$. Then, using Proposition~\ref{UVC}, we have
$MU_\cL(\cL) \subeq MV_{\cL}(\cL) \cB(\cH)\subeq V_\cL(\cL) \cB(\cH)$. Thus $t\mapsto V_tMU_\cL(\cL)$ is continuous. 
Moreover, $t\mapsto U_\cL(L)e^{itC}$ is
continuous for
$L \in \cL$, as
\[
t\mapsto (\1+B)^{-1}e^{itC}=(\1+A+C)^{-1}e^{itC}= [(\1+A+C)^{-1}(\1+C)](\1+C)^{-1}e^{itC}
\]
 is continuous by $(\1+A+C)^{-1}(\1+C)\in \cB(\cH)$.
Thus, by  $e^{itC}\in\cM'$, we have
\[
U_tMU_\cL(L)=e^{itC}V_tMU_\cL(\cL)= V_tMU_\cL(L)e^{itC}
=(V_tMU_\cL(L))\cdot( U_\cL(L)e^{itC}),
\]
hence the maps $t\mapsto U_tMU_\cL(L)$ are continuous for all $L\in\cL$.
Hence by \cite[Lemma~A.1(ii)b]{GrN14}, we have that
$MU_\cL(L) \subeq U_\cL(\cL) \cB(\cH)$. By applying this also to $M^*$ we see that $\cM^V_\cL\subseteq\cM^U_\cL$.

(ii) If $(\cM,V)$ is a cross representation,
then $\cM=\cM^V_\cL$ and so by the previous part we get $\cM=\cM^U_\cL$,
i.e.\ $(\cM,U)$ is a cross representation.

(iii)
If $(\pi,V)$ is a cross representation, then $\pi(\cA)\subseteq \cM^V_{\cL},$ so
by (i) we get $\pi(\cA)\subseteq \cM^U_{\cL},$ i.e.
 $(\pi,U)$ is a cross representation.
\end{prf}
Case (ii), i.e. $\cM=\cM^V_\cL$, means that as we have an inner cross representation
of a von Neumann algebra, hence Theorem~\ref{W-CrossRepsChar} applies and so this representation has the
form in Remark~\ref{W-CrossRepDecomp}(i).

Crossed products for the usual case were first fully defined
 in the paper by Doplicher, Kastler and Robinson \cite{DKR66}. Later in that paper an
 ideal of the crossed product was factored out, to obtain an algebra
which is a crossed product host (in our terminology) which produced only covariant representations
satisfying a desired spectral condition. The question arises whether this can also be done in
our context, starting with a crossed product host, e.g.\ where the action is discontinuous, but cross representations exist.
By Theorem~\ref{thm:5.1}(c) we know that a factor algebra of a crossed product host is again a crossed product host.
The main question is then whether there is a suitable ideal to factor out, so that the resulting
factor (which is a crossed product host), will have only  positive covariant representations.

Let us start with a crossed product host $(\cC, \eta_\cA, \eta_\cL)$
for $(\alpha,\cL)$, where $\cL=C^*(\R)$.
Our aim is to find an ideal in $\cC$ such that the factor algebra is
 a crossed product host producing only positive covariant representations.
Thus we would like to have as our host the $C^*$-algebra
\[ C^*_+(\R):= C_0([0,\infty))
= C^*\{ \hat f \res_{[0,\infty)}\,\mid\, f \in L^1(\R) \}
\cong C_0(\R)/\cL_-,\]
where
\[ \cL_-:=\{ h \in C_0(\R)\,\mid\,  h\res_{[0,\infty)} = 0\}
\cong C_0((-\infty,0)) \trile \cL = C_0(\R).\]
In $\cL$ we also consider the ideal
\[ \cL_+:=\{ h \in C_0(\R) \,\mid\,  h\res_{(-\infty,0]} = 0\}
\cong C_0((0,\infty)) \trile \cL = C_0(\R)\]
and note that $\cL/(\cL_+ + \cL_-) \cong \C$.

If $\eta_\cL( \cL_-)=0$, then $\cC$ is already a crossed product host for positive covariant representations,
so let us assume that $\eta_\cL(\cL_-)\not=\{0\}$. As
$\cC=C^*\big(\eta_\cA(\cA) \eta_\cL(\cL)\big)
=\eta_\cL(\cL) \cC\eta_\cL(\cL)$ follows from the cross property, it is natural to consider
the $C^*$-subalgebra $\cC_-:=C^*\big(\eta_\cA(\cA) \eta_\cL( \cL_-)\big)$,
but in general it need not be an ideal.

\begin{prop} \mlabel{Cideal}
Let  $\alpha:\R\to \Aut(\cA)$ be given, together with
a crossed product host $(\cC, \eta_\cA, \eta_\cL)$
for $(\alpha,\cL)$ where $\cL=C^*(\R)$.
 If
\[ \eta_\cL(\cL_-)\eta_\cA(\cA)\subseteq\br \eta_\cA(\cA) \eta_\cL(\cL_-).\]
then $\cC_-:=C^*\big(\eta_\cA(\cA) \eta_\cL( \cL_-)\big)$ is a closed two-sided ideal of $\cC$.

Moreover, the quotient algebra  $\cC_{(+)}:=\cC/\cC_-$ is also
a crossed product host. Let $\Phi:\cC\to\cC/\cC_-$ be the quotient
map, then  ${\big(\cC/\cC_-,\tilde\Phi\circ\eta_\cA,\tilde\Phi\circ\eta_\cL\big)}$
is a crossed product host with respect to\ $\cL=C^*(\R)$.
As $\tilde\Phi\circ\eta_\cL(\cL_-)=\{0\}$,
$\cC_{(+)}=\cC/\cC_-$ is a crossed product host with respect to\ the host  $C^*_+(\R)$, the set $\Rep(\alpha,{\cal H})_{\eta_\times}$ consists of positive
representations of $(\cA, \R, \alpha)$, i.e.,
$C^*_+(\R)$-representations.
\end{prop}

\begin{prf} If $\eta_\cL(\cL_-)\eta_\cA(\cA)\subseteq\br\eta_\cA(\cA) \eta_\cL(\cL_-).$, then
all monomials in elements of
\[ \eta_\cA(\cA) \eta_\cL(\cL_-)\cup \eta_\cL(\cL_-)\eta_\cA(\cA) \]
can be written as elements of $\br\eta_\cA(\cA) \eta_\cL(\cL_-).$ hence
$\lbr\eta_\cA(\cA) \eta_\cL(\cL_-)\rbr=\cC_-$. Likewise we have
$\cC=\lbr\eta_\cA(\cA) \eta_\cL(\cL)\rbr$ and so
\begin{eqnarray*}
\eta_\cA(\cA) \eta_\cL(\cL)\cdot \eta_\cA(\cA) \eta_\cL(\cL_-)&\subseteq&
\eta_\cA(\cA) \eta_\cL(\cL)\br \eta_\cL(\cL_-)\eta_\cA(\cA).\\[1mm]
&\subseteq& \br\eta_\cA(\cA)  \eta_\cL(\cL_-)\eta_\cA(\cA).
\subeq \br\eta_\cA(\cA)\eta_\cL(\cL_-).\subeq \cC_-\\[1mm]
 \eta_\cA(\cA) \eta_\cL(\cL_-)\cdot\eta_\cA(\cA) \eta_\cL(\cL)&\subseteq&
 \br \eta_\cA(\cA) \eta_\cL(\cL_-) \eta_\cL(\cL).\subeq
 \br\eta_\cA(\cA) \eta_\cL(\cL_-).\subeq \cC_-
\end{eqnarray*}
leads to $\cC\cC_-\subeq\cC_-$ and  $\cC_-\cC\subeq\cC_-$. Thus $\cC_-$ is a closed two sided ideal of
$\cC$. That $\cC_{(+)}:=\cC/\cC_-$ is
a crossed product host follows from Theorem~\ref{thm:5.1}(c), so we only need  to prove that
$\tilde\Phi\circ\eta_\cL(\cL_-)=0$.

For $L\in\cL_-$ and $C\in\cC$, we have
\[
\big(\tilde\Phi\circ\eta_\cL(L)\big)\Phi(C)=\Phi\big(\eta_\cL(L)C\big)=0
\]
as $\eta_\cL(L)C\in\eta_\cL(\cL_-)\cC=\eta_\cL(\cL_-)\lbr
 \eta_\cL(\cL)\eta_\cA(\cA)\rbr \subseteq
\lbr \eta_\cL(\cL_-)\eta_\cA(\cA)\rbr=\cC_-=\ker\Phi$.
\end{prf}

\begin{rem}
If  the crossed product host  $(\cC, \eta_\cA, \eta_\cL)$ is faithfully represented in a
cross representation
$(\pi, U)$ of $(\cA, G, \alpha)$ as in Theorem~\ref{thm:5.1}, then the
condition of Proposition~\ref{Cideal} that
$\eta_\cL(\cL_-)\eta_\cA(\cA)\subseteq\eta_\cA(\cA) \eta_\cL(\cL_-)$ requires that it is also a cross representation
with respect to\ $\cL_-$. This condition can be interpreted as a consistency condition for a quantum constraint condition in the following sense
(cf.\ \cite{Gr06}, \cite{GrL00}).
Note that $\cH_-:=P(-\infty,0)\cH=U_\cL(\cL_-)\cH$. Hence,
if $B = A, A^*$ both satisfy
$\pi(B)U_\cL(\cL_-)\subseteq U_\cL(\cL_-)\pi(\cA)$,
then $\pi(A)$ will preserve both the space $\cH_-$ and its orthogonal
complement $\cH_+:=\cH_-^\perp=P[0,\infty)\cH$. Then the quotient
map consists of restricting $\pi(A)$ to $\cH_+$
and this can only be done if $\pi(A)$ preserves these spaces. The quantum constraint condition
is that of putting the space $\cH_-$ to zero.
\end{rem}

\begin{ex} \mlabel{ex:4.7}
(a)  An example where the condition  from Proposition~\ref{Cideal}
holds is easily obtained.
Let  $H=H^*$ be an unbounded operator on $\ell^2$ such that  $(i\1-H)^{-1}\in \cK(\ell^2)$,
$U_t:=e^{itH}$ and $\alpha_t:=\Ad U_t$ for all $t\in\R$, and assume the host $\cL=C^*(\R)$.
In addition, we assume that the negative part of the spectrum of $H$
is unbounded.
Let $\cH_-:=P(-\infty,0)\cH=U_\cL(\cL_-)\ell^2$ as above, and let
 $\cA:= \cB(\cH_-)+\cD\subset \cB(\cH),$    where $\cD\subset \cB(\cH)$ is the unital
non-degenerate $C^*$-algebra
 $\cD:=\C\1+C^*((i\1-H_+)^{-1})$ where $H_+:=P[0,\infty)H.$
Then the given representation on $\ell^2$ is a  a cross representation
(cf. \cite[Prop.~5.10]{GrN14}),
and $\alpha$ is not strongly continuous
(cf. \cite[Th.~6.1, Ex.~5.11]{GrN14}).
Obviously $\cD$ commutes with both $U_\cL(\cL)$ and $U_\cL(\cL_-)$ hence we only need to verify that
$U_\cL(\cL_-)\cB(\cH_-)\subseteq\br \cB(\cH_-) U_\cL(\cL_-).$  to conclude that
$U_\cL(\cL_-)\cA\subseteq \br\cA U_\cL(\cL_-).$. Now as we started with a cross representation, we have
$\br U_\cL(\cL)\cA.=\br \cA U_\cL(\cL).$
(cf.\ the equivalent condition in Definition~\ref{crossrepDef}).
As $\cB(\cH_-)\subset\cA$, we have
\[ \br U_\cL(\cL_-)\cB(\cH_-).
\subeq \br U_\cL(\cL)\cA.
\subseteq\br \cA U_\cL(\cL).\,.\]
As  $U_\cL(\cL_-)$ consists of compact operators,
\[
\br U_\cL(\cL_-)\cB(\cH_-).\subseteq \cK(\cH_-) =\br \cB(\cH_-) U_\cL(\cL_-).
\]
Thus $U_\cL(\cL_-)\cA\subseteq \br\cA U_\cL(\cL_-).$ so we have satisfied the condition of
Proposition~\ref{Cideal}.
Now $\cC=C^*\big(\eta_\cA(\cA) \eta_\cL(\cL)\big)
=\cK(\cH_-)+C^*((i\1-H_+)^{-1})$ and $\cC_-=C^*\big(\cA U_\cL( \cL_-)\big)
=\cK(\cH_-),$ thus $\cC/\cC_-\cong C^*((i\1-H_+)^{-1})$.

(b) For an example where the condition from Proposition~\ref{Cideal}
fails, consider the translation action of $\R$ on
$\cA=C^*\{\delta_t\mid t\in\R\}\subset C_b(\R)$ where $\delta_t(x):=e^{itx}$ acts by multiplication on $L^2(\R)$.
The translation action is
implemented by the unitaries $U_t:=e^{itP}$ with $P=i{d\over dx}$ the usual momentum operator. 
We first show that $(\cA,U)$ is a cross representation. Now $U_\cL(\cL)=\{f(P)\mid f\in C_0(\R)\}$ and
\[
\delta_t f(P)\delta_{-t}=f(P-t\1)=f_t(P)\quad \hbox{for}\quad f\in C_0(\R),
\quad\hbox{where}\quad
f_t(x):=f(x-t).
\]
Thus, as $f_t\in C_0(\R)$, we have
 $\delta_tU_\cL(\cL)\subseteq U_\cL(\cL)\cA$ and so $(\cA,U)$ is a cross representation for the
host $\cL=C^*(\R)$.
On the other hand, if the support of $f$ is in $(-\infty,0)$,
then this need not be true for the
translated function $f_t$,
i.e.\ we have $\delta_t f(P)=f_t(P)\delta_{t}$ and so by choosing $t$
so that the support of $f_t$ will not be in $(-\infty,0)$, we can see that the condition
of Proposition~\ref{Cideal} is violated: the algebra $\cA \cL_-$ annihilates the
Fourier transform of $L^2(\R_+)$, but $\cL_- \cA$ does not.
As $\cA$ is commutative, we already know
from the Borchers--Arveson Theorem that there are no
non-trivial positive covariant representations.
\end{ex}

A common situation for physics is how to build a positive covariant representation
from one which is not (especially if there is no positive subrepresentation).
One can treat this as a constraint situation by selection of those
elements of the algebra which are compatible with restriction of the generator to its positive part.

\begin{prop} \mlabel{CrossConstrain}
Let $(\pi, U)$ be a covariant representation of $(\cA,\R,\alpha)$.
Let $U_t=e^{itH}$ and let
 $P$ be the spectral measure of $H$. Define
\[  \cO_+^\pi:=P[0,\infty)'\cap\cA=\{A\in\cA\,\mid\,
A P[0,\infty) = P[0,\infty) A\}. \]
Then
\begin{itemize}
\item[\rm(i)]
$\cO_+^\pi$ is an $\alpha$-invariant $C^*$-subalgebra of~$\cA$ preserving $\cH^+:=P[0,\infty)\cH=(U(\cL_-)\cH)^\perp$.
Moreover the restriction $(\pi^+,U^+)$ of $(\pi,U)$ on $\cO_+^\pi$ to $\cH^+$,
i.e.\ $\pi^+(A):=\pi(A)\restriction\cH^+$, $A\in \cO_+^\pi$ and
$U^+_t:=U_t\restriction\cH^+$,
 is a positive covariant representation of  $(\cO_+^\pi, \R, \alpha)$.
 \item[\rm(ii)]
 If $(\pi, U)$ is a cross representation, then  $(\pi^+,U^+)$ is also a cross representation.
 \item[\rm(iii)] If  $(\pi, U)$ is {\it not} a cross representation, let
 \[
\cC_+^\pi:=\big\{A\in\cA\,\mid\, \pi(B)U_\cL(\cL)P[0,\infty)\subseteq U_\cL(\cL)P[0,\infty) \cB(\cH)\quad\hbox{for}\quad B\in\{A,\; A^*\}\big\}.
\]
Then $\cC_+^\pi\subseteq \cO_+^\pi$ is $\alpha\hbox{-invariant}$
and $(\pi^+,U^+)$ restricted to  $\cC_+^\pi$ is a cross representation.
Moreover $\cC_+^\pi$ is the maximal $C^*$-subalgebra of
$\cO_+^\pi$ on which $(\pi^+,U^+)$ is cross.
\end{itemize}
\end{prop}

\begin{prf}
%
(i) From the definition, it is clear that $\cO_+^\pi$ is a $C^*$-algebra. To see that it is
 $\alpha$-invariant, just note that $(\pi, U)$ is a covariant representation where
 $U_t$ commutes with $P[0,\infty)$.
%
 For $\cH^+ = P[0,\infty)\cH$,
$t\mapsto U^+_t:=U_t\restriction \cH^+$ is a positive unitary one-parameter group,
and as $\pi(\cO_+^\pi)$
 restricts to $\cH^+$, it  follows from
\[ \pi(\alpha_t(A))\cH^+ = U_t\pi(A)U_{-t}\cH^+= U^+_t\pi(A)U^+_{-t}\cH^+
\quad \mbox{  for all } A\in\cO_+^\pi  \]
that $(\pi^+,U^+)$ is a positive covariant representation of  $(\cO_+^\pi, \R, \alpha).$

(ii) If  $(\pi, U)$ is a cross representation, then by Definition~\ref{crossrepDef} we have
$\pi(\cA) U_{\cL}(\cL) \subeq \br U_{\cL}(\cL) \pi(\cA).$. Note that $\pi_+(A)=
P[0,\infty)\pi(A)\restriction\cH^+$ for $A\in \cO_+^\pi$.  Thus
\begin{eqnarray*}
\pi^+(\cO_+^\pi)U^+_{\cL}(\cL)
&=&P[0,\infty)\pi(\cO_+^\pi)U_{\cL}(\cL)\restriction\cH^+
\subseteq  P[0,\infty) \br U_{\cL}(\cL) \pi(\cA).\restriction\cH^+\\[1mm]
&=&\br U^+_{\cL}(\cL)P[0,\infty) \pi(\cA).\restriction\cH^+
\subseteq U^+_{\cL}(\cL) \cB(\cH^+).
\end{eqnarray*}
Thus  $(\pi^+,U^+)$ is also a cross representation.

(iii) Note that
\[U_\cL(\cL)P[0,\infty)={\{f(H)\mid f\in C_0([0,\infty))\}},\]
 hence, if $A\in\cC_+^\pi$ and $B \in \{A,A^*\},$ then
 \[
 \pi(B)P[0,\infty)\cH=\pi(B)U_\cL(\cL)P[0,\infty)\cH\subseteq U_\cL(\cL)P[0,\infty)\cB(\cH)\cH=P[0,\infty)\cH,
 \]
i.e.\ $\pi(A)$ and its adjoint preserve the space $\cH^+=P[0,\infty)\cH$,
hence commute with $P[0,\infty)$, and so $\cC_+^\pi\subseteq \cO_+^\pi$.
 From the definition, it is clear that $\cC_+^\pi$ is a $C^*$-algebra. To see that it is
 $\alpha$-invariant, just note that $U_\cL(\cL)P[0,\infty)$ commutes with $U_t$, hence
 \begin{align*}
 \pi(\alpha_t(B))U_\cL(\cL)P[0,\infty
&)=U_t\pi(B)U_\cL(\cL)P[0,\infty)U_{-t}\subseteq U_tU_\cL(\cL)P[0,\infty)\cB(\cH)U_{-t}\\
& =U_\cL(\cL)P[0,\infty)\cB(\cH).
 \end{align*}
 Thus
$(\pi^+,U^+)$ restricts on  $\cC_+^\pi$ to a positive covariant representation.
That it is a cross representation follows from the restriction of
 the defining condition for $\cC_+^\pi$ to  $\cH^+$. To prove maximality, note that the
 cross condition for $(\pi^+,U^+)$ is for $B=A,\,A^*\in \cO_+^\pi$
 \begin{equation}
   \label{eq:*}
\pi^+(B) U^+_\cL(\cL)\subseteq U^+_\cL(\cL)\cB(\cH^+).
 \end{equation}
 Recalling that  $\pi^+(A)=P[0,\infty)\pi(A)\restriction\cH^+$ this condition
\eqref{eq:*} becomes
  \[
\pi^+(B) U_\cL(\cL)P[0,\infty)\restriction\cH^+\subseteq U_\cL(\cL)P[0,\infty)\cB(\cH)\restriction\cH^+
 \]
 which is equivalent to the defining condition for  $\cC_+^\pi$ if we keep in mind that the
 left hand side of that condition vanishes
on $(\cH^+)^\perp$ hence restriction of both sides
 of that condition to
 $\cH^+$ does not change $\cC_+^\pi$. We conclude that $A\in\cC_+^\pi$ which proves the maximality
 claimed for $\cC_+^\pi$.
\end{prf}
\begin{rem}
 Even if  $(\pi,U)$ is a covariant cross representation of
 $(\cA, \R, \alpha)$ with respect to\ $\cL=C^*(\R)$, it is possible that $\cO_+^\pi=\{0\}.$
 A good example is the translation action of $\R$ on $\cA=C_0(\R)$ with  $(\pi,U)$ the
usual covariant representation
 on $L^2(\R)$ (see last sentence in Example~\ref{ex:4.7}(b)).
\end{rem}

 \begin{ex}
 Let $\cA= \cB(\cH)$, where $\cH$ is an infinite dimensional Hilbert space, $H=H^*$ be an unbounded
 selfadjoint operator on $\cH$, and
$U_t:=e^{itH}$ and $\alpha_t:=\Ad U_t$ for all $t\in\R$, and let $\pi$ be the identical representation.
Then $\alpha$ is not strongly continuous (cf.~\cite[Prop.~5.10]{GrN14}).
We will also assume that the spectrum of $H$
has negative parts, $\sigma(H)\cap(-\infty,0)\not=\emptyset$ and choose the host
$\cL=C^*(\R)$. Let
 $P$ be the spectral measure of $H$. Then $\cH^+=P[0,\infty)\cH$,
and the $C^*$-subalgebra of $\cA$
 consisting of all elements which preserve
 $\cH^+$ and its orthogonal complement $\cH^-$,
is the commutant of $P[0,\infty)$ i.e.\
 \[
 P[0,\infty)'= P[0,\infty) \cB(\cH)P[0,\infty)\oplus P(-\infty,0) \cB(\cH)P(-\infty,0)\cong
  \cB(\cH^+)\oplus \cB(\cH^-)\,. \]
The restriction of this to  $\cH^+$ is just the first component,
and if $(i\1-H)^{-1}P[0,\infty)\not\in\cK(\cH^+)$,
then $\big( \cB(\cH^+),U^+\big)$ is not a cross representation
(cf.~\cite[Ex.~5.11]{GrN14} for separable Hilbert spaces).
By Proposition~\ref{CrossConstrain}(iii), $(\pi^+,U^+)$ is always a cross representation on $\cC_+^\pi$,
hence the containment $\pi(\cC_+^\pi)\subeq P[0,\infty)'$ can be proper. As $\cK(\cH^+)\cup U_\R'\subseteq\pi(\cC_+^\pi)$,
we also know that $\cC_+^\pi$ is not trivial.
It is easy to adapt Corollary~\ref{corRCont} to this example to provide useful conditions such as:
\[ \cC_+^\pi=\big\{A \in \cA\,\mid
(\forall s \in \R_+)(\forall B \in \{A,A^*\})\  \lim_{t \to \infty} P[0,t]BP[0,s] = BP[0,s]\big\}.
\]
\end{ex}

\section{Covariant representations of actions on topological groups}
\label{CRATG}

There are many instances of one-parameter groups which act on topological groups, and this
provides an interesting subclass of covariant systems. Positivity conditions for the
one-parameter group have strong structural implications as seen e.g. in \cite[Ex.~4.16]{BGN17}.

Let $G$ be a topological group and
$\alpha \: \R \to \Aut(G)$ be a homomorphism defining a continuous action of
$\R$ on $G$. We then form the topological semidirect product group
$G^\sharp \cong G \rtimes_\alpha \R$
(cf.\ \cite{Ne14}). Continuous unitary representations of $G^\sharp$ are
pairs $(\pi, U)$, where $\pi \: G \to \U(\cH)$ is a continuous unitary representation
and $(U_t)_{t \in \R}$ is a unitary one-parameter group satisfying
\[ U_t \pi(g) U_{-t} = \pi(\alpha_t(g)) \quad \mbox{ for } \quad g \in G, t \in \R.\]
Let $\cA := C^*(G_d)$ denote the $C^*$-algebra of the discrete group underlying $G$.
Then $\alpha$ extends to a $C^*$-action
$(\cA, \R,\alpha)$. We are interested in situations, where the positive
covariant representations of $G^\sharp$
can be described in terms of crossed product hosts for
$(\cA,\R,\alpha)$ with respect to the positive  quotient
$C^*(\R)/C_0((-\infty,0)) \cong C_0([0,\infty)) =: C^*_+(\R)$.

In this section we will analyze the following interesting situation for this context.
 In a positive cross representation $(\pi, U)$ of the action
$\alpha \: \R \to \Aut(G)$,
 the positive unitary one-parameter group $(U_t)_{t \in \R}$
can ``regularize'' the representation $\pi$ of
the group~$G$. For instance, in this representation the group algebra $U(C^*(\R))$ of the one-parameter group times
the discrete group algebra $\pi(C^*(G_d))$
of the group can be a crossed product host which only allows
continuous group representations.
(This happens if the maps $g\mapsto\pi(g)U(B)$ are continuous for all $B\in C^*(\R)$.)
This will produce another method for construction of host algebras for
some infinite dimensional Lie groups.
This then leads to the study of smoothing operators
(cf.\ Theorem~\ref{thm:smoothop} below, or \cite{NSZ17})
which provide a method for constructing host algebras for some infinite dimensional Lie groups. 
In particular, the structure of these host algebras makes it
easier to construct crossed product hosts for the groups involved.

In this context it is interesting to recall that,
in \cite{Ta67} one finds an example of a non-type I $C^*$-algebra
with automorphism group $G$ such that all covariant representations
of $(\cA, G, \alpha)$ are type I. So, even without
spectral conditions, forming crossed products
may have a regularizing effect. This is consistent with our observations
for crossed products with $\cL = C^*(S,C)$, studied in
Subsection~\ref{OlshHost}.

\subsection{Regularization by one-parameter subgroups}
\label{RegPosGp}

Here we want to examine the  phenomenon mentioned above. We start with an example
to illuminate the issue.

\begin{ex}
\label{OscRegHeis}
Let $\cA = \oline{\Delta(X,\sigma)}$ be the Weyl algebra of the symplectic space $(X,\sigma)$
consisting of a  complex pre-Hilbert space $(X,\langle\cdot,\cdot\rangle)$ and the
symplectic form $\sigma(x,y):={\rm Im}\langle x,y\rangle$.
It is the unique unital $C^*$-algebra
with generating unitaries  $(\delta_z)_{z \in X}$ satisfying
\[ \delta_z^*=\delta_{-z} \quad \mbox{ and } \quad
\delta_z\delta_w=e^{-i\sigma(z,w)/2}\delta_{z+w}
\quad \mbox{ for } \quad z,w \in X.\]
Define
a $\T$-action by $\alpha_t (\delta_z) := \delta_{tz}$ for $t \in \T$, $z\in X$.
Then the corresponding Fock representation $(\pi,U,\cF)$ of $(\cA,\T, \alpha)$
(see also Example~\ref{NoGap}) is a cross representation
for $\alpha$  with respect to $\cL=C^*(\T)\cong C_0(\Z)$
(cf.~\cite[Ex.~6.11]{GrN14}). Thus
\begin{equation}
  \label{eq:c-def5}
\cC := C^*(\pi(\cA)U(\cL)) \subeq U(\cL)\cB(\cF)\cap  \cB(\cF)U(\cL)\,.
\end{equation}

As $\cA$ contains a copy of
\[ G = \Heis(X,\sigma) = \T \times X \quad \mbox{ with multiplication } \quad
(a,z)(b,w) = (ab e^{-i\sigma(z,w)/2}, z + w) \]
through the identification
$(a,x) \mapsto a\delta_x$ for $a\in\T, x\in X,$
the action $\eta_\cA \: \cA \to M(\cC)$ also produces a unitary multiplier
action of $G$ on $\cC$.
We will prove that this action is norm continuous, hence that all covariant representations
produced by $\cC$ for $(\cA, \T, \alpha),$ must be regular for $\cA$. In this sense the one
parameter group acting on $G$ regularizes the representations of $G$.

The generator of $U_t$ is the number operator, hence
the $U$-eigenspaces are the $n$-particle  subspaces $\cF_n \subeq \cF$.
Let $P_n \in  \cB(\cF)$ denote the projection onto $\cF_n$.
If we write $\pi(\delta_x)=\exp(i\varphi(x))$, then we get
as in the proof of \cite[Thm.~X.41]{RS75}
\begin{eqnarray*}
\|\varphi(x)P_n\|&\leq& \left(\|a(x)P_n\|+\|a(x)^*P_n\|    \right)\big/\sqrt{2}
\leq 2^{1/2}(n+1)^{1/2}\|x\|, \\
\mbox{ hence}\qquad
\|\varphi(x)^kP_n\|&\leq& 2^{k/2}((n+k)!)^{1/2}\|x\|^k.
\end{eqnarray*}
As
\[ \sum_{k=0}^\infty\|\varphi(x)^kP_n\|\frac{s^k}{k!}\leq
\sum_{k=0}^\infty2^{k/2}((n+k)!)^{1/2}\|x\|^k
\frac{s^k}{k!}<\infty, \]
the series $\sum_{k=0}^\infty\varphi(sx)^k/k!$ on
the space $\cF_n$ is norm convergent,
and converges to an analytic function in the variable $s\geq 0$,
and the limit is  $\pi(\delta_{sx}) P_n$. In fact, from the last sum we see
\[
\big\|(\pi(\delta_x)-\1)P_n\big\|\leq\sum_{k=1}^\infty\|\varphi(x)^kP_n\|\frac{1}{k!}\leq\sum_{k=1}^\infty2^{k/2}((n+k)!)^{1/2}\|x\|^k
\frac{1}{k!}
\]
and it is clear that this goes to zero as $\|x\|\to 0$. Therefore the maps
$X \to  \cB(\cF), x \mapsto \pi(\delta_x) P_n$
are norm-continuous (and in fact analytic).
As $U(\cL)=C^*\big(\{P_n\mid n\geq 0\}\big)$, we get that the maps
\begin{equation}
\label{Lcont}
 x \mapsto \pi(\delta_x)L,\;\; L\in U(\cL)
\end{equation}
are norm continuous. By the cross representation condition, this
implies that the crossed product host
\[ \cC := C^*(\pi(\cA)U(\cL)) \subeq  \cB(\cF) \]
has the property that the  multiplier action of $G = \Heis(X,\sigma)$ on
$\cC$ through $(a,x) \mapsto a\eta_\cA(\delta_x)\in M(\cC)$
is  continuous. 
We conclude that $\cC$ is a host algebra for the Oscillator group
\[  G^\sharp:= \Heis(X,\sigma) \rtimes_\alpha \R\,.  \]
It carries a subset of the representations whose restrictions to the distinguished one-parameter group
is positive. By \cite[Thm.~1.3]{Ch68},
these representations are direct sums of the Fock representation.
\end{ex}

\begin{rem} \mlabel{rem:5.2} For the Heisenberg group
$G = \Heis(X,\sigma)$ we consider the complexification
\[ G_\C = \C^\times \times X_\C \quad \mbox{ with multiplication } \quad
(a,z)(b,w) = (ab e^{-i\sigma(z,w)/2}, z + w),\]
where $\sigma$ also denotes the complex bilinear skew-symmetric extension to~$X_\C$.
Then $\alpha$ extends to a holomorphic action
$\alpha \: \C^\times \to \Aut(G_\C)$ by
$\alpha_c(a,z) = (a, c z)$, so that we can form the complex group
\[  G^\sharp_\C = \Heis(X,\sigma)_\C \rtimes_\alpha \C.  \]
It contains the open subsemigroup
\[ S := \{ (a,z,c) \,\mid\, \Im c > 0 \} = \Heis(X,\sigma)_\C \times \C_+.\]
According to \cite[Thm.~5.4]{Ze15}, a smooth positive covariant
unitary representation of $G^\sharp$
extends naturally to a holomorphic representation
$U^\C \: S \to \cB(\cH).$
By analytic extension, it then follows that the $C^*$-algebra
$\cC_1 := C^*(U^\C(S)) \subeq \cB(\cH)$ is generated by the image of the subset
$\Heis(X,\sigma) \times i(0,\infty) \subeq S$
under $U^\C$. As $U^\C(\C_+)= U(\cL)$ and
$U^\C(g s) = \pi(g) U^\C(s)$ for $g \in G, s \in S,$
it follows that \
\[ \cC_1 = C^*(\pi(\cA)U(\cL)) = \cC.\]
\end{rem}

The important property in  Example~\ref{OscRegHeis} was the norm continuity of the maps
(\ref{Lcont}). Let us generalize this:

\begin{prop} \mlabel{Gcrosscont} {\rm(Host algebras for covariant group representations)}
Let $G$ be a topological group  and let
$\alpha \: \R \to \Aut(G)$ be a homomorphism defining a continuous action of $\R$ on $G$.
 Let $(\pi, U)$ be a covariant representation of $(G,\R, \alpha)$ on $\cH$ which is cross, in the sense that
\[\pi(G)U(\cL) \subeq U(\cL)\cB(\cH),\quad \mbox{ where } \quad \cL=C^*(\R).\]
 If the maps
$\pi^L \: G \to  \cB(\cH),\;\; g \mapsto \pi(g) L,\; L\in U(\cL), $
are norm continuous, then the algebra
\[ \cC := C^*(\pi(G)U(\cL))  \]
is a host algebra for the group
$G^\sharp :=  G \rtimes_\alpha \R$
and it carries a subset of the (continuous) representations of $G^\sharp$.
\end{prop}

\begin{prf} From the cross condition we construct the crossed product host
\[
\cC := C^*(\pi(G)U(\cL)) = C^*(\pi(C^*(G_d))U(\cL))  \subeq U(\cL)\cB(\cH)\cap  \cB(\cH)U(\cL)\,
\]
for the action $\tilde\alpha:\, \R \to \Aut(C^*(G_d))$ by $\tilde\alpha_t(\delta_g)=\delta_{\alpha_t(g)}$
(recall that $C^*(G_d)$ is generated by unitaries $(\delta_g)_{g\in G}$
satisfying the group law in $G$).
As there is an injection from   the covariant representations of $\tilde\alpha$  to
 the representations of $   G_d \rtimes_\alpha \R$, it follows that
$\cC$ is a host algebra for $G^\sharp$. The injection is that a covariant representation $(\tilde\pi, \tilde{U})$
on $\tilde\cH$ (corresponding to a representation $\rho:\cC\to  \cB(\tilde\cH)$)
is mapped to the representation $\gamma$ of  $   G_d \rtimes_\alpha \R$ defined by $\gamma(g,t):=\tilde\pi(\delta_g)\tilde{U}_t$.
To get a representation of $ G^\sharp=  G \rtimes_\alpha \R$, it suffices to prove that the map $g\mapsto \tilde\pi(\delta_g)$ is
strong operator continuous on $G$. As $\tilde{U}(\cL)$ is non-degenerate, $\tilde{U}(\cL)\tilde\cH$ is dense in
$\tilde\cH$, hence it suffices to observe that for
$L\in\cL$ and $v\in\tilde\cH$ the map
\[
g\mapsto \tilde\pi(\delta_g)\tilde{U}(L)v = \rho\left(\pi(g)U(L)\right)v
\]
is continuous, using the assumed continuity of the maps~$\pi^L$.
\end{prf}

\begin{rem}
If $U(\cL)\subseteq\cK(\cH)$, then the maps
$\pi^L \: G \to  \cB(\cH),\;\; g \mapsto \pi(g) L,\;\; L\in U(\cL), $ are automatically norm continuous.
For example, if $G$ is a compact Lie group, and $U(t)=\exp(it\Delta)$,
where $\Delta$
is the Laplacian of $\pi:G\to \cB(\cH)$,
then $U(\cL)\subseteq\cK(\cH)$ is equivalent to the finiteness of the
multiplicities of all irreducible subrepresentations of~$\pi$
(\cite[Prop.~C.5]{GrN14}).
\end{rem}

It is easy to make further examples, e.g.\
Examples~\ref{ex:5.13}, \ref{VirExmp}, \ref{TwLoopExmp} below, and:
\begin{ex}
\label{NUHhost}
Let $\cH$ be an infinite dimensional Hilbert space, and let $\cF(\cH)$ be the symmetric Fock space built on it.
Let $G=\U(\cH)$ be the unitary group with the norm topology.
Then by second quantization there is a representation ${\Gamma:G\to \U(\cF(\cH))}$ on Fock space, which
preserves the  $n$-particle  subspaces $\cF(\cH)_n \subeq \cF(\cH)$.
Let $P_n$ denote the projection onto $\cF(\cH)_n$. These are the eigenspaces for the action
$U:\T\to \U(\cF(\cH))$ by $U_tv=t^nv$ for $t\in\T$, $v\in \cF(\cH)_n$.  As this commutes with
$\Gamma(G)$, conjugation with $U_t$ is trivial on $\Gamma(G)$, and the cross condition is
immediate as
\[\Gamma(G)U(\cL) = U(\cL)\Gamma(G).  \]
Moreover, the restriction of  $\Gamma(G)$ to an $n$-particle  space $\cF(\cH)_n$ is just a
symmetrized  $n$-fold tensor product of the defining representation of $G$, which is norm continuous.
Thus, we obtain norm continuity of the maps $g\mapsto\Gamma(g)P_n$, hence continuity of
$g \mapsto \pi(g) L$ for $ L\in U(\cL)$. By
Proposition~\ref{Gcrosscont}, we now obtain a host algebra
for $G^\sharp=G\times\T$. In fact,  as $U(\T)$ coincides with $\Gamma(Z(G))$ ($Z(G)\cong \T$),
it follows that we have a host algebra for $G=\U(\cH)$, endowed with the norm topology.
\end{ex}

The form of the host algebra obtained in Proposition~\ref{Gcrosscont} makes it easier to
check the cross condition for any further $C^*$-action of $G^\sharp$:

\begin{prop} \mlabel{Gcrosssubgp}
Let $G$ be a topological group  and let
$(\pi_G, U)$ be a covariant representation of a continuous action
$(G, \R,\gamma)$, and assume that $(\pi_G, U)$ is cross for $\cL=C^*(\R)$, and
that for each $L \in U(\cL)$, the map
$\pi_G^L \: G \to  \cB(\cH), g \mapsto \pi_G(g) L$ is
norm continuous, so that we have the host algebra
\[ \cL^\sharp:= C^*(U(\cL)\pi_G(G))
=\br U(\cL)\pi_G(G). \subeq \cB(\cH)\]
for the semidirect product group $G^\sharp = G \rtimes_\gamma \R$.\\[1mm]
Let $(\cA, G^\sharp,\alpha)$ be a $C^*$-action (possibly singular).
Then a   covariant $\cL^\sharp$-representation $(\pi_\cA, V)$ of $(\cA, G^\sharp, \alpha)$ on $\cH'$
is cross with respect to  $\cL^\sharp$ if and only if its restriction to the subsystem  $(\cA, \R, \alpha\restriction\R)$
 is cross with respect to  $\cL=C^*(\R)$.
\end{prop}

\begin{prf}
A covariant $\cL$-representation $(\pi_\cA, V)$ we must have that
$V=\pi'_G\times U'$ for $G^\sharp=G \rtimes_\gamma \R$ for some
$\gamma\hbox{-covariant}$ representation
$(\pi'_G, U')$. Thus $V(\cL^\sharp)=\br U'(\cL)\pi'_G(G).$
(where $U'(\cL)=U'(C^*_+(\R))$ by positivity). Assume first that the
restriction of $(\pi_\cA, V)$ to the subsystem  ${(\cA, \R, \alpha\restriction\R)}$
 is cross with respect to  $\cL$, i.e. $\pi_\cA(\cA)U'(\cL) \subseteq U'(\cL)\cB(\cH').$
Thus we have that
\[
\pi_\cA(\cA)V_{\cL}(\cL^\sharp)=\pi_\cA(\cA)\br U'(\cL)\pi'_G(G).
\subseteq U'(\cL)\cB(\cH')\subseteq V(\cL^\sharp)\cB(\cH'),
\]
where the last step follows because $U'(\cL)\subset V(\cL^\sharp)$. Thus $(\pi_\cA, V)$ is cross for
 $(\cA, G^\sharp, \alpha)$.

 Conversely, assume that  $(\pi_\cA, V)$ is cross for
 $(\cA, G^\sharp, \alpha)$, i.e. $\pi_\cA(\cA)V(\cL^\sharp)\subseteq V(\cL^\sharp)\cB(\cH')$. As
 $U'(\cL)\subset V(\cL^\sharp)$ we get
\[
\pi_\cA(\cA)U'(\cL)\subseteq V(\cL^\sharp)\cB(\cH')\subseteq U'(\cL)\cB(\cH')
\]
because $U'(\cL)$ acts non-degenerately on $V(\cL^\sharp)$. Thus the restriction to $(\cA, \R, \alpha\restriction\R)$
is also cross.
\end{prf}

Thus it suffices to check the cross condition on the one-parameter subgroup in the semidirect product
$G \rtimes_\gamma \R$.

\begin{ex}
We continue Example~\ref{NUHhost} to realize the assumptions of Proposition~\ref{Gcrosssubgp}.
Let $\cH$ be an infinite dimensional Hilbert space, and let $G=\U(\cH)$ be the unitary group with the norm topology. Let $G$ act on the Weyl algebra $\cA=\ccr \cH,\sigma.$, where
$\sigma(x,y):={\rm Im}{\langle x,y\rangle}$
by the usual symplectic action
$\alpha:G\to\Aut(\cA)$ determined by $\alpha\s V.(\delta_x):=\delta_{Vx}$, $V\in G=\U(\cH)$.
On the bosonic Fock space $\cF(\cH)$, we then have the second quantization representation
${\Gamma:G\to \U(\cF(\cH))}$ which is covariant for the Fock representation of the Weyl algebra
$\pi_F:\cA\to \cB(\cF(\cH))$. The restriction of $\Gamma$ to $Z(G)\cong \T$ defines the action
$U:\R\to \U(\cF(\cH))$ by $U_s=\Gamma(e^{is})$ for which Example~\ref{NUHhost} produces the host
algebra
\[ \cL^\sharp:= C^*(U(\cL)\Gamma(G))
=\br U(\cL)\Gamma(G). \subeq \cB(\cH),\quad\hbox{where}\quad \cL=C^*(\R).\]
To establish that  $(\pi_F, \Gamma)$ of $(\cA, G, \alpha)$
is cross with respect to  $\cL^\sharp$, we only need to check by Proposition~\ref{Gcrosssubgp}
that $(\pi_F, U)$ of $(\cA, \R, \alpha)$
is cross with respect to  $\cL$. As $U_s=\exp(is\,\dd\Gamma(\1))$ and $\dd\Gamma(\1)$ is the number operator,
this has already been verified in \cite[Example~6.11]{GrN14}. Thus $(\pi_F, \Gamma)$ of $(\cA, G, \alpha)$
is cross with respect to  $\cL^\sharp$.
\end{ex}

\subsection{Smoothing operators for unitary Lie group representations}
\label{subsec:smop}

For Lie groups, we will have to replace the continuity condition above by a smoothness condition. We want to examine this situation.

\begin{defn}\mlabel{def:1.1a}
Let $G$ be a Lie group (possibly infinite dimensional)
with Lie algebra $\g$ and a smooth exponential function $\exp \: \g \to G$.
We write $\g^*$ for the space of continuous linear functionals on
$\g$ (which carries the structure of a locally convex space),
endowed with the weak-$*$-topology.

Let $U \: G \to \U(\cH)$ be a continuous unitary representation
of $G$ and $\cH^\infty\subeq \cH$ be the subspace  of smooth vectors.
We say that $U$ is {\it smooth} if $\cH^\infty$ is dense, which is always the case
if $G$ is finite-dimensional. If $U$ is smooth, then
we obtain for each $X \in \g$ an essentially self-adjoint operator
\[ -i \dd U(X) \: \cH^\infty \to \cH^\infty,
\quad \dd U(X)v := \frac{d}{dt}\Big|_{t= 0} U(\exp tX)v.\]
We write $\partial U(X) := \oline{\dd U(X)}$ for the closure of the skew-symmetric
operator $\dd U(X)$ on $\cH^\infty$. Its domain
$\cD(\partial U(X))$ is the set of differential vectors for the
unitary one-parameter group $t \mapsto U(\exp tX)$ and for any such vector $v$, we have
\[ \partial U(X)v = \frac{d}{dt}\Big|_{t= 0} U(\exp tX)v.\]
We write
\[ \cD^\infty(\partial U(X)) := \bigcap_{n \in \N} \cD((\partial U(X)^n) \]
and note that this space coincides with the smooth vectors of the one-parameter group
$U^X_t := U(\exp t X)$. Then $\cD^\infty(\partial U(X)) \supeq \cH^\infty$, and mostly
this inclusion is proper.
\end{defn}

We have the following characterization of
{\it smoothing operators} $A \in \cB(\cH)$, i.e.,
operators mapping $\cH$ into $\cH^\infty$
(\cite[Thm.~2.1]{NSZ17}):

\begin{thm} \mlabel{thm:smoothop} {\rm(Characterization Theorem for smoothing operators)}
For a  smooth unitary representation $(U, \cH)$ of a Fr\'echet--Lie group~$G$
with smooth exponential function and $A \in \cB(\cH)$,
the following are equivalent:
\begin{itemize}
\item[\rm(i)] The map $G \to \cB(\cH), g \mapsto U(g) A$ is smooth.
\item[\rm(ii)] $A\cH \subeq \cH^\infty$, i.e., $A$ is a smoothing operator.
\item[\rm(iii)] All operators $A^* \dd U(X_1) \cdots \dd U(X_n)$, with
$X_1, \ldots, X_n \in \g$, $n \in \N$, extend as bounded operators to $\cH$ from $\cH^\infty$.
\end{itemize}
If these conditions are satisfied, then all operators
$\dd  U(X_1) \cdots \dd U(X_n) A$,  $X_1, \ldots, X_n \in \g$, are bounded.
\end{thm}
Smoothing operators are very useful, as we show below.
\begin{prop} \mlabel{prop:5.8} Let $G$ be a finite dimensional Lie group,
$H \subeq G$ be a closed subgroup and $(U,\cH)$ be a continuous
unitary representation of $G$.
If $U(C^*(H))$ contains a dense subspace consisting of smoothing operators for $G$, then
$U(C^*(H)) \subeq U(C^*(G))$.
\end{prop}

\begin{prf} Since $U(C^*(H))$ consists of multipliers of $U(C^*(G))$,
the assertion follows from Lemma~\ref{lem:contop} below.
\end{prf}

\begin{lem} \mlabel{lem:contop}
Let $(U, \cH)$ be a continuous unitary representation
of the finite dimensional Lie group $G$, let $\cL = C^*(G)$ and  let
\[ \cM := \{ A \in \cB(\cH) \,\mid\, A U_\cL(\cL) + U_\cL(\cL) A \subeq U_\cL(\cL) \} \]
be the relative multiplier algebra of the $C^*$-algebra $U_\cL(\cL)\subseteq  \cB(\cH)$.
Every $A \in \cM$ for which the map
\[ U^A \: G \to \cB(\cH), \quad g \mapsto U_g A \]
is continuous is contained in the ideal $U_\cL(\cL) \trile \cM$.
\end{lem}

\begin{prf} Let $(\delta_n)_{n \in \N}$ be a $\delta$-sequence in $C_c(G)$,
i.e., $\int_G \delta_n = 1$ and the supports of the $\delta_n$ shrink to~$\1$.
Then the continuity of $U^A$ implies that
$U_\cL(\cL) \ni U_\cL(\delta_n) A  \to A$, so that the assertion follows from
the closedness of $U_\cL(\cL)$.
\end{prf}

We also recall the following theorems from \cite{NSZ17}:

\begin{thm} {\rm(Subgroup Host Algebra Theorem; \cite[Thm.~4.1]{NSZ17})}
 \mlabel{thm:5.11}
Let $(U, \cH)$ be a unitary representation of the metrizable
  Lie group $G$ and  $\iota_H \colon H \to G$ be a morphism
of Lie groups where $\dim H < \infty$ and $U^H := U \circ \iota_H$
satisfies
\[ \cH^\infty = \cH^\infty(U^H). \]
Then
$C^*\Big(U(G) U^H(C^\infty_c(H))\Big) = C^*\Big(U(G) U^H(C^\infty_c(H))U(G)\Big)$
is a host algebra for a class of smooth representations of $G$.
\end{thm}

\begin{thm} {\rm(\cite[Thm.~4.2]{NSZ17})} \mlabel{thm:5.1b}
Let $(U, \cH)$ be a unitary representation of the
metrizable Lie group $G$. If
\[ \cH^\infty = \cD^\infty(\partial U(X_0)) \quad \mbox{ for some } \quad X_0 \in \g
\quad \mbox{ with } \quad
\sup\Spec(i\partial U(X_0)) < \infty, \]
then $e^{i\partial{U}(X_0)}$ is a smoothing operator and
$\cL:= C^*\big(U(G) e^{i\partial{U}(X_0)}\big) = C^*\big(U(G) e^{i\partial{U}(X_0)} U(G)\big)$
is a host algebra for a class of smooth representations
$(\rho_G,\cK)$ of $G$ satisfying
\begin{equation}
  \label{eq:spec-esti}
\sup\Spec(i\partial{\rho_G}(X_0)) \leq \sup\Spec(i\partial{U}(X_0)).
\end{equation}
\end{thm}

\begin{rem} (a) If the group $G$ in Theorem~\ref{thm:5.1b}
is finite dimensional, then Proposition~\ref{prop:5.8} implies that
the holomorphic semigroup $e^{\C_+ \partial U(X_0)}$ is contained in $U(C^*(G))$.
The subalgebra of $U(C^*(G))$ generated by this semigroup coincides with the
image of $C^*(\R)$ under the integrated representation of
the one-parameter group $U_{X_0}(t) := U(\exp t X_0)$.

(b) (\cite{NSZ17}) Assume that $G$ is finite dimensional and let
$X_1, \ldots, X_n$ be a basis of its Lie algebra $\g$.
Then
\[ \cH^\infty = \cD^\infty(\Delta) \quad \mbox{ for } \quad
\Delta := \oline{\sum_{j = 1}^n \dd U(x_j)^2}\]
(\cite[Cor.~9.3]{Nel69}) and since $\Delta$ has non-positive spectrum,
it follows that the contraction semigroup $(e^{t\Delta})_{t > 0}$, which is an
abstract version of the heat semigroup on $L^2(G)$, consists
of smoothing operators.

(c) By Theorems~\ref{thmsbsmoothonegen} and \ref{thm:5.2} below,
there is a rich source of smoothing operators in some representations of  $G$.
\end{rem}

Let $G$ be a Lie group (possibly infinite dimensional)
 with an exponential function and let
$\alpha \: \R \to \Aut(G)$ be a homomorphism defining a continuous action of $\R$ on $G$.
Let $(\pi, U)$ be a positive covariant representation of $(\cA,\R, \alpha)$ which is
smooth as a representation of $G^\sharp = G \rtimes_\alpha \R$.
If $U_t = e^{itH}$, then we write $U_z := e^{izH}$, $\Im z \geq 0$, for its holomorphic
extension to $\oline{\C_+}$, where $\C_+ = \{ z \in \C \,\mid\, \Im z > 0\}$ denotes the
upper half plane.
Note that $C^*\{U_z\mid z\in\C_+\}=U(C^*(\R))$ (cf.\ Theorem~\ref{HostSpec} below).

\begin{thm} \mlabel{thm:5.6}
Let $G$ be a Lie group with smooth exponential function and
 let $(\pi, U)$ be a positive covariant representation of
$(G, \R,\alpha)$ for which
\begin{itemize}
\item[\rm(S)] the operators $(U_z)_{z \in \C_+}$ are smoothing operators,
i.e., $G \to \cB(\cH), g \mapsto \pi(g)U_z$ is smooth.
\end{itemize}
Then
\[ \cC:= C^*(\pi(G)U_{\C_+}) \subeq \cB(\cH).\]
is a crossed product host for $(\cA, \cL, \alpha)$, where
$\cA = C^*(G_d)$ and $\cL = C^*_+(\R)$. Moreover, all representations of
$\cC$ define smooth representations of the Lie group $G$.
\end{thm}

\begin{prf}
Observe first that $U(\cL)=U(C^*_+(\R))=C^*(U_{\C_+})$.
In view of Theorem~\ref{thm:smoothop}, (S) implies that
\[ \dd\pi(U(\g)) U_z \subeq \cB(\cH) \quad \mbox{ for } \quad \Im z > 0.\]
Clearly, $U_{\C_+} \subeq \cC$ and $\pi(G) \cC\subeq \cC$, so that we get a homomorphism
$\eta_G \: G \to M(\cC)$ into the multiplier algebra of $\cC$.
Since $\cC$ contains a dense subspace spanned by elements of the form
\[ \pi(g) U_z A, \quad A \in \cC, g\in G, z \in \C_+,\]
the multiplier action has a dense subset of smooth vectors.
In particular, it is strongly continuous.
From this property one already derives that every non-degenerate representation
$(\pi, \cH)$ of $\cC$ defines a smooth representation $\tilde\pi \circ \eta_G \: G \to \U(\cH)$
of $G$ on $\cH$.

Next we observe that, for $z \in \C_+$ and $g \in G$, the map
\[ \R \to \cC, \quad t \mapsto U_t \pi(g) U_z
= \pi(\alpha_t(g)) U_{z+t} \]
is continuous if the action of $\R$ on $G$ is continuous, and smooth if this action is smooth.
This implies that $\cC$ satisfies the cross condition for
$\cL = C^*_+(\R)$ and the discrete group algebra $\cA = C^*(G_d)$.
\end{prf}
The algebra $\cC:= C^*(\pi(G)U_{\C_+})$ in the theorem is therefore a host algebra
for a class of representations of $G^\sharp = G \rtimes_\alpha \R$, as  representations of
$G^\sharp$ are all covariant representations of $(G, \R,\alpha)$.
The next example specifies the class of representations via Theorem~\ref{thm:5.1b}.
\begin{ex} Suppose that $\alpha$ defines a smooth action, so that
$G^\sharp = G \rtimes_\alpha \R$ is a Lie group.
For $X_0 := (0,1) \in \g \times \R$ we then have
$U_z = e^{z \partial U^\sharp(X_0)}$.
Now the holomorphy of the map $\C_+ \to \cB(\cH), z \mapsto U_z$
implies that
$U^\sharp(\exp \R X_0) e^{i\partial U^\sharp(X_0)}$ and
$e^{\C_+\partial U^\sharp(X_0)}$ generates the same norm-closed subspace of $\cB(\cH)$.
Therefore the algebra $\cL$ in Theorem~\ref{thm:5.1b} coincides with the
$C^*$-algebra $\cC = C^*(\pi(G) U_{\C_+})$
in Theorem~\ref{thm:5.6}.
\end{ex}
\begin{rem}
The form of the host algebra obtained in Theorem~\ref{thm:5.6} makes it easier to
check the cross condition for any further $C^*$-action of $G^\sharp$.
In particular, we can adapt Proposition~\ref{Gcrosssubgp} to condition (S)
to also simplify the cross condition for $C^*$-actions of $G^\sharp$.
\end{rem}
\begin{ex} \mlabel{ex:5.14} The assumptions of the preceding theorem are in particular
satisfied for the unitary representation
of the Heisenberg group $G = \Heis(\cH,\sigma)$
on Fock space $\cF(\cH)$ and the automorphism group
$\alpha_t(z,v) = (z, tv)$ for $t \in \T$, corresponding on $\cF(\cH)$
to the number operator (Remark~\ref{rem:5.2}).
In this case $\cC$ can either be obtained from the holomorphic
semigroup representation of $S$ or as in Example~\ref{OscRegHeis}
as a crossed product host.
\end{ex}

The hypotheses of Theorem~\ref{thm:5.6} are realized by the initial Example~\ref{OscRegHeis}. Further examples
where Theorem~\ref{thm:5.6} produces crossed product hosts are given in Example~\ref{ex:5.13},
 and in Section~\ref{sec:8}
(Examples~\ref{VirExmp} and \ref{TwLoopExmp}).

\section{General spectral conditions---Review}
\mlabel{SpecCondRev}

Next we want to examine spectral conditions  for Lie groups $G$ other than $\R$.
The generalization is in two directions; higher dimensions, and non-abelian groups.
We want to analyze the covariant representations $(\pi, U)\in\Rep(\alpha, \cH)$
of $(\cA, G, \alpha)$, where $U$ satisfies a spectral condition
of a type we will specify below.
Lie group representations
satisfying such spectral conditions are of fundamental importance in
physics~\cite{Bo96, Ot95, H92, LM75},
harmonic analysis~\cite{Ol82, Ol90, HO96, Ne99, Ne10, Ne12, Ne14b,  Ne17}
and operator algebras (cf.~\cite[Ch.~8]{Pe89}).
In this section we will review some
Lie group representations satisfying spectral conditions in order to use it
in the next sections.
Proofs and further details are in~\cite{Ne99}.
This theory extends to some infinite dimensional Lie groups
(cf.~\cite{Ne08, MN12, Ze15}).

Below we will need the following terminology:
\begin{defn}\mlabel{defSpec}
Let $G$ be a Lie group (possibly infinite dimensional)
with Lie algebra $\g$ and a smooth exponential function $\exp \: \g \to G$.
We write $\g^*$ for the space of continuous linear functionals on
$\g$, endowed with the weak-$*$-topology.

(a) Let $U \: G \to \U(\cH)$ be a smooth unitary representation
of $G$ and $\cH^\infty\subeq \cH$ be the subspace  of smooth vectors.
Each smooth unit vector $v \in \cH^\infty$ defines a
continuous linear functional $\Phi(v) \in \g^*$ by
\[ \Phi(v)(X):=-i\la \dd U(X)v,\,v\ra \quad \mbox{ for } \quad X \in \g\]
and the {\it momentum set of $U$} is defined by
\[ I_U:=\hbox{weak-$*$-closed convex hull of\;}
\big\{\Phi(v)\,\mid\,v\in\cH^\infty,\,\|v\|=1\big\}\subeq \g^*.\]
This is an $\Ad^*(G)$-invariant weak-$*$-closed convex subset of $\g^*$.

(b) Given a subset $C \subeq \g^*$, we say that $U$ satisfies the
{\it $C$-spectral condition} if $I_U \subeq C$.
We write
\[
{\rm Rep}_C(G,\cH):=\big\{
U\in{\rm Rep}(G,\cH)\,\mid\, I_U \subeq C \big\}\]
for the set of smooth representations on $\cH$ satisfying the
$C$-spectral condition.

(c) The fact that the operators $i\dd U(X)$ on $\cH^\infty$
are essentially selfadjoint implies that
  \begin{equation}
    \label{eq:spec}
\inf\big({\rm spec}(-i\partial U(X))\big)=\inf\ I_U(X) \quad \mbox{ for } X \in \g.  \end{equation}
 The smooth representation $U$ is said to
be {\it semibounded} if the open convex cone
\[ W_U := \{ X \in \g\,\mid\, Y \mapsto \inf I_U(Y) \ \
\mbox{is bounded below in a neighborhood of} \ X\}\]
is not empty. If $\g$ is finite-dimensional, this cone $W_U$ coincides
with the interior of
\[
B(I_U):=\{ X \in \g\,\mid\, \inf\ I_U(X) >-\infty\}.
\]
For finite dimensional groups $G$, the cone
$W_U$ is non-empty (i.e.\ $U$ is semibounded) if and only if
the convex set $I_U$ contains no affine lines (\cite[Prop.~V.1.15]{Ne99}).

For a general locally convex space $V$, we use the following notation:
for any subset $S \subeq V^*$ (the topological dual),
we write
\[ B(S):=\{ v \in V\,\mid\,\inf \la S, v \ra >-\infty\}
\quad\hbox{and}\quad
S^\star:=\{v \in V\,\mid\,\la S, v \ra \subseteq [0,\infty)\}\]
and note that these are both convex cones, the {\it dual cone} $S^\star$ is closed,
and $B(S)=B\big(\overline{{\rm conv}(S)}\big)$.
Clearly, if $S_1\subset S_2\subeq V$, then $S_2^\star\subeq S_1^\star$ and $B(S_2)\subeq B(S_1)$.
If $W\subeq V$ is a convex cone, then $B(W)=W^\star$ in $V^*$,
and $W^{\star\star}=\overline{W}$.

(d) For the set $W^\star\subeq \g^*$, we derive from \eqref{eq:spec} that for
a smooth unitary representation $U \: G \to \U(\cH)$ to satisfy the
 $W^\star$-spectral condition, means that all the operators
$-i\partial U(X)$, $X \in W$, have non-negative spectrum.
Such representations are also called {\it $W$-dissipative}
(cf.\ \cite[Def.~X.3.11]{Ne99}).
\end{defn}

\begin{rem} (a) By definition, a unitary representation $(U,\cH)$ of $G$
is $W$-dissipative if and only if it satisfies the
$W^\star$-spectral condition (Definition~\ref{defSpec}).
If the cone $W$ is open, $W$-dissipativity implies semiboundedness.
If, conversely, $(U,\cH)$ is semibounded and we extend $U$ by
$U^\sharp(z,g) := z U(g)$ to the trivial central extension $G^\sharp := \T \times G$,
then the dissipative cone
\[ W_{U^\sharp}^+ := \{ X \in \g^\sharp \,\mid\, - i\partial U^\sharp(X) \geq 0 \} \]
has interior points. This simple construction essentially reduces the
investigation of semibounded representations to $W$-dissipative representations
with respect to open invariant cones $W \subeq \g$.

(b) If $(U,\cH)$ is a semibounded unitary representation
and $\iota_H \colon H \to G$ a morphism
of Lie groups for which the image of the derived homomorphism
$\L(\iota) \:  \fh \to \fg$ intersects $W_U$, then
$U^H := U \circ \iota_H$ is also semibounded.
This holds in particular for restrictions to
Lie subgroups $H \subeq G$ for which $\fh \cap W_U \not=\eset$.
\end{rem}

\subsection{A host algebra for $C$-spectral representations}
\label{smoothHost}

We start by quoting some results from \cite{NSZ17}.
Semibounded representations have a number of special useful properties, e.g.

\begin{thm}\mlabel{thmsbsmoothonegen}
{\rm(Zellner's Smooth Vector Theorem;  \cite[Thm.~3.1]{NSZ17})}
Let $G$ be a Lie group (possibly infinite dimensional)
with a smooth exponential function.
If $U:G\rightarrow\U(\cH)$ is a semibounded unitary
representation and $X_0\in W_U$, then
$\cH^\infty=\cD^\infty(\partial U(X_0))$
and all operators of the form
$\int_\R f(t) U(\exp t X_0)\, dt$, $f \in C^\infty_c(\R)$, are smoothing.
\end{thm}

From these smoothing operators,
one can even obtain  host algebras which are full for $C$-spectral
representations of infinite dimensional Lie groups. The following theorem
is the primary source for host algebras constructed from
semibounded representations:

\begin{thm} \mlabel{thm:5.2} {\rm(\cite[Thm.~4.3, Prop.~4.1]{NSZ17})}
Let $(U, \cH)$ be a semibounded unitary representation of the metrizable
Lie group $G$ with smooth exponential function.
Then the operators $e^{i\partial U(X)}$, $X \in W_U$, are smoothing operators
and $\cL :=
C^*\big(U(G) e^{i\partial{U}(W_U)}\big)= C^*\big(U(G) e^{i\partial{U}(W_U)} U(G)\big)$
is a host algebra for a class of smooth representations
$(\rho,\cK)$ of $G$ satisfying
$I_{\rho} \subeq I_U$. For every $X \in W_U$, we have
$\cL = C^*\big(U(G) e^{i\partial{U}(X)}\big)$.
\end{thm}

From the preceding theorem one obtains the following general existence
result on host algebras for $C$-spectral representations:

\begin{cor} \mlabel{cor:4.9} {\rm(\cite[Cor.~4.1]{NSZ17})}
{\rm(A full host algebra for $C$-spectral representations)}
Let $G$ be a metrizable
Lie group with smooth exponential function and
let $C \subeq \g^*$
be a weak-$*$-closed $\Ad^*(G)$-invariant subset
  which is semi-equicontinuous in the sense that its support function
\[ s_C \colon \g \to \R \cup \{\infty\}, \quad
s_C(x) := \sup \la C, -x \ra = -\inf \la C, x \ra \]
is bounded in a neighborhood of some point $X_0 \in \g$.
Then there exists a
host algebra $(\cL,\eta)$
of~$G$ whose representations correspond to those
semibounded unitary representations $(U, \cH)$ of~$G$ for which
$s_U(X) := \sup\Spec(i\partial U(X)) \leq s_C(X)$, i.e., $I_U \subeq C$.
\end{cor}

From Theorem~\ref{thmsbsmoothonegen} we get in particular:
\begin{ex} \mlabel{ex:5.13}
If $(\pi,\cH)$ is a semibounded representation of the
Lie group $G$ and $X \in W_\pi$, then we may consider the
$\R$-action on $G$ given by $\alpha_t^X(g) := \exp(tX) g \exp(-tX)$.
With $U_t := \pi(\exp tX)$ we then obtain a covariant
representation $(\pi, U)$ of $(G,\R, \alpha^X)$.
Combining Theorems~\ref{thmsbsmoothonegen}
and \ref{thm:5.11} implies that condition~(S) in Theorem~\ref{thm:5.6} is satisfied
and the last statement in Theorem~\ref{thm:5.2} implies that the crossed product host
$\cC = C^*(\pi(G)U_{\C_+})$ from there
coincides with the host algebra $\cL = C^*(\pi(G) e^{i \partial(W_U)})$
from Theorem~\ref{thm:5.2}. As $U_t = \pi(\exp tX)$, the multiplier action $\eta_{G^\sharp}$ of
$G^\sharp = G \rtimes_\gamma \R$ on $\cC$ satisfies
\[
\{(\exp tX,-t)\in G^\sharp\,\mid\,t\in\R\}\subseteq{\rm Ker}(\eta_{G^\sharp})\,,
\]
 hence $\eta_{G^\sharp}$ coincides with the multiplier action  $\eta_G$ of
 $G$ on $\cL=\cC$.
\end{ex}

\subsection{Analysis of host algebras for $C$-spectral representations}
\label{OlshHost}

In this subsection $G$ denotes a connected finite dimensional Lie group
and $\g = \L(G)$ its Lie algebra.
Specific groups we have in mind for $G$
are $\R^n$, $\SL_2(\R)$ and the conformal group $\SO_{2,n}(\R)$
of $n$-dimensional Minkowski space.

\begin{defn}
Fix an open  convex $\Ad(G)$-invariant nonempty
cone $  W \subeq \g$.
If $W$ contains no affine lines, then, for each element $X \in W$,
the operator $\ad X$ on $\g_\C$  is diagonalizable
with $\Spec(\ad X) \subeq i \R$, i.e.\ the cone
$W$ is {\it elliptic} (cf.~\cite[Prop.~VII.3.4]{Ne99}).
We call the cone $W$ {\it weakly elliptic} if
$\Spec(\ad X) \subeq i \R$ for all $X \in W$,
but $\ad X$ need not be diagonalizable.
\end{defn}
 Let $U \: G \to \U(\cH)$ be a continuous unitary representation
of $G.$
If $\dd U$ is injective (e.g.\ if $\ker U$ is discrete), then
$I_U^\perp = {\rm ker}\,\dd U = {0}$ is trivial, hence
by \cite[Prop.~VII.3.4(c)]{Ne99}, any open cone
$W \subeq B(I_U)$ is elliptic.
In particular, if $U$ is also semibounded, then $W_U$ is elliptic.
This produces a large set of examples of elliptic  cones, but here are
two particularly important examples for applications.
 \begin{ex}
 \label{Lightcone}
 (Abelian groups)
  In physics an important class of representations consists of
 continuous unitary representations $U:\R^4\to \U(\cH)$ such that the
 joint spectrum of the generators for the  coordinate subgroups $\R$ in $\R^4$
 is in the closure
\[ \oline{V_+} =\big\{\b x.\in\R^4\mid x_0^2-x_1^2-x_2^2-x_3^2
 \geq 0,\; x_0\geq 0\big\}\]
of the open forward light cone
\[ V_+ := \{ x = (x_0, {\bf x}) \in \R^4 \mid x_0 > 0, x_0^2 > {\bf x}^2 \}.  \]
With $G=\R^4$, we  get $\g\cong\R^4\cong\g^*$
 and we take $W := V_+ \subset\R^4\cong\g$
to be the interior of the dual light cone, and thus $W^\star=\oline{V_+}$.
As $\R^4$ is abelian, $W$ is certainly elliptic.
Then the  $W^\star$-spectral condition means that all the operators
$-i\partial U(X)$, $X \in W$, have non-negative spectrum,
and hence that the joint spectrum of $U$ is in $W^\star=\oline{V_+}$
(cf.\ \cite[Def.~II.6.3]{Bo96}).
\end{ex}
 \begin{ex}
 \label{Sympcone}
 (Non-abelian groups)
Let  $(X,\sigma)$ be a finite dimensional  symplectic space and consider
the associated symplectic group $G=\Sp(X,\sigma)\subset \GL(X)$.
For any element of its Lie algebra $A\in\sp(X,\sigma)\subset \gl(X),$ we define
the Hamiltonian quadratic form $H_A:X\to\R$ by 
\[
H_A(x):=\hlf\sigma(Ax,x).
\]
Then we write $H_A\gg 0$ if and only if $\sqrt{H_A}$ is a Hilbert norm on $X$, i.e.\
$H_A$ is a positive definite quadratic form.
By \cite[Thm.~3.1.19]{AM78}, this is equivalent to requiring $X$ to have
a complex structure $J\in \Sp(X,\sigma)$ commuting with $A$ such that
\[\langle x,y\rangle := \hlf\sigma(Ax,y)-\frac{i}{2}\sigma(Ax,Jy)\]
defines a Hilbert  inner product on $X$.
Define an open invariant convex cone
$W\subset\sp(X,\sigma)$ by
\[
W:=\{A\in\sp(X,\sigma)\,\mid\,H_A\gg 0\}
\]
(cf.\ \cite[\S 6.2]{Ne10}).
Then $W$ is nonempty (in canonical coordinates it
contains the Hamiltonian function of the harmonic oscillator) and $B(W)=W^\star.$
Moreover, the closed invariant convex cone $\overline{W}$ is pointed, hence
$W$ is elliptic. Since
$A \mapsto H_A$ is a linear bijection from $\sp(X,\sigma)$ onto
the space of quadratic forms on~$X$ and $\oline W$ corresponds to the
positive semidefinite ones, this cone is pointed.
That $\oline W$ is  the unique non-zero closed invariant convex
cone in $\sp(X,\sigma)$ up to sign follows from \cite[Thm.~VIII.3.21]{Ne99}
and the fact that $C_{\rm min} = C_{\rm max}$ (in the notation from loc.\ cit.)
holds for the Lie algebra $\sp(X,\sigma)$.

For a representation $U \: G \to \U(\cH)$ to satisfy
the  $W^\star$-spectral condition, means that all the operators
$-i\partial U(X)$, $X \in W$, have non-negative spectrum, which is a desirable property
for quantum mechanics. 

For extensions of this correspondence to infinite dimensional
symplectic spaces we refer to \cite{Ne10}, and in particular to \cite{NZ13},
where oscillator groups $\Heis(X,\sigma) \rtimes_\alpha \R$ with non-trivial
positive covariant representations are characterized.
\end{ex}

In the rest of this section we will summarize the theory involved with
obtaining host algebras $(\cL,\eta)$ for the representations satisfying the
 $W^\star$-spectral condition, i.e.\ such that  $\Rep(G,{\cal H})_\eta={\rm Rep}_{W^\star}(G,\cH)$ for each Hilbert space
 $\cH$. Such a host algebra can be obtained directly from
a factor algebra of the group algebra $C^*(G)$, but we describe here a route
for obtaining it
which generalizes to some infinite dimensional Lie groups (hence for which $C^*(G)$
is undefined); see~\cite{Ne08} and the applications to host algebras in \cite{NSZ17}.

For any connected Lie group
$G$, and a  nonempty open convex invariant
cone $ W \subeq \g$ not containing affine lines (or more generally if
$W$ is weakly elliptic), there is a very useful associated object,
the  Olshanski semigroup of $G$ and $W$, which we now define.

\begin{defn}\mlabel{def:OlShSgp}

(a) First we recall that, for any Lie group $G$, there exists a
pair $(G_\C, \eta_G)$, consisting of a complex Lie group $G_\C$ and a
morphism $\eta_G \: G \to G_\C$ of real Lie groups, with the following
universal property: For every smooth homomorphism $\alpha \: G \to H$ of $G$ to a complex
Lie group $H$, there exists a unique morphism $\alpha_\C \: G_\C \to H$ of complex
Lie groups with $\alpha_\C \circ\eta_G = \alpha$
(\cite[Rem.~XIII.5.7]{Ne99}). The pair $(G_\C,\eta_G)$ is called the
{\it universal complexification of $G$}. In view of its universal property, it is unique
up to isomorphism.

In order to define the complex conjugation  $\sigma$ on $G_\C,$
we need to consider the construction of the universal complexification
for the case where $G$ is connected.
Let $q_G \:  \tilde G \to G$ be the universal covering group of $G$
and write $(\tilde G)_\C$ for the simply connected Lie group whose
Lie algebra is the complexification
$\g_\C = \g \otimes_\R \C$ of $\g$.
Then the simply connectedness of $\tilde G$ shows that
there exists a unique morphism $\eta_{\tilde G} \: \tilde G \to (\tilde G)_\C$
whose differential is the inclusion $\g \into \g_\C$. There also exists an
antiholomorphic involutive automorphism $\sigma$ of $(\tilde G)_\C$ whose differential
is the complex  conjugation $x + i y \mapsto x - i y$ on $\g_\C$.
That $((\tilde G)_{\C}, \eta_{\tilde G})$ is a universal
complexification of $\tilde G$ follows directly from the simple connectedness
of $(\tilde G)_\C$.

In general the map
$\eta_{\tilde G}$ is not injective, as the example $G = \SL_2(\R)$ with
$(\tilde G)_\C = \SL_2(\C)$ shows.
Let $N \subeq (\tilde G)_\C$ be the smallest closed complex subgroup of
$(\tilde G)_\C$ containing $\eta_{\tilde G}(\ker q_G) \subeq Z((\tilde G)_\C)$.
Then $G_\C:= (\tilde G)_\C/N$ is a complex Lie group and $\eta_{\tilde G} \: \tilde G
\to (\tilde G)_\C$ factors through a morphism
$\eta_G \: G \to G_\C$ of Lie groups.
Then $(G_\C,\eta_G)$ is a universal complexification of $G$.
From the construction it follows that $G_\C$ always carries an antiholomorphic
involutive automorphism $\sigma \: G_\C \to G_\C$ with
$\sigma \circ \eta_G = \eta_G$, and this implies that
$\eta_G(G)$ coincides with the identity component of the group $G_\C^\sigma$
of $\sigma$-fixed points in $G_\C$.
In general the kernel of $\eta_G$ is non-trivial and it may even have positive
dimension, so that $\dim_\C G_\C < \dim_\R G$.

(b) Let $W \subeq \g$ be an open weakly elliptic cone.
Then the {\it Olshanski semigroup} $\Gamma_G(W)$ corresponding to
$G$ and $W$ is defined as follows.

{\bf Case 1:} First we assume that $\eta_G$ is injective, so that we may
consider $G$ as a closed subgroup of $G_\C$. Then Lawson's Theorem
(\cite[Thm.~XIII.5.6]{Ne99}) asserts that
\[ \Gamma_G(W) := G \exp(iW) \subeq G_\C \]
is an open subsemigroup of $G_\C$ stable under the antiholomorphic
involution $s^* := \sigma(s)^{-1}$ and for which the polar map
\[ G \times W \to \Gamma_G(W), \quad (g,X) \mapsto g \exp(iX) \]
is a diffeomorphism. Here we use that the closure
$\oline W$ of a weakly elliptic cone is also weakly elliptic because of the
``semicontinuity'' of the spectrum of operators on finite dimensional spaces.

{\bf Case 2:} If $G$ is simply connected, then
$\ker \eta_G$ is discrete, so that Case 1 applies to
the subgroup $G_1 := G/\ker \eta_G \subeq G_\C$. Then we define
$\Gamma_G(W)$ as the simply connected covering of $\Gamma_{G_1}(W) \subeq G_\C$.
Basic covering theory implies that it inherits the involution $*$ and a
diffeomorphic polar map $G \times W \to \Gamma_G(W)$.
We also write $\exp \: i W \to \Gamma_G(W)$ for the canonical
lift of the exponential function $\exp\res_{iW} \: i W \to \Gamma_{G_1}(W) \subeq G_\C$.

{\bf Case 3:} In the general case we put
$\Gamma_G(W) := \Gamma_{\tilde G}(W)/\ker q_G,$
and one verifies that this inherits the involution and the polar decomposition
from $\Gamma_{\tilde G}(W)$.
\end{defn}

\begin{rem} The  Olshanski semigroup $S=\Gamma_G(W)$ has a range of desirable properties
(cf.\ \cite[Thm.~XI.1.12]{Ne99}). It is a complex manifold with a holomorphic
semigroup multiplication and an antiholomorphic involution, i.e., $(S,*)$ is  a
{\it complex involutive semigroup}.

In terms of the polar decomposition $s = g \exp(iX)$, the involution is given by
\[ s^* = \exp(iX)g^{-1} = g^{-1}\exp(i\Ad(g)X).\]
Furthermore, we have an action of $G$ on $S$ by left and right multiplications
(actually by multipliers of $(S,*)$):
\[ h (g \exp(iX)) = hg \exp(iX) \quad \mbox{ and } \quad
g \exp(iX) h = gh \exp(i\Ad(h)^{-1}X).\]
Each $X \in W$ generates a holomorphic one-parameter semigroup
\[ \gamma_X \: \C_+ := \{ z \in \C \mid \Im z > 0\} \to S, \quad
\gamma_X(a + i b) = \exp(a X) \exp(i bX). \]
The ``boundary values'' of $\gamma_X$ are given by the one-parameter group
of $G$ generated by $X$.
\end{rem}
\begin{rem}
(1) Given the range of conditions we had to assume to construct the Olshanski semigroup,
it is natural to consider the existence question, i.e. when these conditions hold.
The conditions are that $G$ is a connected Lie group, and that
 $W \subeq \g$ is a nonempty invariant open weakly elliptic cone.
 The technical characterization for a Lie algebra to have such a cone $W$, is
 given in \cite[Thm's~VIII.3.6 and VIII.3.10]{Ne99}, but a more constructive
 discussion is given below in Remark~\ref{rem:OlshHost}(c).
Examples are plentiful, some above, and more in the literature (cf. \cite{Ne99, Ne10}).

 (2) If $U$ is a semibounded representation, we defined above the open
 invariant cone $W_U$  in $\g$. If $U$ is injective, then
 $W_U$ is elliptic (\cite[Prop.~VII.3.4(c)]{Ne99})
and there exists a corresponding Olshanski semigroup
$S = \Gamma_G(W_U)$. If $U$ is not injective, it defines an injective semibounded representation
$\overline U$  on the factor
group $Q := G/{\rm ker}(U)$ and we can construct an Olshanski
semigroup~$\Gamma_Q(W_{\overline U})$.
 \end{rem}

The Olshanski semigroup  $\Gamma_G(W)$ has a close relation to the
$W$-dissipative representations.

\begin{thm} \mlabel{Wreps-OS} Let $S = \Gamma_G(W)$
be a complex Olshanski semigroup and
$(U,\cH)$ be a continuous unitary representation of $G$ such that $W \subeq B(I_U)$
(e.g.\ if $U$ is $W$-dissipative).
Then
\[ \hat U \: S \to \cB(\cH), \quad \hat U(g \exp(iX)) := U(g)e^{i\partial U(X)}\quad \mbox{ for }
\quad g \in G, X \in W, \]
defines a holomorphic non-degenerate $*$-representation of $S$ ``extending'' $U$
in the sense that
\[ U(g) \hat{U}(s) = \hat{U}(gs) \quad \mbox{ for }
\quad  g \in G,\, s \in S, \]
and \eqref{eq:spec} implies that
\begin{equation}
  \label{eq:norm}
\big\|\hat{U}(g\exp(iX))\big\|=\big\|e^{i\partial U(X)}\big\|
=e^{-\inf(I_U(X))}.
\end{equation}
Note that $e^{i\partial U(X)}$ always is a positive operator.
\end{thm}

The main claim is proven in  \cite[Thm~XI.2.3]{Ne99}.  Not all
holomorphic $*$-representations of  $\Gamma_G(W)$ correspond to
$W$-dissipative representations. To formulate the exact class, we have to take a
closer look at the norms of the operators $\hat U(s)$.

\begin{defn}\mlabel{absval}
(a) Let $C \subeq \g^*$ be a closed convex $\Ad^*(G)$-invariant subset
with $W\subseteq B(C)$. Then
$\lambda_C:S\to[0,\infty)$ defined by
\[
\lambda_C\big(g\exp(iw))\big):=e^{-\inf(C(w))}\quad \mbox{ for } \quad
g \in G,\, w \in W
\]
is an {\it absolute value} in the sense that
$\lambda_C(s^*)=\lambda_C(s)$ and $\lambda_C(st)\leq\lambda_C(s)\lambda_C(t)$ for
all $s,\,t\in S$ (cf.~\cite[XI.3.2]{Ne99}). It is locally bounded, and its explicit
form further implies that $\lambda_C(s^*s) = \lambda_C(s)^2$ for $s \in S$.

(b) Given an absolute value $\lambda:S\to[0,\infty)$, we say that a holomorphic $*$-representation
$V:S\to \cB(\cH)$ is {\it $\lambda$-bounded} if $\|V(s)\|\leq\lambda(s)$ for all $s\in S$.
\end{defn}

\begin{thm} \mlabel{OSreps}
Let $S =\Gamma_G(W)$ be a complex  Olshanski semigroup.
\begin{itemize}
\item[\rm(i)] Let $V:S\to\al B.(\al H.)$ be a holomorphic $*$-representation.
Then there is a unique continuous representation $U \: G \to \U(\cH)$
such that
\[ V(g \exp(iX)) = U(g)e^{i\partial U(X)}\quad \mbox{ for }
\quad g \in G, X \in W\,. \]
It satisfies $W \subeq B(I_U)$.
\item[\rm(ii)] A continuous unitary representation $(U,\cH)$ of $G$ extends to a
holomorphic $*$-representation $\hat{U} \: S \to  \cB(\cH)$ if and only if $W \subeq B(I_U)$.
For a closed convex $\Ad^*(G)$-invariant subset
$C \subeq \g^*$ with $W\subseteq B(C)$, the relation $I_U \subeq C$ holds
if and only if $\hat U$ is $\lambda_C$-bounded.
\end{itemize}
\end{thm}

Part (i), resp. (ii), is proven in \cite[Prop.~XI.3.1]{Ne99}, resp.~\cite[Thm~XI.3.3]{Ne99}.
So by (i) we have a bijection between the non-degenerate holomorphic representations of $S$
and a subclass of representations of $G$, and part (ii) identifies this subclass.
It also gives the condition on $\hat{U}$ which corresponds to the $C$-spectral condition
for $U$. We next construct the appropriate host algebra for these representations.
A generalization of this theorem to infinite dimensional Lie groups
and the corresponding Olshanski semigroup can be found in \cite[Thm.~4.9]{Ze15}
for locally convex groups and in \cite{MN12} for Banach--Lie groups.

\begin{defn}\mlabel{def:C-spocHost}
Let $S=\Gamma_G(W)$ be a complex Olshanski semigroup and
$C \subeq \g^*$ be a closed convex $\Ad^*(G)$-invariant subset
with $W\subseteq B(C)$.
 Let $\C[S]:=\C^{(S)}$ be
the space of finitely supported functions $f:S\to\C.$
This is the span of $\{\delta_s\mid s\in S\}$ where $\delta_s(t)=1$ if $t=s$ and zero otherwise.
Then $\C[S]$ is a $*$-algebra
with the product $\delta_s\cdot\delta_t:=\delta_{st}$ and involution
$f^*(s):=\overline{f(s^*)}.$ It has a submultiplicative norm
$\|\cdot\|_C$ given by
\[
\|f\|_C:=\sum_{s\in S}|f(s)|\,\lambda_C(s) \quad \mbox{ for } \quad f\in\C[S].\]
The completion of $\C[S]$
with respect to this norm is a Banach $*$-algebra,
which we denote by  $\cB := \ell^1(S,C)$.
Given any $C^*$-seminorm $p:\al B.\to[0,\infty)$, i.e.\
\[ p(ab) \leq p(a)p(b), \quad p(a^*) = p(a) \quad \mbox{ and }
\quad p(a^*a) = p(a)^2
\quad \mbox{ for }\quad a,b \in \cB,\]
let the $C^*$-algebra $\cB_p$ be  the completion of
$\cB/p^{-1}(0)$ with respect to $p$. Then there is a natural morphism
$\xi_p:S\to\cB_p$ by $\xi_p(s):=\delta_s+p^{-1}(0)\in \cB/p^{-1}(0)$. Define
\[
\cP_{\rm hol}:=\{p:\al B.\to[0,\infty)\,\mid\,p\;\hbox{is a $C^*$-seminorm,}\;\;
\xi_p:S\to\cB_p\;\;\hbox{is holomorphic}\}.
\]
Since every $C^*$-seminorm $p$ on $\ell^1(S,\C)$ satisfies
$p(f) \leq \|f\|_C$ for every $f$
(cf.~\cite[Cor.~III.1.21]{Ne99}), we obtain a well defined seminorm
\[ p_C(f) := \sup_{p \in \cP_{\rm hol}} p(f) \leq \|f\|_C \]
on $\cB=\ell^1(S,C)$ which is easily seen to be a $C^*$-seminorm
and since it corresponds to the holomorphic morphism
\[ \xi_C \: S \to \bigoplus^\infty_{p \in \cP_{\rm hol}} \cB_p, \quad
s \mapsto (\xi_p(s)), \]
it is the unique maximal element in $\cP_{\rm hol}$.
Define $C^*(S,C):=\cB_{p_C}$, i.e.\ the completion of  $\ell^1(S,C)/p_C^{-1}(0)$
with respect to $p_C$ and observe that it is isomorphic to the $C^*$-algebra
generated by the image of $\xi_C$. In particular, we
have the holomorphic $*$-morphism $\xi_C \: S \to C^*(S,C)$
whose range spans a dense subalgebra and which satisfies
\[ \|\xi_C(s)\| \leq \lambda_C(s) \quad \mbox {for } s \in S. \]
Since the absolute value $\lambda_C$ on $S$ is $G$-biinvariant,
the left and right actions of $G$ on $S$ lead to a
homomorphism $\zeta_G \: G \to \U(M(C^*(S,C)))$ determined by
 $\zeta_G(g)\cdot\sum\limits_{i=1}^na_i\,\delta_{s_i}:=\sum\limits_{i=1}^na_i\,\delta_{g\cdot s_i}.$
 Thus we obtain the triple ${(C^*(S,C),\xi_C,\zeta_G)}$, and in the remark below we will see that
 $(C^*(S,C),\zeta_G)$ is a host algebra for $G$.
\end{defn}

\begin{rem}\mlabel{UdefCorr}
(a) The triple ${(C^*(S,C),\xi_C,\zeta_G)}$ constructed above is universal
in the following sense. Let $\xi\: S \to \cA$
be any  holomorphic $*$-homomorphism into a $C^*$-algebra $\cA$
(i.e., it is multiplicative $\xi(ab) = \xi(a)\xi(b)$  and involutive, i.e.\ $\xi(s^*) = \xi(s)^*$)
which satisfies $\|\xi(s)\| \leq \lambda_C(s)$ for every $s \in S$. Then there exists a
unique morphism $\hat\xi\: C^*(S,C) \to \cA$ with
$\hat\xi\circ \xi_C = \xi$ (cf.\ \cite[Thm.~3.5]{Ne08}).
This universal property of the pair
$(C^*(S,C),\xi_C)$ determines it up to isomorphism
by standard arguments. 

(b) The definition of $C^*(S,C)$ from above seems to differ from the definition in
 \cite[Def.~IV.2.5]{Ne99}  of a
$C^*$-algebra $C^*(S,\lambda_C)$ which was based on the universal representation
specified by a set of holomorphic positive definite functions.
We claim that $C^*(S,C)\cong C^*(S,\lambda_C)$.
Now $C^*(S,\lambda_C)$  has a holomorphic
$*$-morphism $\xi_{\lambda_C} \: S \to C^*(S,\lambda_C)$
with total range,
and it has the universal property that, for every
${\lambda_C}$-bounded holomorphic representation
$V \: S \to \cB(\cH)$,
there is a unique representation
$\hat{V} \: C^*(S,{\lambda_C}) \to \cB(\cH)$ with
$\hat{V} \circ \xi_{\lambda_C} = V$
(\cite[Thm.~IV.2.7]{Ne99}).
Since every $C^*$-algebra $\cA$ has a faithful $*$-representation
$\cA \into \cB(\cH)$,
the pair $(C^*(S,{\lambda_C}), \xi_{\lambda_C})$ has the analogous universal
property in (a). That is,  for every holomorphic $*$-morphism
$\xi\: S \to \cA$ into a $C^*$-algebra with $\|\xi(s)\| \leq
{\lambda_C}(s)$ there exists a unique $C^*$-morphism
$\hat\xi \: C^*(S,{\lambda_C}) \to \cA$ with
$\hat\xi \circ \xi_{\lambda_C} = \xi$.
The universal property now implies the
existence of a unique isomorphism
$\Phi \: C^*(S,{\lambda_C}) \to C^*(S,C)$ with
$\Phi \circ \xi_{\lambda_C}  = \xi_C$.
\end{rem}

\begin{theorem}\mlabel{HostSpec}
Let $S=\Gamma_G(W)$ be a complex Olshanski semigroup and
$C \subeq \g^*$ be a closed convex $\Ad^*(G)$-invariant subset
with $W\subseteq B(C)$.
Then the following assertions hold:
\begin{itemize}
\item[\rm(i)] There is a surjective $*$-homomorphism $\gamma:C^*(G)\to C^*(S,C)$ such that
$\gamma\circ\eta_G(g)=\zeta_G(g)\circ\gamma$ for all $g\in G$, where
$\eta_G:G\to \U(M(C^*(G)))$ is defined by
\[
\big(\eta_G(g)f\big)(h):=
f(g^{-1}h)
\quad\hbox{ for }\quad g,h\in G,\,f\in L^1(G).
\]
\item[\rm(ii)] The pair $(C^*(S,C),\zeta_G)$ is a host algebra
for the $C$-spectral representations of $G$ in the sense that,
 for each complex Hilbert space
${\cal H}$, the map
\[ \zeta^* \: \Rep(C^*(S,C),{\cal H}) \to \Rep(G, {\cal H}), \quad\hbox{given by}\quad
\zeta^*(\pi):= \tilde\pi \circ \zeta_G, \]
is injective with range $\Rep_C(G,\cH)$.
\item[\rm(iii)] For any representation $\pi\in\Rep(C^*(S,C),{\cal H})$,
we have that
$\pi(C^*(S,C))=\zeta^*(\pi)\big(C^*(G)).$
Note that $\pi\circ \xi_C$ is the unique
holomorphic extension of $\zeta^*(\pi):G\to \U(\cH)$.
\end{itemize}
 \end{theorem}

 \begin{prf} (i) According to
\cite[Thm~XI.6.24]{Ne99}, the
surjective $*$-homomorphism $\gamma:C^*(G)\to C^*(S,C)$ is given by
\[
\gamma(f):=\lim_{n\to\infty}\int_Gf(h)\,\xi_C\big(\exp( \hbox{${i\over n}$}\,X)\,h\big)\,d\mu(h)
\quad\mbox{ for any }\quad X\in W\quad\mbox{and all}\quad f\in C_c^\infty(G),\]
where $\xi_C:S\to  C^*(S,C)$ is the canonical map induced by
$s \mapsto \delta_s$, and $\mu$ is a left Haar measure on $G$. Then
 \begin{eqnarray*}
\gamma(\eta_G(g)\,f)&=&\lim_{n\to\infty}\int_G(\eta_G(g)\,f)(h)\,\xi_C
\big(\exp\big(\textstyle{\frac{i}{n}X}\big)\,h\big)\,d\mu(h)\\
&=&\lim_{n\to\infty}\int_Gf(g^{-1}h)\,\xi_C\big(
\exp\big(\textstyle{\frac{i}{n}X}\big)\,h\big)\,d\mu(h) \\[1mm]
&=&\lim_{n\to\infty}\int_Gf(h)\,\xi_C
\big(\exp\big(\textstyle{\frac{i}{n}X}\big)\,gh\big)\,d\mu(h) \\[1mm]
&=&\lim_{n\to\infty}\int_Gf(h)\,\xi_C
\big[g\,
\exp\big(\textstyle{\frac{i}{n}\Ad(g^{-1})X}\big)\,h\big]\,d\mu(h)
=\zeta_G(g)\,\gamma(f)
\end{eqnarray*}
using the independence of $\gamma$ on $X\in W$.

(ii) This part follows directly from part (i) and \cite[Prop.~X.6.17]{Ne99}.

(iii) This part is proved in Prop.~XI.6.22 and the proof of
Thm.~XI.6.24 in \cite{Ne99}.
 \end{prf}

As $C^*$-algebras always have injective representations, this
theorem guarantees the existence
of  $C$-spectral representations of $G$ given the hypotheses for $W$ and $C$,
provided we know that $\zeta_G$ is injective. For pointed cones $W$
this is always the case by the Gelfand--Ra\"\i{}kov Theorems
(\cite[Thms.~XI.5.1/2]{Ne99}).

 \begin{rem}  \mlabel{rem:OlshHost}
(a) The surjective morphism $\gamma \: C^*(G) \to C^*(S,C)$ induces a
morphism
\[ \tilde\gamma \: M(C^*(G)) \to M(C^*(S,C))\]
 (\cite[Prop.~10.3]{Ne08}),
and the relation
$\gamma \circ \eta_G(g) = \zeta_G(g) \circ \gamma$ now implies that
\[ \tilde\gamma \circ \eta_G = \zeta_G \: G \to \U(M(C^*(S,C))).\]

(b) By Theorem~\ref{HostSpec}(i), this host algebra for the
$C$-spectral representations of $G$ could be obtained by
taking the quotient of the group algebra $C^*(G)$ by a suitable ideal.

(c) Construct the universal  $C$-spectral representation
$U_C:G\to \U(\cH_C)$ by a direct sum of the GNS representations
of all states on $C^*(G_d)$ whose GNS representations restrict on $G\subset C^*(G_d)$
to  $C$-spectral representations. This is
a $C$-spectral representation, hence it has a unique
holomorphic extension $\wt U_C$ to $S$.
By Theorem~\ref{HostSpec}(ii),
$\wt U_C$ must be the universal representation of
$C^*(S,C)$, hence $C^*(S,C)$ is isomorphic to the
$C^*$-algebra generated by
$\wt U_C(S)\subset \cB(\cH_C).$
By Theorem~\ref{HostSpec}(iii)
it then follows that
\[
C^*(S,C)\cong\zeta^*(\wt U_C)\big(C^*(G)),\qquad\hbox{hence}\qquad
\ker\gamma=\ker\zeta^*(\wt U_C).
\]
This gives an alternative definition of $C^*(S,C)$ as the image of $C^*(G)$
under the universal $C\hbox{-spectral}$ representation.
(An explicit characterization of $\ker\gamma$ is given
in \cite[Prop.~X.6.17]{Ne99}.)
 \end{rem}

\begin{defn}\mlabel{def:C-quotient}
Conversely, if one only assumes that  $C \subeq \g^*$ is a closed convex
$\Ad^*(G)$-invariant subset,
the question arises if there is a host algebra which is full for the
$C$-spectral representations. This can be obtained as
a quotient
\begin{equation}
  \label{eq:c-quotient}
C^*(G)_C := C^*(G)/\cI_C,
\end{equation}
where $\cI_C \trile C^*(G)$ is the common kernel of all
$C$-spectral representations.
\end{defn}

\begin{rem} (a) All $C$-spectral representations
of $G$ contain
$N := \oline{\la \exp C^\bot \ra},$ $C^\bot\subset\g,$ in their kernel,
so that $\cI_C$ contains the kernel of the quotient map
$C^*(G) \to C^*(G/N)$, and we may identify $C^*(G)_C$ with
$C^*(G/N)_{\tilde C}$, where
$\tilde C = C \cap \fn^\bot \subeq \fn^\bot \cong (\g/\fn)^*$.

(b) If $G$ is simply connected, then $\fn = C^\bot$ and we simply have
$\tilde C = C$, now considered as a subset of $(\g/\fn)^*$.
This procedure leads to a situation where $C$ spans $\g^*$, i.e.\ $C^\bot = \{0\}$.
Then $B(C)$ contains no affine lines, but it need not possess interior points.

(c) If $C \subeq \g^*$ is an invariant linear subspace, then $B(C) = C^\bot$
and $C \cong (\g/C^\bot)^*$. Therefore, if $G$ is simply connected,
 the $C$-spectral condition
becomes vacuous for the quotient group $G/N$, and $C^*(G)_C \cong C^*(G/N)$.

(d) Whether $C^*(G)_C$ can be written as some $C^*(S,C)$
depends on whether $B(C)$ has interior points. If this is the case, then
$W := B(C)^\circ/C^\bot$ is an elliptic cone in the quotient algebra $\g/C^\bot$
(\cite[Prop.~VII.3.4(c)]{Ne99}), so that we have  the corresponding
complex Olshanski semigroup $S= \Gamma_{G_1}(W)$ and $C^*(S,C)$ is a host algebra for
$G$ which is a quotient of $C^*(G/N)$, where $N \trile G$ is the normal subgroup
corresponding to the ideal $\fn = C^\bot \trile \g$.
This argument shows that we only need elliptic cones to obtain the semigroup.
\end{rem}

\begin{ex} Suppose that $S = \Gamma_G(W)$ is a complex Olshanski semigroup,
where $W \subeq \g$ is an open invariant convex cone.
Let $\hat U \: S \to \cB(\cH)$
be a non-degenerate holomorphic representation
and $U \: G \to \U(\cH)$ be the corresponding unitary representation of $G$.
We fix an element $X_\circ \in W$ and put $U_z := e^{z \partial U(X_\circ)}$ for $\Im z \geq 0$.
This makes sense because the boundedness of $\hat U$ implies that all operators
$-i\partial U(X)$, $X \in W$, are positive.
Then condition (S) in Theorem~\ref{thm:5.6} is satisfied because $\hat U$ is holomorphic and, for every
$z \in \C_+$, the map $G\to S, g \mapsto g \exp(z X_\circ)$, is real analytic, hence in
particular smooth. This produces a crossed product host which is a host algebra for $G$.
\end{ex}
\begin{rem}
In Example~\ref{ex:5.13} we already saw that using this idea for a semibounded representation
$U$ manages to produce a host algebra of the form $C^*(\pi(G) e^{i \partial U(W_U)})$ even in an
infinite dimensional setting. In our context, it clearly is a
homomorphic image of $C^*(S,I_U).$
\end{rem}
Above Proposition~\ref{ResLap} we saw that, for any representation $U:G\to \U(\cH)$
of a finite dimensional Lie group $G,$ the resolvent of the Laplacian is contained
in $U(C^*(G)).$ In general it is not true that the resolvent of
$-i \partial U(X)$ is in $U(C^*(G))$ for a non-zero $X\in\g$. However,
we now show that a
special property of $C\hbox{-spectral}$ representations is that, for any
$X$ in the open cone
\[ W_U = \{ X  \in \g \,\vert\, -i \partial U(X) \quad\hbox{is bounded below} \}^\circ
= B(I_U)^\circ,  \]
the resolvent $(\1-\partial U(X))^{-1}$ is in $U(C^*(G))$.

\begin{prop}\mlabel{ResCspec}
Let $G$ be a connected finite dimensional Lie group,
let $C \subeq \g^*$ be a closed convex $\Ad^*(G)$-invariant subset
and let $U:G\to \U(\cH)$ be a $C\hbox{-spectral}$ representation.
Then, for any non-zero $X\in B(C)^\circ$,
its resolvent  $(\1-\partial U(X))^{-1}$ is in $U(C^*(G))$.
Moreover, its multiplication action on $U(C^*(G))$ is non-degenerate.
\end{prop}

\begin{prf} Let $\oline U$ denote the corresponding representation
of the quotient group $\oline G := G/\ker U$.
Then $X \in W_U$ is equivalent to its image in $\oline\g$ being contained in
$W_{\oline U}$. We may therefore assume that $U$ is injective.
Then $W := W_U \supeq B(C)$ is an elliptic open invariant cone in $\g$
(\cite[Prop.~VII.3.4(c)]{Ne99})
and there exists a corresponding Olshanski semigroup
$S = \Gamma_G(W)$. As $U$ is a $C$-spectral representation,  we have
$\hat U(C^*(S,C))=U(C^*(G))$ by Remark~\ref{rem:OlshHost}(c).
Thus $\hat U(S) \subeq  U(C^*(S,C))=U(C^*(G))$ and by
Theorems~\ref{OSreps} and \ref{Wreps-OS}
 this is given by
\[ \hat U(g \exp(iX)) = U(g)e^{i\partial U(X)}\quad \mbox{ for }
\quad g \in G, X \in W\,. \]
Thus the one-parameter semigroup $(e^{it\partial U(X)})_{t>0}$
is contained
in $U(C^*(G))$ and hence the $C^*$-algebra generated by it is in $U(C^*(G)).$
The function $x\to e^{tx}$ is a $C_0$-function on ${\rm spec}(i\partial U(X))\subset
\R_-$ which separates the points, hence the $C^*$-algebra generated by $ e^{it\partial U(X)}$
consists of all $C_0$-functions of $i\partial U(X)<0$. In particular,
the resolvent $(\1- \partial U(X))^{-1}$ is in $U(C^*(G))$ as
$x\to (1-x)^{-1}$ is a $C_0\hbox{-function}$ on $\R_-.$

If $(\1-\partial U(X))^{-1}$  acts degenerately on  $U(C^*(G)),$  i.e.
 $(\1- \partial U(X))^{-1}L=0$  for some non-zero $L\in U(C^*(G)),$
 then the non-zero space $\ran(L)$ is contained in $\ker(\1- \partial U(X))^{-1},$
which contradicts  that $(\1- \partial U(X))^{-1}$ is a resolvent.
\end{prf}

\begin{cor}
\label{subgpalgcont}
Let $G$ be a connected finite dimensional Lie group,
$H \subeq G$ be a closed connected subgroup, and
$(U,\cH)$ be a semibounded unitary representation of $G$.
If $W_U \cap \fh \not=\eset$, then
$(U\res_H)(C^*(H)) \subeq U(C^*(G))$.
\end{cor}

\begin{prf} Passing to the quotient group $G/\ker U$, we
may w.l.o.g.\ assume that $U$ is injective. Here we use that we have for
every closed normal subgroup $N$ a natural surjection $C^*(G) \to C^*(G/N)$.
Then the momentum set $I_U$ generates the linear space $\g^*$
because $I_U^\bot = \ker \dd U = \{0\}$ is trivial.
Hence $W_U$ is an open elliptic cone.
By assumption, we may apply Proposition~\ref{ResCspec} with
$X \in \fh  \cap W_U$ to conclude that
$(\1 - \partial U(X))^{-1} \in U(C^*(G))$.
As $U(C^*(H))$ is contained in the multiplier algebra of $U(C^*(G))$  and
$(\1 - \partial U(X))^{-1}U(C^*(H))$ is dense in $U(C^*(H))$, the assertion
follows.
\end{prf}
For a closed subgroup $H \subset G$  of a locally compact group, in general
$C^*(H)\not\subeq C^*(G),$ just $C^*(H)\subeq M(C^*(G))$, so for many
representations  $(U,\cH)$ of $G$ we do not have  this property
of semibounded representations that $(U\res_H)(C^*(H)) \subeq U(C^*(G))$.

\begin{rem} The preceding two results can alternatively be obtained
without reference of Olshanski semigroups by using the more recent
technique of smoothing operators (cf.~Theorem~\ref{thmsbsmoothonegen}).

(a) To obtain Proposition~\ref{ResCspec}, we first observe that
the existence of an element in $B(C)^\circ = W_U$
implies that $U$ is semibounded.
For $X \in W_U$, let $U^X_t := U_{\exp tX}$.
Then Theorem~\ref{thmsbsmoothonegen} implies that
$U^X(C^\infty_c(\R))$ consists of smoothing operators,
so that $U^X(C^*(\R)) \subeq U(C^*(G))$ follows from Proposition~\ref{prop:5.8}.
In particular,
the resolvent $(\1- \partial U(X))^{-1}$ is in $U(C^*(G))$ as
$x\to (1-ix)^{-1}$ is a $C_0$-function on $\R$.

(b) For Corollary~\ref{subgpalgcont}, we can argue as follows:
Let $X \in W_U \cap \fh$
and consider the unitary one-parameter group $U^X_t := U_{\exp tX}$.
By Theorem~\ref{thmsbsmoothonegen}, $U^X$ and $U$ have the same space
of smooth vectors. As $X \in \fh$, the representation
$U\res_H$ also has the same space of smooth vectors.
Therefore $U\res_H(C^\infty_c(H))$ consists of smoothing operators
and the assertion follows from Proposition~\ref{prop:5.8}.
\end{rem}

\section{Spectral conditions for covariant representations}
\label{SpecCondCovRep}

In this section we study how
spectral conditions for non-abelian Lie groups relate to cross representations.

\subsection{Covariant $C$-spectral representations which are cross}
\label{CrossParam}

Having obtained the host algebra  $(C^*(S,C),\zeta_G)=(\cL,\eta)$
for the $C$-spectral representations of $G$ in the previous section, we can now apply
the general theory developed in \cite{GrN14}.

\begin{defn}\mlabel{def:1.1d2}
Let $(\cA, G, \alpha)$ be an automorphic $C^*$-action,
where $G$ is a connected finite dimensional Lie group,
and fix an
open invariant weakly elliptic
cone $  W \subeq \g$
with complex Olshanski semigroup
$S = \Gamma_G(W)$. Let
 $C \subeq \g^*$ be a closed
convex invariant subset with $W\subseteq B(C)$,
 and fix the host algebra $(\cL,\eta)$ with
$\cL = C^*(S,C)$.

(a) We say that a covariant representation
$(\pi, U)$ of $(\cA, G,\alpha)$ on $\cH$ satisfies the
{\it $C$-spectral condition} if the $G$-representation $U$ does, i.e.,
if $I_U \subeq C$ (cf.\ Definition~\ref{defSpec}).
Then we call $(\pi, U)$ a {\it $C$-spectral representation}
and put
\[
{\rm Rep}_C(\alpha,\cH):=\big\{
(\pi,U)\in{\rm Rep}(\alpha,\cH)\,\mid\, I_U \subeq C \big\}.
\]

(b) If $(\pi,U)\in{\rm Rep}_C(\alpha,\cH),$ then  its cyclic components
also satisfy the $C$-spectral condition. This allows us to define the {\it universal  $C$-spectral representation}
$(\pi_C,U_C)\in{\rm Rep}_C(\alpha,\cH_C)$ as follows.
Let $\fS_\alpha$
denote the set of those states
$\omega$ on $\cA\rtimes_{\alpha} G_d$ which produce
a covariant representation
$(\pi_\omega,U_\omega)\in{\rm Rep}(\alpha,\cH_\omega)$ and
put
\[ \fS_C:=\big\{\omega\in\fS_\alpha\,\mid\, (\pi_\omega,U_\omega)\in{\rm Rep}_C(\alpha,\cH_\omega)\big\}.\]
Then
\[
\pi_C:= \bigoplus_{\omega\in\fS_C}\pi_\omega,\quad U_C:= \bigoplus_{\omega\in\fS_C}U_\omega
\quad\hbox{on}\quad  \cH_C=\bigoplus_{\omega\in\fS_C}\cH_\omega.
\]
The  universal  $C$-spectral representation coincides with
the universal  covariant $\cL$-representation as in \cite[Def.~5.4]{GrN14}).

(c) A {\it $C$-spectral
crossed product for $\alpha$} is a
crossed product host algebra $\cC$ of $(\cA,G, \alpha)$, where
$\cL = C^*(S,C)$ is the  host algebra for the $C$-spectral representations
of $G$.
\end{defn}

Our main task will be to analyze the conditions for the existence of a non-zero
 $C$-spectral crossed product for a given $\alpha \: G \to \Aut(\cA)$.

\begin{rem} \mlabel{rem:1.2}
(a) For a $C$-spectral representation $(\pi, U)$, we have
$U(C^*(S,C))=U(C^*(G))$ by Remark~\ref{rem:OlshHost}(c), hence it is a cross representation
with respect to $C^*(S,C)$ if and only if it is a cross representation with respect to $C^*(G)$ (similar to the case in the
first paragraph of Section~\ref{RepPosCros}).

(b) For the trivial case $C = \g^*$, we have ${\rm Rep}_C(\alpha,\cH)={\rm Rep}(\alpha,\cH).$
If $\alpha$ is strongly continuous, then ${\rm Rep}(\alpha,\cH)$ is bijective to the set of representations of the crossed product
$\cA\rtimes_{\alpha} G$ on $\cH$ (cf.~\cite[Prop.~7.6.4]{Pe89}), hence
 ${\rm Rep}(\alpha,\cH)\not=\emptyset$ for some non-zero $\cH.$
Moreover $\cC := \cA \rtimes_\alpha G$ is a $C$-spectral
crossed product for $\alpha$.
 However, there are subsets $C\not=\{0\}$ for which
${\rm Rep}_C(\alpha,\cH)=\emptyset$ for all non-zero~$\cH$
(cf.~Example~\ref{CounterPos} below).

(c)  Assume that $C \subeq \g^*$ is a closed convex $\Ad^*(G)$-invariant subset
for which $B(C)$ has non-empty interior, i.e., $C$ contains no affine lines.
Recall that, if $C$ generates $\g^*$, i.e., $C^\bot = \{0\}$, then
the open cone $B(C)^\circ$ is elliptic (\cite[Prop.~VII.3.4(c)]{Ne99}).
For any elliptic open invariant cone $W \subeq B(C)$,
we have constructed the host algebra  $(C^*(S,C),\zeta_G)=(\cL,\eta)$
for the $C$-spectral representations of $G$, and hence the
covariant $C$-spectral representations are just the covariant $\cL$-representations, i.e.\
${\rm Rep}_C(\alpha,\cH)=\Rep_\cL(\alpha, \cH)$.
The structure of ${\rm Rep}_{W^\star}(\alpha,\cH)$  has been
analyzed for  $G=\R^n$  in \cite{Bo96}, and for  strongly continuous actions $\alpha$ in
\cite[Ch.~8]{Pe89} (using Arveson spectral subspaces).
\end{rem}

\begin{ex}
\mlabel{CounterPos}
Let $G = \R$, $\cA=C_0(\R)$, and consider the action
$\alpha \: G \to \Aut(\cA)$ by translations.
This is strongly continuous, and the crossed product
$\cA \rtimes_\alpha G$ is a transformation group algebra, well-known to be
isomorphic to $\cK(L^2(\R))$ (cf.~\cite[Thm~II.10.4.3]{Bla06}).
Thus each representation of $\cA \rtimes_\alpha G$  is a direct sum of the unique irreducible
representation, given by $(\pi(f)v)(x)= f(x)v(x)$ and $(U_tv)(x)=v(x-t)$ for all $f\in C_0(\R)$,
$v\in L^2(\R)$, $t,\,x\in\R$. All covariant representations of $\alpha$ are therefore direct sums
of this representation. However, the selfadjoint
generator of $U_\R$ is $i\frac{d}{dx}$, and this has spectrum $\R$.
Thus in every covariant representation of $\alpha$ the generator of the implementers has spectrum
$\R$. So if we choose $C=[0,\infty)$, then ${\rm Rep}_C(\alpha,\cH)=\emptyset$ for all
non-zero~$\cH$.
 \end{ex}

Example~\ref{CounterPos}
shows that, even if a full crossed product exists for the host algebra $C^*(G)$,
and there is a  host algebra $\cL$ which is full for the $C$-spectral representations of $G$, then
a non-zero $C$-spectral crossed product for $\alpha$ need not exist.

\begin{theorem}\mlabel{CrossSpecParam2}
Let $G$ be a connected Lie group, let $C \subeq \g^*$ be a closed
convex invariant subset with
$B(C)^\circ \not=\eset$. Fix
the host algebra $(\cL,\eta)$ with $\cL = C^*(G)_C$.
Assume we have an automorphic $C^*$-action $(\cA,G, \alpha)$
and a $C$-spectral representation $(\pi,U)\in{\rm Rep}_C(\alpha,\cH)$
and consider for $X \in \g$ the representation $U^X_t := U_{\exp tX}$ of $\R$
and the corresponding $\R$-action $\alpha_X$ on $\cA$.
Then the following are equivalent:
\begin{itemize}
\item[\rm(i)] For every $X \in B(C)^\circ$,
$(\pi,U^X)$ is cross for $\alpha_X$ with respect to $C^*(\R)$.
\item[\rm(ii)]  There exists an  element $X \in B(C)^\circ$
for which $(\pi,U^X)$ is cross for $\alpha_X$ with respect to $C^*(\R)$.
\item[\rm(iii)] $(\pi,U)$ is a cross representation of  $(\cA,G, \alpha)$
with respect to $\cL$ (or equivalently with respect to $C^*(G)$).
\end{itemize}
\end{theorem}

\begin{prf} (i) $\Rarrow$ (ii) is trivial.

(ii) $\Rarrow$ (iii): We use Proposition~\ref{ResCspec} from which we have
$R:=(\1-\partial U(X))^{-1}\in U(C^*(G))$.
By assumption, we have
$\pi(\cA) U_{\cL}(R) \subeq
\lbr U_{\cL}(R) \cB(\cH)\rbr= {U_{\cL}(\cL) \cB(\cH)}$
by Proposition~\ref{ResCspec}.
Note that ${U_{\cL}(\cL) \cB(\cH)}=\lbr U_{\cL}(R) \cB(\cH)\rbr$ by the fact that
multiplication by $R\in\cL$ is non-degenerate on $\cL$.
By normality of $ U_{\cL}(R)$ we obtain from the assumption that
 $\pi(\cA) U_{\cL}(C^*(R)) \subeq{U_{\cL}(\cL) \cB(\cH)}$. Let
$(E_n)_{n\in\N}\subset C^*(R)$
 be an approximate identity (countable as $ C^*(R)$ is separable),
 then $\|(E_n-\1)L\|\to 0$ for all $L\in\cL$ by nondegeneracy (cf. \cite[Thm.~A.2]{GrN14}). In particular
 \[
\pi(A) U_{\cL}(L) \in \lbr\pi(A) U_{\cL}(\{E_nL\,\vert\,n\in\N\})\rbr
\subset {U_{\cL}(\cL) \cB(\cH)}\quad\hbox{for all}\quad A\in\cA,\;
L\in\cL.
 \]
Thus ${(\pi,U)}$ is a cross representation.

(iii) $\Rarrow$ (i): If ${(\pi,U)}$ is a cross representation,
fix any $X \in B(C)^\circ$
and put $R:=(\1- \partial U(X))^{-1}.$
Then, by $R\in\cL$, it follows immediately that
$\pi(\cA) U_{\cL}(R) \subeq{U_{\cL}(\cL) \cB(\cH)}=\lbr U_{\cL}(R) \cB(\cH)\rbr$,
and now (i) follows from Theorem~\ref{thm:5.1}(b)(ii).
\end{prf}

This preservation of the cross property under restriction  to one-parameter subgroups
does not hold for general cross representations (cf.~\cite[Ex.~6.10]{GrN14}).
It requires the $C$-spectral condition.

\begin{rem} (1) If $X \in B(C)^\circ$
is such that the action $\alpha_X \: \R \to \Aut(\cA)$
is strongly continuous, then Example~\ref{ex:1.1} implies that
$(\pi, U^X)$ is cross for $\alpha_X$ with respect to $C^*(\R)$, so that
Condition~(i) in Theorem~\ref{CrossSpecParam2} is satisfied.

(2)  If the universal covariant $C^*(S,C)$-representation
$(\pi_u, U_u)$ of $(\cA, G, \alpha)$ satisfies the conditions
in Theorem~\ref{CrossSpecParam2} above, it is cross,
hence a crossed product host of $(\cA,G,\alpha)$ and $C^*(S,C)$
exists, which is full for the $C\hbox{-spectral}$ representations.
In this case  the
covariant $C$-spectral representations comprise the
representation theory of a $C^*$-algebra, hence we can analyze it with
standard $C^*$-tools.
\end{rem}

\begin{ex}
\mlabel{FockPosCross}
(a)  We return to Example~\ref{NoGap}, where
$(\cH,\sigma)$ consists of a non-zero complex Hilbert space
 $\cH$ and the symplectic form
$\sigma(x,y):={\rm Im}{\langle x,y\rangle}$.
The Weyl algebra $\cA=\ccr \cH,\sigma.$ carries an  action
$\alpha:\Sp(\cH,\sigma)\to\Aut(\cA)$ determined by $\alpha\s T.(\delta_x):=\delta_{T(x)}$.
Consider the  Fock representation $\pi_F:\cA\to \cB(\cF(\cH))$ (details and notation in
Example~\ref{NoGap}). This is covariant for the subgroup $\U(\cH)\subset \Sp(\cH,\sigma)$,
where the unitary implementers are the second quantized unitaries, i.e.
\[
\pi_F\big(\alpha\s U.(\delta_x)\big)=\Gamma(U)\pi_F(\delta_x)\Gamma(U)^{-1}
\quad \mbox{ for } \quad x \in \cH,\]
where
\[ \Gamma(U)\big(v_1\otimes_s\cdots\otimes_s v_n\big):=
Uv_1\otimes_s\cdots\otimes_s Uv_n \quad \mbox{ for } \quad
U\in\U(\cH),\ v_1,\ldots v_n\in\cH, \;n\in\N.\]
Consider now the case in Example~\ref{Lightcone}, where
we have a continuous unitary representation \break
${U:\R^4\to \U(\cH)}$ such that the
 joint spectrum of the generators for the  coordinate subgroups in $\R^4$
 is in the closed forward light cone
\[ \oline{V_+} =\big\{\b x.\in\R^4\mid x_0^2-x_1^2-x_2^2-x_3^2
 \geq 0,\; x_0\geq 0\big\}.\]
With $G=\R^4$  we  get $\g\cong\R^4\cong\g^*$
 and we take $W := V_+\subset\R^4\cong\g$
to be the open light cone,
and by commutativity of $G$, the cone $W$ is elliptic and invariant. As  $W$ is a convex cone,
we have $B(W)=W^\star=\oline{V_+}$. Then $B(V_+)=W^{\star\star}=\overline{W}\supset W$.
Thus we have a  $\oline{V_+}$-spectral condition for  $U:\R^4\to \U(\cH)$
(i.e. all the operators
$-i\partial U(X)$, $X \in W$, have non-negative spectrum)
and have satisfied the geometric hypotheses of Theorem~\ref{CrossSpecParam2}
by Proposition~\ref{prop:3.29}.

The second quantized representation
\[
\R^4\to  \U(\cF(\cH))\quad\hbox{by}\quad \b x.\mapsto\Gamma(U(\b x.))
\]
is still  a  $\oline{V_+}$-spectral representation, and it implements the action of $\R^4$
on $\cA$ by $\b x.\mapsto \alpha_{U(\b x.)}$ in the Fock representation.
If we fix any $X \in W$, then $-i\partial U(X)$  has non-negative spectrum, so
by Example~\ref{NoGap} the Fock representation is cross with respect to the host $C^*(\R)$
for the one-parameter automorphism group $t\mapsto \alpha_{U(tX)}$.
Then by  Theorem~\ref{CrossSpecParam2}, the Fock representation
$(\pi_F,\,\Gamma\circ U)$ is a cross representation with respect to $C^*(\R^4)$
for the full action $\alpha:\R^4\to{\rm Aut}(\cA)$.

(b) In general, we have not just a representation of $\R^4$, but a representation of the
restricted Poincar\'e group $\cP_0$ which is a semidirect product of the
proper orthochronous Lorentz group (i.e. the identity component) with $\R^4$,
where the joint spectrum of the representation on $\R^4$ is in the positive light cone.
In its larger Lie algebra, $V_+$ has empty interior, as it is contained
in a proper subspace, and there are no open invariant weakly elliptic
cones containing $V_+$ in the Lie algebra of $\cP_0$
(this follows from \cite[Thm.~VIII.3.6]{Ne99} because the Lie algebra
$\so_{1,3}(\R) \cong \fsl_2(\C)$ is not hermitian).
Hence we cannot satisfy the criteria of
Theorem~\ref{CrossSpecParam2}. However, the  Poincar\'e group is contained in the
conformal group $\SO_{2,4}(\R),$ which is a hermitian Lie
group and its Lie algebra contains a closed pointed generating invariant cone
$\tilde W$ with nonempty interior intersecting the translation Lie algebra
$\R^4$ precisely in the
closure of $V_+$ (\cite{HNO94}). Assuming that we have instead a
$\tilde{W}^\star\hbox{-representation}$   $U:\SO_{2,4}(\R)\to \U(\cH),$
then we can now proceed as above to conclude that for
any $X \in \tilde{W}^\circ$, as $-i\partial U(X)$  has non-negative spectrum. Then
by Example~\ref{NoGap} the
Fock representation is cross with respect to the host $C^*(\R)$
for the one-parameter automorphism group $t\mapsto \alpha_{U(tX)}$.
Then by  Theorem~\ref{CrossSpecParam2}, it follows that the Fock representation
$(\pi_F,\,\Gamma\circ U)$ is a cross representation with respect to
$C^*(\SO_{2,4}(\R))$
for the action ${\alpha: \SO_{2,4}(\R)\to \Aut(\cA)}$.
Presently we do not know if restriction also leads to a cross
representation for $C^*(\cP_0)$. This $C^*$-algebra is certainly
``more singular'' than $C^*(\SO_{2,4}(\R))$.
\end{ex}

For semibounded representations,
the cross property restricts to closed subgroups which are compatible with the associated open cone:
\begin{theorem}\mlabel{CrossSubgp}
Let $G$ be a connected finite dimensional Lie group, let
$H \subeq G$ be a closed connected subgroup, and
assume an automorphic $C^*$-action $(\cA,G, \alpha)$
and a covariant representation $(\pi,U)\in{\rm Rep}(\alpha,\cH)$
such that $(U,\cH)$ is a semibounded unitary representation of $G$ for which
 $W_U \cap \fh \not=\eset$. If  $(\pi,U)$ is a cross representation with respect to
 the host $C^*(G)$ then its restriction    $(\pi,U\res_H)$   to $H$ is a
 cross representation for $(\cA,H, \alpha\res_H)$  with respect to
 the host $C^*(H).$
\end{theorem}

\begin{prf} By Corollary~\ref{subgpalgcont} we have that $(U\res_H)(C^*(H)) \subeq U(C^*(G)).$
Thus, as  $(\pi,\, U)$ is a cross representation
with respect to $C^*(G),$ we have (abbreviating notation):
\[
\pi(\cA)U(C^*(H)) \subeq \pi(\cA)U(C^*(G))
\subseteq U(C^*(G)) \cB(\cH)=U(C^*(H)) \cB(\cH), \]
where the last equality follows from the fact that $C^*(H)$ acts non-degenerately
on $C^*(G)$ by multiplication. This means that
$(\pi,\, U\res_{H})$ is a cross representation with respect to $C^*(H)$
for the restricted action $\alpha:\,H\to{\rm Aut}(\cA)$.
\end{prf}
Without the semiboundedness of $U$ this restriction result fails
(cf.\ \cite[Ex.~6.10]{GrN14}).

\subsection{The connection of ${\cal C}$ with the crossed
product $\cA\rtimes_{\alpha} G$}
\label{crossprod1}

When a conventional crossed product exists, the question arises about
how it is connected with the crossed product host $\cC$ for $C$-spectral representations.
By Theorem~\ref{HostSpec}(iii), we suspect that $\cC$ is a factor algebra of
$\cA\rtimes_{\alpha} G$, and this is what we now show.

\begin{thm}
  \mlabel{thm:1.10}
Let $G$ be a connected Lie group and
let  $C \subeq \g^*$ be a closed
convex invariant subset with $B(C)^\circ \not=\eset$.
We fix the host algebra $(\cL,\eta)$ with $\cL = C^*(G)_C$.
Assume we have an automorphic $C^*$-action $(\cA,G, \alpha)$
which is strongly continuous and that $\fS_C\not=\emptyset.$ Then
 \[
 \cC\cong(\pi_C\times U_C)\big(\cA\rtimes_{\alpha} G\big)
=\lbr \pi_C(\cA)\cdot U_C(C^*(G)_C)\rbr.  \]
 Moreover, the structure maps $\eta_\cA \: \cA \to M(\cC)$ and
$\eta_G \: G \to \U(M(\cC))$ are the compositions with
$\pi_C\times U_C$ of those on $\cA\rtimes_{\alpha} G$.
\end{thm}

\begin{prf}
As $\alpha$ is strongly continuous,  its crossed product $\cA\rtimes_{\alpha} G$
is defined. In fact, the crossed product is characterized up to isomorphism by the fact that it is a crossed product host
 $({\cal L}_c, \eta_\cA,\eta_G)$ for the covariant representations
${\rm Rep}(\alpha,\cH)$ for each $\cH$,
and that it satisfies ${\cal L}_c=\lbr \eta_\cA(\cA)\cdot\eta_G(C^*(G))\rbr$
(cf.~\cite[Thm~2.6.1]{Wil07} and \cite{Rae88}).

Thus, given a representation $\pi\in\Rep({\cal L}_c,{\cal H})$, then for the associated covariant pair
$\eta^*(\pi)=\big(\tilde\pi \circ \eta_\cA, \tilde\pi \circ \eta_G\big) \in \Rep(\alpha, {\cal H})$,  we have
\[
\pi({\cal L}_c)=\lbr \big(\tilde\pi\circ\eta_\cA(\cA)\big)
\cdot\big(\tilde\pi\circ\eta_G(C^*(G))\big)\rbr
=C^*\!\Big(\big(\tilde\pi\circ\eta_\cA(\cA)\big)\cdot
\big(\tilde\pi\circ\eta_G(C^*(G))\big)\Big).
\]
As $\fS_C\not=\emptyset,$
 the universal  $C$-spectral representation
 $(\pi_C,U_C)\in{\rm Rep}_C(\alpha,\cH_C)$ is nontrivial.
 Then $\{0\}\not= \cC \cong C^*\!\Big(\pi_C(\cA)\cdot
U_C(C^*(G)_C)\Big)=\pi_\cC(\cC)$ which fixes
 the defining representation $\pi_\cC$ of ${\cal C}$.
We have $\eta^*(\pi_\cC)=\big(\pi_C, U_C\big) \in \Rep(\alpha, {\cal H}_C)$,
 and this pair defines the
representation $\pi_C\times U_C$
of $\cA\rtimes_{\alpha} G$, for which we again have that
$\eta^*(\pi_C\times U_C)=\big(\pi_C, U_C\big)$. Then
\begin{eqnarray*}
(\pi_C\times U_C)\big(\cA\rtimes_{\alpha} G\big)&=&\lbr \pi_C(\cA)\cdot U_C(C^*(G)_C)\rbr
=C^*\!\Big(\pi_C(\cA)\cdot U_C(C^*(G)_C)\Big)
=\pi_\cC(\cC),
\end{eqnarray*}
where we use that the representation $U_C$ of $C^*(G)$ factors through the
quotient $C^*(G)_C$. Thus, as $\pi_\cC$ is the defining representation  of ${\cal C}$,
it follows that ${\cal C}$ is a factor algebra of the crossed product
$\cA\rtimes_{\alpha} G$, and that the ideal factored out is $\ker(\pi_C\times U_C)$.
\end{prf}

\subsection{Generalizing the Borchers--Arveson Theorem to non-abelian groups}
\label{BATnonab}

Recall that for one-parameter $W^*$-dynamical systems, we have the important
Borchers--Arveson Theorem (cf.~\cite[Thm.~3.2.46]{BR02}):
\begin{thm} \mlabel{BA-thm} {\rm(Borchers--Arveson)}
Let $(\cM,\R, \alpha)$ be a $W^*$-dynamical system on a von Neumann algebra
$\cM\subseteq \cB(\cH)$. Then the following are equivalent.
\begin{itemize}
\item[\rm(i)] There is a positive strong operator continuous unitary one-parameter group
$U:\R\to\U(\cH)$ such that $\alpha_t=\Ad U_t$ on $\cM.$
\item[\rm(ii)] There is a positive strong operator continuous unitary one-parameter group
$U:\R\to\U(\cM)$ such that $\alpha_t=\Ad U_t$ on $\cM.$
\item[\rm(iii)]
For the Borel subset $S\subseteq\R$, let
$\cM^\alpha(S)$ denote the Arveson spectral subspace. Then
\[
\bigcap_{t\in\R}\br\cM^\alpha[t,\infty)\cH.=\{0\}.
\]
\end{itemize}
If these conditions hold, then we may take  $U:\R\to\cM$ to be
$U_t=\int_{\R}e^{-itx}dP(x),$ where $P$ is the projection-valued measure uniquely determined by
\[
P[t,\infty)\cH=\bigcap_{s<t}\br\cM^\alpha[s,\infty)\cH..
\]
\end{thm}
The inner one-parameter group $(U_t)_{t \in \R}\subset\cM$
given by the theorem is minimal in the following sense:
for any other unitary one-parameter group
$(\tilde U_t)_{t \in \R}$ with non-negative
spectrum implementing the same automorphisms, i.e.\
$\Ad(U_t) = \Ad(\tilde U_t)=\alpha_t$ for $t \in \R$,
the corresponding one-parameter group
$Z_t := \tilde U_{-t} U_t\in\cM'$ has non-negative spectrum.

Moreover, if there is a cyclic vector $\Omega\in\cH$ and a positive  one-parameter unitary group
$V:\R\to \U(\cH)$ such that $V_t\Omega=\Omega$ and $\Ad(U_t) = \Ad(V_t)=\alpha_t$ for $t \in \R$,
then $U_t=V_t$ for all $t \in \R$ (cf. \cite[Pro.~8.4.13]{Pe89}, \cite[Cor.~5.6]{BGN17}).

Our aim in this subsection is to generalize the  Borchers--Arveson Theorem
to the current (possibly noncommutative) setting, where $\R$ is replaced
by a Lie group~$G$. There is already a generalization
to finite dimensional abelian Lie groups (cf.~\cite{Bo66})
and there are also specific instances dealing
with conformal groups $\SO(2,d)$; see \cite{Koe02},
where the lifting is called the Borchers--Sugawara construction
(\cite{Su68}). This construction is also used in Conformal Field Theory
to  construct a projective representation of $\Diff(\bS^1)$ from a
positive energy representation of loop groups
(see \cite{KR87} for the algebraic side and \cite{GW84} for the analytic side).
In the context of infinite tensor products this kind of lifting problems
has been studied in \cite{St90}.

An invariant cyclic vector  provides the following generalization:

\begin{prop}
  \mlabel{BA-thmGen}
Let $G$ be a connected finite dimensional Lie group,
let  $C \subeq \g^*$ be a closed
convex invariant subset with $B(C)^\circ \not=\eset$, and
let $(\cA,G,\alpha)$ be an automorphic $C^*$-action.
Let $\omega$ be an $\alpha$-invariant state on $\cA$ such that
$(\pi_\omega,\,U^\omega)\in{\rm Rep}_C(\alpha,\cH_\omega)$ is a $C$-spectral representation
where $U^\omega : G\to \U(\cH_\omega)$ is the GNS unitary representation. Then
$U^\omega_G\subset\pi_\omega(\cA)'',$ and $U^\omega$ restricted to each one-parameter subgroup
$G_X:=\exp (\R X)\subset G$ for  $X \in B(C)^\circ$
is the Borchers--Arveson minimal group
for $t \to  \alpha_{\exp(tX)}.$
\end{prop}

\begin{prf}
The restriction of $U^\omega$ to the one-parameter subgroup
$G_X = \exp(\R X)$ for  $X \in B(C)^\circ$
 coincides with the minimal Borchers--Arveson
unitary group in $\pi_\omega(\cA)''$ which implements
$t \to  \alpha_{\exp(tX)}$,
as it remains the GNS unitary representation, even for
its restrictions (cf.~\cite[Cor.~5.6]{BGN17}).
As $B(C)^\circ \subeq \g$ is open, the union $\mathop{\bigcup}\limits_{X \in B(C)^\circ}G_X$
generates $G$ as a group and hence $U^\omega_G$, the group generated by
the subgroups $U^\omega_{G_X}$, is contained in $\pi_\omega(\cA)''$.
\end{prf}

Thus in the case that the $C$-spectral representation $(\pi_\omega,\,U^\omega)$ is the
GNS representation of an $\alpha\hbox{-invariant}$ state $\omega,$
a generalization of the Borchers--Arveson theorem holds. In the general case,
where we do not have such a state, the situation is more complicated.


\begin{thm}  {\rm(Non-abelian Borchers--Arveson Theorem)}  \mlabel{BA-nonab}
Let $G$ be a connected finite dimensional Lie group,
let  $C \subeq \g^*$ be a closed
convex invariant subset with $B(C)^\circ \not=\eset$, and
let $(\cA,G,\alpha)$ be an automorphic $C^*$-action.
For a $C$-spectral covariant representation
$(\pi,U)$ of $(\cA, G, \alpha)$, we consider the von Neumann algebra
\[ \cM := \pi(\cA)'' \quad \mbox{ and }\quad
Z := Z(\U(\cM)) = \U(\cM \cap \cM'),\]
where the abelian group $Z$ is endowed with the strong operator topology.
Then the following assertions hold:
\begin{itemize}
\item[\rm(i)] The representation $U \: G \to \U(\cH)$ is of inner type,
i.e., $U_G \subeq \U(\cM)\U(\cM')$. This is equivalent to
all automorphisms $\Ad(U_g) \in \Aut(\cM)$ being inner.
\item[\rm(ii)] The group
\[ \hat G := \{ (g, V) \in G \times \U(\cM) \mid V^{-1} U_g \in \U(\cM') \}\]
defines by $q \: \hat G \to G, (g,V) \mapsto g$
a central extension of $G$ by $Z$ which has continuous
local sections. If $Z$ is finite dimensional, i.e., if $\cM$ has finite
dimensional center, then $\hat G$ is a Lie group.
\item[\rm(iii)] By $\hat U_{(g,U)} := U$ and $\hat V_{(g,U)} := U^{-1} U_g$,
 we obtain continuous unitary representations
\[ \hat U \: \hat G \to \U(\cM), \quad
\hat V \: \hat G \to \U(\cM') \quad \mbox{ with }\quad
\hat U_{(g,U)} \hat V_{(g,U)} = U_g\quad \mbox{ for } \quad g \in G.\]
If $\cM$ is a factor, i.e., $Z = \T$, then both representations
are semibounded with $B(C)^\circ \subeq W_{\hat U} \cap W_{\hat V}$.
\item[\rm(iv)] Let $q_G \: \tilde G \to G$ denote the simply connected
covering of $G$ and suppose that $\cM$ is a factor.
If the second Lie algebra cohomology $H^2(\g,\R)$ vanishes, then
the central extension $\hat G$ of $G$ by $Z = \T$ splits over $\tilde G$
and there exist continuous homomorphisms
\[ \tilde U \: \tilde G \to \U(\cM), \quad
\tilde V \: \tilde G \to \U(\cM') \quad \mbox{ with }\quad
\tilde  U_g \tilde  V_g = U_{q_G(g)}\quad \mbox{ for }  \quad g \in \tilde G.\]
Then $\tilde U$ is unique up to a continuous homomorphism $\chi \: \tilde G \to \T$.
\end{itemize}
\end{thm}

This theorem generalizes the Borchers--Arveson Theorem in several
ways. First, in (i) it asserts the innerness of the automorphisms $\alpha_g$
for $\cM$,
and from (ii) we derive that they can be implemented by
unitary operators $\hat U_g$, unique up to a multiplicative factor
in the central group~$Z$, so that
$\hat U_g \hat U_h = f(g,h)\hat U_{gh}$, where
$f \: G \times  G \to Z$ is a group cocycle. Instead of dealing with cocycles,
we prefer to encode these ambiguous lifts in the central extension $\hat G$.
The minimality in the Borchers--Arveson Theorem is reflected in the second part of
(iii), asserting that both representations are semibounded.
The example of the projective oscillator representation of $\Sp_{2n}(\R)$
on the Fock space $\cF(L^2(\R^n))$
(cf.~Example~\ref{Sympcone}) and the fact that the ground state
energy of the oscillator Hamiltonian is non-zero, shows
that in (iv) we  cannot expect that $\tilde U$ provides
the minimal Borchers--Arveson implementation on one-parametergroups generated by
elements $X \in B(C)^\circ$. Presently we do not have a non-abelian version
of a minimality condition that specifies a unique lift.
Another issue that we leave open is the determination of
the extension group $\Ext(G,Z)$ in terms of Lie algebra cohomology
if $Z$ is infinite dimensional (see \cite{St90} for an interesting example of
this type).

\begin{prf} (i)
Pick a linear basis $X_1,\ldots, X_n$ of $\g$ contained in $B(C)$. Then
\cite{Bo66} implies the existence of unitary one-parameter groups
$\tilde U^{X_j} : \R \to \U(\cM)$ with
\[ \pi(\alpha_{\exp(tX_j)}(A)) = \tilde U^{X_j}_t \pi(A) \tilde U^{X_j}_{-t}
\quad \mbox{ for } \quad A\in \cA, t \in \R, j = 1,\ldots, n.\]
Since $G$ is connected, this implies that, for each
$g \in G$, there exists an element $\tilde U_g \in \U(\cM)$ with
\[ \pi(\alpha_g(A)) = \Ad(\tilde U_g)\pi(A)
= \tilde U_g \pi(A) \tilde U_g^*\quad \mbox{ for } \quad A \in \cA.\]
Then $\tilde V_g := \tilde U_g^{-1} U_g \in \cM' = \pi(\cA)'$ for each $g \in G$,
so that (i) follows from $U_g = \tilde U_g \tilde V_g$.

(ii) We consider the smooth map
\[ \Phi \: \R^n \to G, \quad
{\bf t} = (t_1, \ldots, t_n)
\mapsto \exp(t_1 X_1) \cdots \exp(t_n X_n) \]
and observe that, by the Inverse Function Theorem,
there exists a $0$-neighborhood $U \subeq \R^n$ for which
$\Phi\res_U \: U \to \Phi(U)$ is a diffeomorphism onto an open
neighborhood of the unit element of $G$. Here we use that the differential
in $0$ is the linear isomorphism
$\R^n \to \g$ specified by the chosen basis. As the map
\[ \tilde\Phi \: \R^n \to \U(\cM), \quad
{\bf t} = (t_1, \ldots, t_n)
\mapsto \tilde U^{X_1}_{t_1} \cdots \tilde U^{X_n}_{t_n} \]
is continuous, we obtain a continuous local section
\[ \sigma \: \Phi(U) \to \hat G, \quad
\Phi({\bf t}) \mapsto (\Phi({\bf t}), \tilde\Phi({\bf t})) \]
of the group extension~$\hat G$ of $G$ by $Z$,
hence it is a locally trivial topological fiber bundle.

If the center $\cM\cap \cM'$ is finite dimensional, i.e., isomorphic
to $\C^d$ with pointwise multiplication, then $Z \cong \T^d$ is a
torus group, and this implies that $\hat G$ actually is a finite dimensional
Lie group (\cite[\S VII.4]{Var07}).

(iii) The first assertion follows directly from the definition.
For the second assertion, let $X \in B(C)$ and
$m_X := \inf \Spec(-i\partial U(X))$. Then
$U_{\exp(tX)} e^{-it m_X}$ has a positive generator and
the Borchers--Arveson Theorem provides a minimal implementing unitary
group $\tilde U^X$ in $\cM$ with non-negative generator, for which
\[ \tilde V^X_t := \tilde U^X_{-t} U_{\exp(tX)} e^{-it m_X} \in \U(\cM') \]
also has a positive generator.

We write the Lie algebra of $\hat G$ as
$\hat\g = \g \oplus_\omega \R$ with the Lie bracket
\[ [(x,z), (x',z')] = ([x,x'], \omega(x,x')),\]
where $\omega \in Z^2(\g,\R)$ is a Lie algebra cocycle.
Then the lifts of $X \in \g$ to $\hat\g$ are of the form $\hat X := (X,c)$, $c \in \R$,
and the unitary one-parameter group $\hat U_{\exp(t \hat X)}$, which also
implements $\alpha_{\exp(tX)}$ on $\cM$, must have the form
$\tilde U^X_t e^{itm}$ for some $m \in \R$, where
$\tilde U^X$ is the minimal Borchers--Arveson implementation.
Here we use that $\cM$ is a factor, so that
the generators of the implementing one-parameter groups are unique up to an
additive constant.
We conclude that,
for $X \in B(C)$ and any lift $\hat X \in \hat\g$, both operators
$-i\partial \hat U(\hat X)$ and  $-i\partial \hat V(\hat X)$ are bounded from below.
This proves (iii).

(iv) If $H^2(\g,\R)$ vanishes, then \cite[Ex.~7.17]{Ne02} implies that
the pullback of the central $\T$-extension $\hat G$ of $G$
to the simply connected group $\tilde G$ splits, i.e., there exists a
continuous homomorphism $\sigma \: \tilde G \to \hat G$ with
$q \circ \sigma = q_G$. Then $\sigma(g) = (q_G(g), \tilde U_g)$,
where $\tilde U \: \tilde G \to \U(\cM)$ is a continuous
unitary representation. The first assertion now follows with
$\tilde V_g := \tilde U_g^{-1} U_{q(g)}$.
If $\hat U \:\tilde G \to \U(\cM)$ is another homomorphism splitting
$\hat G$, then $\chi(g) := \tilde U_g \hat U_g^{-1} \in Z$
defines a  homomorphism $\chi \: \tilde G \to Z$. Conversely,
any such homomorphism $\chi$ leads by $\hat U_g := \tilde U_g \chi(g)$
to another lift.
\end{prf}

The vanishing condition for $H^2(\g,\R)$ in Theorem~\ref{BA-nonab}(iv)
is satisfied for many Lie groups arising in physics,
such as the conformal Lie algebra $\so_{2,d}(\R)$ and the Poincar\'e Lie algebra
$\R^d \rtimes \so_{1,d-1}(\R)$. Hence a lift exists for the simply connected covering
of the corresponding groups. This is used in the version
of the Borchers--Sugawara Theorem obtained in \cite{Koe02}.

One could expect that the semiboundedness requirement implies that the
central extension $\hat G$ is trivial. The following example
shows that this is not the case, not even for factors of type~I.

\begin{ex} Let $\Osc := \Heis(\C) \rtimes_\beta \R$ denote the
four-dimensional oscillator group
\[ \Osc = (\T \times \C) \times \R, \qquad
(a,z,t)(b,w,s) = (ab e^{-i\Im(\oline z e^{it} w)/2}, z + w, t + s) \]
(cf.~Example~\ref{Sympcone}). We now form the groups
\[ \hat G := \Osc \times \Osc \quad \mbox{ and } \quad
G := \hat G/Z, \quad Z := \{ ((a,0,0),(a^{-1},0,0)) | a \in \T\}.\]
As $Z$ is contained in the commutator group of $\hat G$, the central $\T$-extension
$\hat G \to G$ is non-trivial.

For any pair $\bn = (n_1, n_2) \in \N_0^2$, there exists a semibounded representation
$(U^\bn, \cH^{\bn})$ of $\hat G$ with $U^\bn((a,0,0),(b,0,0)) = a^{n_1} b^{n_2}.$
It can be obtained on $\cF(\C)^{\otimes n_1} \otimes \cF(\C)^{\otimes n_2}$
from the Fock representation
of $\Osc$ on $\cF(\C)$. Now
$V := U^{(2,1)}$ and
$W := U^{(1,2)}$ have the property that
\[ V((a,0,0),(a^{-1},0,0)) = a^2 a^{-1} = a \quad \mbox{ and } \quad
 W((a,0,0),(a^{-1},0,0)) = a a^{-2} = a^{-1},\]
so that the tensor product representation
$(V \otimes W, \cH^{(2,1)} \otimes \cH^{(2,1)})$ factors through a
semibounded representation $(U, \cH^{(2,1)} \otimes \cH^{(2,1)})$
of the quotient group $G$.

Here
$\cM := \cB(\cH^{(2,1)}) \otimes \1$ is a von Neumann algebra with
$\cM' = \1 \otimes \cB(\cH^{(1,2)})$ and
$U_G \subeq \U(\cM)\U(\cM').$ The corresponding lifting problem
leads precisely to the non-trivial central extension $\hat G$.
\end{ex}

\begin{rem} Let $\cM \subeq \cB(\cH)$ be a von Neumann algebra,
$G$ a topological group and \break
$U \: G \to \U(\cH)$ be a continuous unitary representation of
inner type, i.e., $U_G \subeq \U(\cM)\U(\cM')$.
Above we discussed the problem to factorize $U$ into two representations
$\hat U \: G \to \U(\cM)$ and
$\hat W \: G \to \U(\cM')$.
Here we briefly discuss this problem on a general level.

(a) If $\cM$ is multiplicity free, i.e., $\cM'$ is abelian,
then $\cM' \subeq \cM'' = \cM$ shows that
$U_G \subeq \U(\cM)$ and we can put $\hat U_g := U_g$ and $\hat W_g = \1$.

(b) If $\cM$ is a von Neumann algebra in  standard form,
then we have a continuous representation of its automorphism
group $\Aut(\cM)$ on $\cH$ implementing the natural action on~$\cM$
(cf.\ \cite{BGN17}).
In particular, we have a unitary representation of the group
$G := \PU(\cM) := \U(\cM)/Z \subeq \Aut(\cM)$
of inner automorphisms.
In terms of an antiunitary involution $J$ on $\cH$ with $J\cM J = \cM'$ and
$JMJ = M^*$ for $M \in \cM \cap \cM'$,
the representation of $G$ is given by
\[ U \: G \to \U(\cH), \quad u Z\mapsto u J u J \quad \mbox{ for }
\quad u \in \U(\cM).\]

We identify the quotient $\PU(\cM)$ with
the group of inner automorphisms of $\cM$, but here
one has to be aware of the fact that this identification
is not compatible with the usual group topologies
(cf.\ \cite{Han77}, \cite{Co73})\begin{footnote}{The factor $\cM$ of type III
given by \cite[Cor.~1.5.8(c)]{Co73}
has the remarkable property that, for every faithful normal semifinite weight,
its modular automorphism group consists of inner automorphisms
and yet it is not implemented by any one-parameter unitary group in~$\cM$.
As explained in \cite[p.~21]{Ta83}, this happens only for
$W^*$-algebras with nonseparable predual. A similar example
can be found in \cite[p.~214]{St90}.}
\end{footnote}.
In this case $\hat G \cong \U(\cM)$ and a factorization of
the representation $U$ of $G$ exists if and only if the central extension
by~$Z$ splits.

As we shall see below, this is never the case when $\cM$ is a factor,
which corresponds to the case $Z = \T$, so that
$\U(\cM)$ is a central $\T$-extension of $\PU(\cM)$.
One indicator for the non-triviality of the central extension is
the intersection
\[ Z^1 := Z \cap (\U(\cM), \U(\cM)) \]
of $Z$ with the commutator group\begin{footnote}{Trivial central
extensions $Z \times G \to G, (z,g) \mapsto g$ are never perfect because
 all commutators are of the form $(e,g)$.}.
\end{footnote}of $\U(\cM)$.

(c) Consider the case when $\cM$ is of type I, i.e., $\cM \cong \cB(\cK)$ for some complex
Hilbert space~$\cK$.

For $\cM = M_n(\C)$, we have $\PU(\cM) = \PU_n(\C) = \PSU_n(\C)$, so that
$Z_1 = C_n = \{ z \1 \mid z^n =1 \}$ is the finite group of $n$th roots of unity.
To lift $\PU_n(\C)$ we therefore
have to pass to its $n$-fold simply connected covering group $\SU_n(\C)$.

If $\cK$ is infinite dimensional, then
the group $\PU(\cK) = \PU(\cM)$
is simply connected (as a consequence of Kuiper's Theorem),
but the central extension
$\U(\cK) \to \PU(\cK)$ is  non-trivial because the group
$\U(\cK)$ is perfect.
This follows from the fact that
every element of $\T \1$ is a commutator, which follows from
a simple construction with shift operators (see \cite{Hal67} for details).

If $\cM$ is of type II$_1$, then it has a trace, so that
the central extension $\fu(\cM) \to \pu(\cM)$ of Banach--Lie algebras splits
by $\su(\cM) := \{ X \in \fu(\cM) \,\mid\, \tr X = 0\}$. However,
for every $n \in \N$, there exists a morphism
$M_n(\C) \to \cM$ of von Neumann algebras, showing that
$C_n \subeq Z^1$
because $C_n \subeq (\U_n(\C), \U_n(\C)) = \SU_n(\C)$. This implies that
$Z^1$ is dense in $Z$. Accordingly, the central extension
$\U(\cM) \to \PU(\cM)$ is non-trivial.

If $\cM$ is not of finite type and $\cH$ is separable,
then \cite[Prop.~III.1.3.6]{Bla06} implies the existence
of a unital morphism $\cB(\ell^2) \to \cM$ of von Neumann algebras
and we obtain the non-triviality of the central extension
$\U(\cM) \to \PU(\cM)$ from the corresponding assertion for $\cB(\ell^2)$.
\end{rem}

\section{Examples}
\mlabel{sec:8}

\subsection{A full crossed product host with a spectrum condition}
\mlabel{subsec:8.1}

We want to give an example of a full crossed product host
when the usual crossed product does not exist, in particular
for a discontinuous  action  $\alpha:G\to\Aut\cA$.
Unfortunately, our previous example  \cite[Ex.~5.9]{GrN14} contains an error,
 hence fails.
We will not do this for $\cL=C^*(G)$, but choose instead a host $\cL$
which produces a spectral restriction.

 Start with Example 6.12 in \cite{GrN14}, and for simplicity
 take the smallest nontrivial symplectic space $(X,\sigma) = (\R^2,\sigma)$
with $\sigma((p,q),(p',q')) = pq'-p'q$.
Let $G := \Sp(X,\sigma)$ denote the corresponding symplectic group,
and consider the Weyl algebra $\cA=\ccr X,\sigma.$ with the (discontinuous)  action  $\alpha:G\to\Aut\cA$
 by Bogoliubov transformations $\alpha_g(\delta_x):=\delta_{g(x)}$.
Let  ${(\pi_F,V)}\in\Rep_\cL(\alpha,{\cal H}_F)$ with $\cL={C^*(G)}$ be the Fock representation with the usual
(second quantized) unitary
implementers for $\alpha$. By \cite[Ex.~6.12]{GrN14}, this is a cross representation.
By the Stone--von Neumann Theorem, any
regular representation of $\cA$ is a multiple of the Fock representation,
hence a cross representation if we take the same unitary implementers $V_g$
for $\alpha_g$ in each summand. Denote these implementers by $\tilde{V}_g$.

Now  $G = \Sp(X,\sigma)$ contains the scale transformations
$(q,p) \to (e^tq, e^{-t}p)$ for $t  \in   \R$,
hence any covariant representation ${(\pi,U)}\in\Rep_\cL(\alpha,{\cal H})$
must be a regular representation of $\cA,$ hence  ${(\pi,\tilde{V})}$
is a cross representation. A full crossed product host exists
if and only if each corresponding representation ${(\pi,U)}$ is a cross representation.
We can only prove below that ${(\pi,U)}$ is cross, if it is semibounded
with respect to a certain cone,
which implies that we obtain a crossed product host which is full for this
class of semibounded representations.

Observe that the Fock representation  ${(\pi_F,V)}$ is
 a positive representation for e.g. the harmonic oscillator evolution $t \to r(t) \in  \Sp (X, \sigma)$
with Hamiltonian  $P^2+Q^2.$ Now $r(t)$ has generator contained in the
open invariant elliptic cone $W\subset \sp (X, \sigma)$ defined in
Example~\ref{Sympcone},
and in fact the Fock representation  ${(\pi_F,V)}$ is positive on $W,$
i.e. $t\mapsto V_{\exp(tA)}$ is positive for all $A\in W.$
Thus we consider the $W^\star$-representation, which produces
the host $\cL_{W^\star}$ associated with this cone.
 We want to show that there is a full crossed product host for
 $(\alpha,\cL_{W^\star})$, i.e.\ that every $\cL_{W^\star}\hbox{-covariant}$ representation of
 $\alpha$ is cross (cf.~Theorem~\ref{thm:exist}).
Consider an $\cL_{W^\star}$-covariant representation $(\pi, U).$
Then $\pi$ is regular, hence ${(\pi,\tilde{V})}$ is cross, and it is also an
$\cL_{W^\star}$-covariant representation. As the Fock implementers leave
the Fock vacuum invariant, and this is cyclic, it follows from Proposition~\ref{BA-thmGen}
that ${V}_G\subset\pi_F(\cA)''$ and that for any one-parameter group
$t\mapsto\exp(tA)\in G,$ $A\in W,$ the unitary group $t\mapsto V_{\exp(tA)}$ coincides with
the Borchers--Arveson minimal group for $t\mapsto\alpha_{\exp(tA)}.$
As $\pi$ is a multiple of $\pi_F,$ we can define a normal map
$\tilde\pi:\pi_F(\cA)''\to \pi(\cA)''$ by $\tilde\pi(\pi_F(A)):=\pi(A)$ for $A\in\cA$
and extending to $\pi_F(\cA)''$ (in other words, we just take the appropriate multiple
of $\pi_F$ to produce $\pi$). As ${V}_G\subset\pi_F(\cA)''$ we have that
$\tilde\pi(V_G)=\tilde{V}_G.$ However, normal maps take minimal Borchers--Arveson groups
to minimal Borchers--Arveson groups (cf. \cite[Lemma~4.19]{BGN17}),
hence  $t\mapsto \tilde{V}_{\exp(tA)}$ coincides with
the Borchers--Arveson minimal group for $t\mapsto\alpha_{\exp(tA)},$ $A\in W.$
As $(\pi, U)$ is an $\cL_{W^\star}$-covariant representation,
 $t\mapsto U_{\exp(tA)}$ is a positive unitary group for any $A\in W.$
 Thus, using the fact that  ${(\pi,\tilde{V})}$ is cross with respect to $\cL_{W^\star}$
 (and that its restriction to $t\mapsto\exp(tA)$ is cross with respect to $C^*(\R)$),
 it follows from Theorem~\ref{CrossRepsSame}(ii) that
 $(\pi, U)$ is cross with respect to $C^*(\R)$ for each subgroup $t\mapsto\alpha_{\exp(tA)},$ $A\in W.$
Thus, by Theorem~\ref{CrossSpecParam2}, $(\pi, U)$ is cross with respect to
 $\cL_{W^\star}$ for the whole action $\alpha:G\to\Aut\cA$.   \\[3mm]

\begin{rem} We connect this discussion with the semigroup picture. We have an
Olshanski semigroup
\[ S_1 = \Gamma_{G_1}(W_1) = \Heis_{2n}(\C) \rtimes \Gamma_{G}(W),\]
where $G = \Mp_{2n}(\R)$ and $G_1 = \Heis_{2n}(\R) \rtimes \Mp_{2n}(\R)$, and
$W \subeq \sp_{2n}(\R)$ is the canonical open invariant cone.
We know that our representation  extends holomorphically
to a representation  $\hat U$ of $S_1$ (Theorem~\ref{OSreps}).
Then $\cA \cL_C$ is generated by
$\hat U(S_1)$,  and these operators are left
continuous for the action of $G$. Hence we have a cross representation.
\end{rem}


 \subsection{The translation action on $\cA = C_b(\R)$}
 \label{transCb}

A fundamental example is that of the translation action of $\R$. This is strongly continuous on
$C_0(\R)$, hence every covariant representation is cross. On the other hand the translation action
on $C_b(\R)$ is discontinuous, so a very natural question is whether this is cross for standard
representations, e.g.\ on $L^2(\R)$. This is the question we now answer.
\begin{prop}
 For the translation action $\alpha$ of $\R$ on $\cA = C_b(\R)
\subeq \cB(L^2(\R))$, the covariant representation of $(\cA,\R,\alpha)$ on $L^2(\R)$ is
not a cross-representation with respect to  the host $\cL=C^*(\R)$.
\end{prop}

\begin{prf}  For $\cL=C^*(\R)$,
we start by determining what is in $\cA_{\cL}$.
To fix notation, put
\[ (\alpha_tf)(x):=f(x+t), \qquad \pi:\cA\to  \cB(L^2(\R)),
\quad (\pi(f)\psi)(x):=f(x)\psi(x)\]
and
\[ U:\R\to U(L^2(\R)), \quad (U_t\psi)(x):=\psi(x+t).\]
 Then $U_t:=\exp(itP)$, where $P=-i\frac{d}{dx}$ on $C^\infty_c(\R)$. Denote its
 spectral measure by $E$.
As usual, $Q$ denotes the multiplication operator by $x$.

 As $\alpha$ is strongly continuous on $C_0(\R)$, we have $C_0(\R)\subset\cA_{\cL}$.
 Moreover, the functions $e_r(x):=e^{ixr},$ $r\in\R$,
 are eigenvectors of $P$ with eigenvalue $r$, hence they shift the spaces
 $E([a,b])$ to $E([a+r,b+r])$ and hence the matrix for $\pi(A)=\pi(e_r)=e^{irQ}$ below Lemma~\ref{Hdecomp}
 will have $c_0\hbox{-rows}$ and columns, showing that  $e_r\in\cA_{\cL}$.
To verify the defining condition for $\cA_{\cL}$  explicitly,  pick
$L \in \cL$ and write  $U_\cL(L)=f(P)$ for some $f\in C_0(\R)$. Then
\[ \pi(A) U_\cL(L)=e^{irQ}f(P)=f(P-r\1)e^{irQ}
\in \{h(P)\mid h\in C_0(\R)\}\cdot e^{irQ} \subset U_\cL(\cL) \cB(L^2(\R)) \]
which shows explicitly that $e_r\in\cA_{\cL}$.
Hence $\cA_{\cL}$ contains all almost periodic functions.

We now show that if $\pi(A)=e^{iQ^3}$, then $A\not\in\cA_{\cL}$
(a similar argument works for $e^{iQ^n}$ where $n\geq 3$). By Corollary~\ref{contSpace}(iii)
it suffices to show that
\[
\lim_{t\to 0}\left\| e^{itP}\pi(A)E[1,2]  -\pi(A)E[1,2]  \right\|\not=0\,,
\]
 which we will now prove. Define $\psi_k\in E[1,2]L^2(\R)$ by
 \[
\psi_0(x):=\sqrt{2\pi}(\cF^{-1} \chi\s[1,2].)(x)=\int_1^2 e^{ipx}dp
=\frac{e^{ix}}{ix}(e^{ix}-1),\qquad
\psi_k(x):=\psi_0(x+k)
\]
i.e.\ $\psi_k=e^{ikP}\psi_0$, hence  $\|\psi_k\|=\sqrt{2\pi}$ for all $k$.
Though $\lim\limits_{t\to 0}\left\|\big( e^{itP}-\1\big)\pi(A)\psi_k  \right\|=0$
for all $k$, we will show that
\[
\lim_{t\to 0}\sup_k\left\|\big( e^{itP}-\1\big)\pi(A)\psi_k  \right\|\not=0,
\]
which suffices for our proof. 
\begin{eqnarray*}
&&\!\!\!\!\!\!\!\!\!\!\!\!
\left\|\big( e^{-itP}-\1\big)\pi(A)\psi_k  \right\|^2=\int_\R\left| (\pi(A)\psi_k)(x-t)- (\pi(A)\psi_k)(x)\right|^2dx\\[1mm]
&=&\int_\R\left| \frac{e^{i(x-t)^3 }}{x+k-t}e^{-i(x+k-t)}\big(e^{-i(x+k-t)}-1\big)
-\frac{e^{ix^3 }}{x+k}e^{-i(x+k)}\big(e^{-i(x+k)}-1\big)
\right|^2dx   \\[1mm]
&=&\int_\R\left| e^{i(-3x^2t+3xt^2-t^3+t) }\cdot\frac{ e^{-i(x+k-t)}-1  }{x+k-t}
-\frac{ e^{-i(x+k)}-1  }{x+k}\right|^2dx   \\[1mm]
&=&\left\|\big( Ze^{itP}-\1\big)\varphi_k  \right\|^2
\end{eqnarray*}
where $Z:=\exp i\xi(Q)$ with $\xi(x):=-3x^2t+3xt^2-t^3+t$, and $\varphi_k:=e^{-ikP}\varphi_0$
with $\varphi_0(x):=(e^{-ix}-1)\big/x$. Now
\begin{equation}
  \label{eq:8.1}
\left\|\big( e^{itP}-\1\big)\pi(A)\psi_k  \right\|=\left\|\big( Ze^{itP}-\1\big)\varphi_k  \right\|
\geq\big|\| Z(e^{itP}-\1)\varphi_k  \|-\|Z\varphi_k-\varphi_k\|                 \big|.
\end{equation}
Note that $\varphi_k\in L^2(\R)$, so  for each $\varepsilon>0$ there is a $t_\varepsilon>0$
such that $|t|<t_\varepsilon$ implies that for all $k$ we have
${\|(e^{itP}-\1)\varphi_k\|}={\|(e^{itP}-\1)\varphi_0\|}<\varepsilon$.
For $0 < |t| < t_\eps$ we thus obtain
$\| Z(e^{itP}-\1)\varphi_k  \|\leq\|Z\|\varepsilon=\varepsilon$. Regarding the second term in \eqref{eq:8.1}, observe
\[
\|Z\varphi_k-\varphi_k\|^2=\int_\R\left| e^{i\xi(x) }-1\right|^2
\left|\frac{ e^{-i(x+k)}-1  }{x+k}\right|^2dx,
\quad \mbox{ and } \quad |\varphi_0(x)|=\left|\frac{\sin(x/2)}{x/2}\right|\]
 has a maximum of $1$ at $x=0$, hence the translate
$|\varphi_k(x)|$ has a maximum of $1$ at $x=-k$. By continuity we have an interval $[-L,L]$ (resp. $[-k-L,-k+L]$)
such that $|\varphi_0(x)|>\hlf$ for $|x|\leq L$ (resp. $|\varphi_k(x)|>\hlf$ for $|x+k|\leq L$),
and we assume $\varepsilon^2<L$. So
\begin{eqnarray*}
\|Z\varphi_k-\varphi_k\|^2&\geq&\int_{-k-L}^{-k+L}\left| e^{i\xi(x) }-1\right|^2
|\varphi_k(x)|^2dx
\geq\frac{1}{4}\int_{-k-L}^{-k+L}\left| e^{i\xi(x) }-1\right|^2dx
\\[1mm]
&=&L-\hlf\int_{-k-L}^{-k+L}\cos\xi(x) dx
=L-\hlf\int_{-L}^{L}\cos(\xi(y-k))\, dy\,.
\end{eqnarray*}
Now
\[ \int_{-L}^{L}\cos(\xi(y-k))\, dy=
{\rm Re}\int_{-L}^{L}e^{i6ykt}\exp i(-3y^2t+3yt^2)\,dy
\cdot e^{i(-3k^2t-3kt^2-t^3+t)}.\]
By the Riemann--Lebesgue Lemma,
the last integral goes to zero as
$k\to\infty$, hence
\[
\limsup\limits_{k\to\infty}\|Z\varphi_k-\varphi_k\|^2\geq L\quad \hbox{and so}\quad
\sup\limits_{k}\|Z\varphi_k-\varphi_k\|\geq \sqrt{L}>\varepsilon.
\]
 We conclude
from \eqref{eq:8.1} that
\[
\sup_k\left\|\big( e^{itP}-\1\big)\pi(A)\psi_k  \right\|\geq\sqrt{L}-\varepsilon
\]
from which it follows that
$\lim\limits_{t\to 0}\sup\limits_k\left\|\big( e^{itP}-\1\big)\pi(A)\psi_k  \right\|\not=0\,.$
Thus $A \not\in \cA_{\cL}$,
so that $(\pi,U)$ is not a cross representation.
\end{prf}

\begin{rem} One would like to characterize $\cA_{\cL}$  precisely.
 Proving $e^{iQ^2}\not\in\cA_{\cL}$ is possible, but needs a different argument:
 \[
U_\cL(L) \pi(A)e^{itP} =f(P)e^{iQ^2}e^{itP}=f(P)e^{it(P+2Q)}e^{iQ^2} =f(P)e^{itP}e^{2it}e^{iQ^2} e^{i2tQ}
 \]
 so as $f(P)e^{itP}e^{2it}$ converges in norm to $f(P)$ as $t\to 0$, for small $t$ the displayed expression is arbitrarily
 close to $f(P)e^{iQ^2} e^{i2tQ}$. The latter is discontinuous in norm with respect to\ $t$ which we can show by
 evaluating it on a suitable sequence of unit vectors $\psi_n$. Thus $e^{iQ^2}\not\in\cA_{\cL}$.
 \end{rem}

 \subsection{Host algebras for more infinite dimensional groups}
 \label{CPHinfdim}

Above, in Section~\ref{subsec:smop} we analyzed $\R$-actions on Lie groups (equivalently,
on their discrete group algebras) and via smoothing operators in positive  covariant representations,
obtained crossed product hosts which provide hosts for the associated semidirect product groups
$G \rtimes_\alpha \R$. Here we want to do two more infinite dimensional examples.
 \begin{ex} \label{VirExmp}
 (Positive energy representations of the Virasoro group)
Let $\Vir$ denote the simply connected Virasoro group.
This is a central extension of the group $\Diff_+(\bS^1)$ by $\R$.
On the Lie algebra level we have
\[ \vir = \R \bc \oplus_\omega \cV(\bS^1), \]
where $\omega \: \cV(\bS^1) \times \cV(\bS^1) \to \R$ is a certain
$2$-cocycle; see \cite[\S 8.2]{Ne10} for details.

Let $\bd = (0,\partial_\theta) \in \vir$ denote the element
corresponding to the generator $\partial_\theta \in \cV(\bS^1)$
of the rigid rotations of $\bS^1 \cong \R/2\pi \Z$. Then
$\ft := \R \bc \oplus \R \bd \subeq \vir$
is a $2$-dimensional abelian subalgebra corresponding to a subgroup
$T := \exp(\ft) \cong \R^2$ in $\Vir$.

There is a surprisingly simple classification of open invariant cones
$W \subeq \vir$ (\cite[Thm.~8.15]{Ne10}). They are determined by
the intersection $C := W \cap \ft$, and the arising cones $C$ are, up to sign,
precisely the ones with
\[ \R^\times_+ \bc + \R^\times_+ \bd \subeq C \subeq
\R \bc + \R^\times_+ \bd.\]

We say that a continuous unitary representation $(U,\cH)$ of $\Vir$
is a {\it positive energy representation} if
$-i\partial U(\bd) \geq 0$. If $(U,\cH)$ is semibounded
and $W_U \subeq \vir$ is the corresponding open invariant cone, then
$\bd \in \oline{W_U} \cup -\oline{W_U}$ implies that
either $U$ or the dual representation $U^*$ satisfies
the positive energy condition.
Conversely, every positive energy representation
is smooth  and semibounded by \cite[Thm.~2.18]{Ze17}.

We also note that, for a semibounded representation $(U,\cH)$,
the element $\bd$ need not be contained in $W_U$,
but either $\bd + \bc$ or $-(\bd + \bc)$ is.
We are therefore in the situation of Example~\ref{ex:5.13},
which shows that the one-parameter group $T_1 := \exp(\R(\bd + \bc)) \subeq T$
leads to the host algebra
\[ \cC := C^*(U_{\Vir} U(C^*(T_1))) \]
for $\Vir$. This algebra is a crossed product host for
the discrete group algebra $\cA = C^*(\Vir_d)$ and the $\R$-action
on $\Vir$ given by
\[ \alpha_t(g)
= \exp(t(\bd + \bc)) g \exp(-t(\bd + \bc))= \exp(t\bd) g \exp(-t\bd),\]
which correspond to the canonical action of the rigid rotations of $\bS^1$ on~$\Vir$.

These host algebras can be used to obtained direct integral decompositions
of semibounded, resp., positive energy representations into irreducible ones.
This provides an alternative to the complex analytic methods used in
\cite{NS14} for the same purpose.
Since $G$ does not carry an analytic Lie group structure, it would be meaningless to strengthen (S) in Theorem~\ref{thm:5.6} to an analytic version.
\end{ex}

\begin{ex}
\label{TwLoopExmp}
(Positive energy representations of twisted loop groups)
Let $K$ be a simple simply connected compact Lie group,
$\phi \in \Aut(K)$ be a finite order automorphism,
$\phi^N = \id_K$, and let
\[ \cL_\phi(K) := \{ \xi \in C^\infty(\R, K) \,\mid\, (\forall t \in \R)\, \xi(t+1)
= \phi^{-1}(\xi(t))\} \]
with pointwise multiplication be the corresponding twisted loop group
(\cite{PS86}).
Then
$(\alpha_t f)(x) := f(x-t)$ defines a smooth action of $\R$ on $\cL_\phi(K)$
which factors through an action of the circle group $\R/\Z N$.
The twisted loop group has a natural simply connected central extension
$G$ by the circle group $\T$. Accordingly,
\[ \g \cong \R \bc \oplus_\omega \cL_\phi(\fk) \quad \mbox{ with } \quad
\cL_\phi(\fk) = \{ \xi \in C^\infty(\R, \fk) \,\mid\, (\forall t \in \R)\, \xi(t+1)
= \phi^{-1}(\xi(t))\}, \]
$\omega(\xi,\eta) = \frac{1}{2\pi} \int_0^1 \kappa(\xi'(t), \eta(t))\, dt$
is the corresponding Lie algebra cocycle,
$\kappa$ is a suitably normalized Cartan--Killing form on $\fk$,
and $\alpha$ naturally lifts to $G$, which leads to a Lie group
$G^\sharp = G \rtimes_\alpha \R$.

Let $\fm := \fk^\phi \subeq \fk$ be the fixed point algebra of $\phi$ in $\fk$.
Then the fixed point algebra
$\g^\alpha = \R \bc \oplus \fm$ is a finite dimensional compact Lie algebra,
and if $\ft_\fm \subeq \fm$ is maximal abelian, then
$\ft = \R \bc \oplus \ft_\fm \subeq \g$ is a finite dimensional abelian Lie algebra
in $\g$ and
\[ \ft^\sharp = \R \bc \oplus \ft_\fm \oplus \R \bd, \quad \bd = (0,0,1), \]
 is a finite dimensional abelian Lie algebra in the Lie algebra
$\g^\sharp$ of $G^\sharp$. Averaging over the torus group
$\Ad(T^\sharp)$ shows that every open invariant cone in $\g^\sharp$ intersects
$\ft^\sharp$ non-trivially (cf.\ the proof of Theorem~5.2 in \cite{Ne14b}).

If $(U,\cH)$ is a semibounded representation of $G^\sharp$,
it follows that condition (S) in Theorem~\ref{thm:5.6} is satisfied for every $X \in \ft^\sharp \cap W_U$,
and this applies to all irreducible projective positive energy representations
of $(\cL_\phi(K),\R, \alpha)$. For more details and refinements
concerning the untwisted case we refer to
\cite[\S 2.3]{Ze17} and \cite{JN20}.
\end{ex}

\section{Conclusions and Discussion}

We do a brief review of the work above, and list open problems and
 possible future directions of enquiry.
We extended our earlier work on crossed product hosts in \cite{GrN14} first by
showing stability of the cross property under perturbation by a natural class of operators,
which allowed us to verify that the Fock representation of the Weyl algebra is cross for
any dynamics defined by a positive operator on the one-particle space. Second, we characterized
precisely the normal cross representations for inner one-parameter actions  on $W^*$-algebras.
Turning into a new direction, we added the requirement of a spectral condition to our study
of cross representations which brought in useful tools. For the one-parameter case, we found that
if the Borchers--Arveson minimal representation is cross, then so are all others. For
a one-parameter action on a topological group, we obtained a regularization property
for cross representations, for which the crossed product hosts are actually host algebras
for the semidirect product group, hence produce continuous group representations.
The action generated by the number operator on the
Heisenberg group in a Fock representation is one example.
Using smoothing operators we also did the Lie algebra version.
For the generalization of the spectral condition
to (possibly non-abelian) Lie groups we saw  two constructions of host algebras; one via smoothing
operators
and one via  Olshanski semigroups. The class of group representations controlled by
these host algebras are the semibounded ones.
For a covariant representation with semibounded unitary implementing representation $U$, we
found that it is cross if and only if the restriction of it to any of the one-parameter subsystems
with generator in $W_U$ is cross. This allows us to prove the cross property for
 the Fock representation of the Weyl algebra
with the symplectic action of either the Minkowski translation group, or
the conformal group based on a one-particle unitary representation where the joint spectrum of the
generators of translations are in the forward light cone.
 We obtained a non-abelian extension of the Borchers--Arveson Theorem to
the general Lie case when there is an invariant cyclic vector.
In the general case there are lifting obstructions that we describe in terms of
central extensions.

There are numerous open problems hence  possible future directions for this project.
\begin{itemize}
\item{}  Most desirable would be a full determination or characterization
of $C^*$-actions which admit cross representations, analogous  to Theorem~\ref{W-CrossRepsChar}
for inner one-parameter actions of von Neumann algebras.
For the case of a discontinuous action
of a locally compact group,  \cite[Cor.~8.4]{GrN14} is a first step.
\item{} For a $C^*$-action which admits
different crossed product hosts, analyze the relations between
these fully. In this direction we already have for a fixed host algebra
that two crossed product hosts for which there is containment of their sets of
associated representations, there is a factoring \cite[Thm.~5.8]{GrN14}.
So for the sets of cross representations associated with crossed product hosts the
partial order of containment is relevant.
But how are the crossed product hosts for maximal sets of cross representations associated to them
related?
What relation is there between crossed product hosts for different hosts?
\item{} Extend the tools associated with the spectral
condition to infinite dimensional Lie groups.
Natural examples would involve the Fock representation for an action of the
gauge group on the CAR algebra
obtained from a unitary representation of the group on the one-particle space.
Other natural examples come from projective unitary representations
of automorphism groups of the CAR algebra in quasi-free representations
(cf.~\cite{PS86, La94, JN19}).
What host algebra would we take for the group, and given that, is this Fock representation cross?
\item{} In our Example~\ref{FockPosCross} we left open whether the Fock representation is cross
 for the action of the Poincar\'e group alone on the Weyl algebra. This is important for physics to determine.
 \item{} We had a number of examples on the cross condition of the Fock representation
 for actions on the Weyl algebra (e.g. Examples~\ref{NoGap}  and ~\ref{FockPosCross}).
 Determine the cross property for the Fock representation of the analogous actions on the resolvent algebra.
 \item{} In Subsection~\ref{BATnonab} we discussed the extension of the Borchers--Arveson Theorem to
 non-abelian Lie groups and obtained one if there is an invariant cyclic vector. In the general case
 lifting obstructions are described in terms of central Lie group extensions and it would
be interesting to find more direct conditions for their vanishing. For finite dimensional
Lie groups, this could be based on the structure theory for Lie groups with semibounded
representations exposed in \cite{Ne99}.
 \item{}  The usual crossed products have been studied for a long time, and bring with them many
 useful tools and insights (cf.~\cite{Wil07}), e.g.\ inducing of covariant representations, direct integrals
 of covariant representations, structure theories for states and covariant representations and
 connections between ideal structures  (e.g. Gootman--Rosenberg Theorem). It is clear that
 by construction, many of these will adapt to crossed product hosts, and doing this explicitly is a
 fruitful future direction.
\end{itemize}

\section*{Acknowledgment}

Hendrik Grundling would like to thank Friedrich--Alexander University Erlangen--Nuremberg
who generously funded his visit 3-14/9/2018 to Erlangen
for joint work with Karl-Hermann Neeb on this manuscript.

\end{document}